\documentclass[a4paper,11pt]{article}
\usepackage[colorlinks=true,linkcolor=blue,citecolor=blue]{hyperref}
\usepackage{lineno,hyperref}
\usepackage{psfrag,graphicx,subfigure,epsfig}
\usepackage[dvipsnames]{xcolor}
\usepackage[sumlimits,namelimits]{amsmath}
\usepackage{amssymb,amsmath,amsthm,tocvsec2}

\usepackage[numbers]{natbib}
\numberwithin{equation}{section}

\usepackage{todonotes}

\def \btau{{\mbox{\boldmath$\tau$}}}
\newcommand{\bm}{\boldsymbol}
\newcommand{\Aum}{{\bm{u}}_m}
\newcommand{\Auc}{{\bm{u}}_c}
\newcommand{\Aumh}{{\bm{u}}_{mh}}
\newcommand{\Auch}{{\bm{u}}_{ch}}
\newcommand{\Austar}{{\bm{u}}_{c\star}}

\newcommand{\Au}{{\bm{u}}}
\newcommand{\Av}{{\bm{v}}}

\newcommand{\Avh}{\bm{v}_h}

\newcommand{\Anc}{{\bm{n}}_c}
\newcommand{\Anm}{{\bm{n}}_m}

\newcommand{\Aumbhn}{\bar{\bm{u}}_h^{n+1}}

\newcommand{\bx}{{\bm{x}}}
\newcommand{\BX}{{\bm{X}}}

\newcommand{\BV}{{\bm{V}}}

\newcommand{\hv}{{\hat {v}}}
\newcommand{\wh}{{\omega}_h}

\newcommand{\mK}{\mathbb K}
\newcommand{\Ag}{{\bm{g}}}
\newcommand{\Aj}{{\bm{j}}}

\title{A decoupled  numerical method  for two-phase flows of different densities and viscosities in superposed fluid and porous layers}

\author{
	{\sc Yali  Gao}\footnote{School of Mathematics and Statistics, Northwestern Polytechnical University, Xi'an Shaanxi, 710129, P.R.China
		Email: \emph{gaoylimath@nwpu.edu.cn}} ,\ \
	{\sc Daozhi Han}\footnote{Department of Mathematics and Statistics, Missouri University of Science and Technology, Rolla, MO 65401.
		Email: \emph{handaoz@mst.edu}} ,\ \
	{\sc Xiaoming He}\footnote{ Department of Mathematics and Statistics, Missouri University of Science and Technology, Rolla, MO 65401.  Email:
		\emph{hex@mst.edu}}\ \
	and \ {\sc Ulrich R\"{u}de}\footnote{Department of Computer Science 10, University Erlangen-Nuremberg, D-91058 Erlangen, Germany. Email: \emph{ulrich.ruede@fau.de}}
}

\begin{document}

\maketitle

\newtheorem{algorithm}{Algorithm}[section]
\newtheorem{theorem}{Theorem}[section]
\newtheorem{remark}{Remark}[section]
\newtheorem{lemma}{Lemma}[section]

\begin{abstract}
In this article we consider the numerical modeling and simulation via the phase field approach of two-phase flows of different densities and viscosities in superposed fluid and porous layers.
The model consists of the Cahn-Hilliard-Navier-Stokes equations in the free flow region and the Cahn-Hilliard-Darcy equations in porous media that are coupled by seven domain interface boundary conditions. We show that the coupled model satisfies an energy law. Based on the ideas of   pressure stabilization and artificial compressibility, we propose an unconditionally stable time stepping method that   decouples the computation of the phase field variable, the  velocity and pressure of free flow, the velocity and pressure of porous media,  hence significantly reduces the computational cost.  The energy stability of the scheme effected with the finite element spatial discretization is rigorously established. We verify numerically that our schemes are convergent and energy-law preserving.
Ample numerical experiments are performed to illustrate the features of two-phase flows in the coupled free flow and porous media setting.

\end{abstract}


\section{Introduction}
\label{1}
Multi-phase flow in superposed fluids and porous media has many applications in science and engineering. A prime example is the mixing of shallow groundwater and surface water in the hyporheic zone--a region of sediment and porous space beneath and alongside a stream bed. The hyporheic zone is a natural habitat for aquatic organisms  and plays a major role in maintaining the self-purification function of streams.  It is important to understand the  hydrodynamic and biogeochemical processes of multiphase nature in this zone, cf. \cite{MBCardenas_WRR_2015}. Other applications of multi-phase flow in superposed fluids and porous media include contaminant transport in karst aquifers \cite{JMatusick_PZanbergen_GRA_2007}, oil recovery in petroleum engineering  \cite{GlSw2004},
water management in PEM fuel cell technology \cite{KTuber_DPocza_CHebling_JPS_2003},  cardiovascular modeling and simulation \cite{LFormaggia_AQuarteroni_AVeneziani_2009} etc. Therefore it is of great importance to develop numerical models and efficient algorithms for simulating multi-phase flow in coupled free flow and porous media.

As a solid fundamental work, the modeling of single-phase flow (such as water) in superposed free flow and porous media is usually based on the Stokes-Darcy type systems, see \cite{MdAAlMahbub_XMHe_NJNasu_CQiu_HZheng_1,YCao_MGunzburger_XMHe_XWang_2,MDiscacciati_EMiglio_AQuarteroni_1, JHou_MQiu_XMHe_CGuo_MWei_BBai_1, WJLayton_FSchieweck_IYotov_1, BRiviere_IYotov_1, SKFStoter_PMuller_elat_CMAME_2017,STlupova_RCortez_1} and many others. There have been abundant numerics for investigating this type single-phase model, such as finite element methods \cite{YCao_MGunzburger_XHu_FHua_XWang_WZhao_1,JCamano_GNGatica_ROyarzua_RRuizBaier_PVenegas_1,AMarquez_SMeddahi_FJSayas_2}, discontinous Galerkin method \cite{PChidyagwai_BRiviere_1,VGirault_BRiviere_1,KLipnikov_DVassilev_IYotov_1,BRiviere_1},
domain decomposition methods \cite{YBoubendir_STlupova_2,YCao_MGunzburger_XMHe_XWang_1,MDiscacciati_LGerardo-Giorda_1,MDiscacciati_AQuarteroni_AValli_1,
MGunzburger_XMHe_BLi_SINN_2018,XMHe_JLi_YLin_JMing_1, CQiu_XMHe_JLi_YLin_1, DVassilev_CWang_IYotov_1},
multigrid methods \cite{TArbogast_MGomez_1,MCai_MMu_JXu_2,MMu_JXu_1},  and so on.

The study of multi-phase flow in this context is very challenging, and up to our knowledge no sharp interface model is available for two-phase flow in superposed fluids and porous media. In recent years, diffuse interface model has become popular in numerical modeling of multi-phase flow \cite{JLowengrub_LTruskinovsky_1998,AGLamorgese_DMolin_RMauri_MJM_2011}. In this approach the sharp interface of two immiscible fluids is replaced by a diffusive interface of finite thickness where different fluids mix due to chemical diffusion. The diffuse interface approach could describe topological transitions of interfaces and avoid the cumbersome procedure of interface tracking in numerical simulations, cf. \cite{HGLee_JLowengrub_JGoodman_PFI_2002, HGLee_JLowengrub_JGoodman_PFII_2002}. A hybrid of the sharp interface model  in  porous media and the  diffuse interface model in the free flow is  proposed in \cite{JChen_SSun_XWang_1}. In \cite{DHan_DSun_XWang_1} the authors systematically derive  a diffuse interface model, the Cahn-Hilliard-Stokes-Darcy model, for two-phase flow of matched/similar densities in the setting of coupled free flow and porous media. Well-posedness and numerical solvers for this model have been studied in \cite{DHan_XWang_HWu_1} and \cite{WChen_DHan_XWang_2017}, respectively. Generalization of the model to include inertia effect is done in \cite{YGao_XHe_LMei_XYang_2018}.  A diffuse interface model for two-phase flow of arbitrary densities and viscosities in superposed fluids and porous media remains open.

There are mainly two types of approach on developing diffuse interface models for two-phase flow of different densities in a single domain. The first approach defines a mass-averaged velocity that leads to a quasi-incompressible Cahn-Hilliard fluid model \cite{JLowengrub_LTruskinovsky_1998}. The mass-averaged velocity is non-solenoidal inside the diffusive interface, and the resulting model is a high order,  nonlinear, strongly coupled system that is difficult for numerical simulation \cite{GLLW2017}. The second approach adopts a volume-averaged velocity which is a solenoidal vector field everywhere, cf. \cite{FBoyer_3, HDing_PDMSpelt_CShu_1, HAbels_HGarcke_GGrun_2012}. Due to the divergence-free velocity, efficient legacy numerical solvers for incompressible fluid are applicable. However, the classical continuity equation for density is no longer valid if volume-averaged velocity is employed.

In this article we take an approximation approach for developing numerical models for superposed two-phase free flow and porous media, in the sense that we utilize the mass-averaged velocity but neglect the compressibility effect of the velocity field inside the thin diffusive interface. Such an approach has appeared in \cite{JShen_XYang_1} for numerical modeling of two-phase flow of variable densities in a single domain. On the domain interface boundary between free flow and porous media, we impose the Beavers-Joseph-Saffman-Jones interface boundary
condition \cite{GBeavers_DJoseph_1} and the Lions interface boundary condition which states that the free-flow stress in the normal direction including the total pressure (pressure plus dynamic pressure)  is balanced by the pressure in porous media. Under these conditions, we show that our model, the Cahn-Hilliard-Navier-Stokes-Darcy system, satisfies an energy law.

The design of accurate and long-time stable time-stepping method for the Cahn-Hilliard-Navier-Stokes-Darcy system is very challenging for a number of reasons. The first challenge is the stiffness inherent to diffuse interface models (large transition over thin  layers). There is a large body of literature on developing unconditionally stable time-marching algorithms for diffuse interface models. These methods include the convex-splitting strategy \cite{ElSt1993, JShen_CWang_XWang_SMWise_SINA_2012, BLWW2013, Grun2013, DHan_JSC_2016, RLi_YGao_JChen_LZhang_XMHe_ZChen_1,YLiu_WChen_CWang_SMWise_NM_2017,JYang_SMao_XMHe_XYang_YHe_1}, the stabilization method \cite{ShenYang2010, YYan_WChen_CWang_SMWise_CICP_2018}, the Invariant Energy Quadratization approach \cite{XYang_DHan_JCP_2017, QChen_XYang_JShen_JCP_2017, CXu_CChen_XYang_XMHe_1, XYang_LJu_CMAME_2017, JZhao_XYang_YGong_QWang_CMAME_2017,XYang_JZhao_XMHe_1}, and the Scalar Auxiliary Variable approach \cite{FLin_XMHe_XWen_1,JShen_JXu_JYang_JCP_2018, JShen_JXu_JYang_SIRV_2019}. The second issue is the coupling between the nonlinear Cahn-Hilliard equation and the fluid equations, and the coupling between fluid velocity and pressure. Operator-splitting  is typically utilized to decouple the computation, cf. \cite{DKay_RWelford_1,  DHan_JSC_2016, JShen_XYang_SINN_2015, WChen_DHan_XWang_2017}. The third challenge is the coupling between free flow and porous media via domain interface boundary conditions. Various domain decomposition approaches have been proposed to minimize the computational cost \cite{WJLayton_FSchieweck_IYotov_1, JChen_SSun_XWang_1, WChen_DHan_XWang_2017, YGao_XHe_LMei_XYang_2018}.

It is of great importance to develop decoupled numerical algorithms while maintaining the unconditional stability for solving the Cahn-Hilliard-Navier-Stokes-Darcy system (CHNSD).  A decoupled algorithm is proposed in \cite{WChen_DHan_XWang_2017} for solving the Cahn-Hilliard-Stokes-Darcy model in which the decoupling between the Cahn-Hilliard equation and fluid equations hinges upon the presence of time derivative in the Darcy equations. In our CHNSD model the governing equations for flow in porous media is the classical Darcy equations without the time derivative term. To accomplish the decoupling between phase field variable and Darcy velocity we resort to the technique of pressure stabilization from \cite{DHan_JSC_2016} originally designed for solving the Cahn-Hilliard-Darcy equations. Furthermore, it is desirable to separate the computation of velocity and pressure when solving the Navier-Stokes equations.    Due to the presence of the nonlinear Lions domain interface boundary condition, we adopt a special method of artificial compressibility \cite{VDeCaria_WLayton_MMcLaughlin_CMAME_2017,XHe_NJiang_CQiu_IJNME_2019} which avoids boundary conditions  in the update of the pressure.  We rigorously establish the unconditional long-time stability of the proposed algorithm and  verify numerically that the fully discrete schemes are convergent and energy-law preserving. Ample numerical experiments are performed to illustrate the distinctive features of two-phase flows in superposed fluids and porous media.

The rest of the article  is as follows. In Section \ref{sec(2)}, we propose the Cahn-Hilliard-Navier-Stokes-Darcy model for two-phase flows of arbitrary densities in superposed fluid and porous media, we show that the model satisfies an energy law, and we also develop an unconditionally stable coupled time-stepping method for solving the model. In Section \ref{sec(4)}, {we provide the fully discrete, decoupled numerical scheme  and establish its energy stability.}
 Numerical results {are} reported in Section \ref{sec(5)}.

\section{The Cahn-Hilliard-Navier-Stokes-Darcy model}\label{sec(2)}
In this section, we propose the Cahn-Hilliard-Navier-Stokes-Darcy model (CHNSD) for two-phase flows of {different densities and viscosities} in a fluid layer overlying porous media. We refer to \cite{DHan_DSun_XWang_1,DHan_XWang_HWu_1} for a phase field model for two-phase flows of matched density in the coupled setting where the linear flow regime (Stokes equations) is assumed in the free flow region.
 We provide the weak formulation of this model and then show that the model obeys a dissipative energy law on the PDE level. We also introduce an unconditionally stable coupled time-stepping method for solving the CHNSD system.

\subsection{The model}

\begin{figure}[h]
\centering
\setlength{\abovecaptionskip}{0pt}
\setlength{\belowcaptionskip}{5pt}
\includegraphics[height=2.5in]{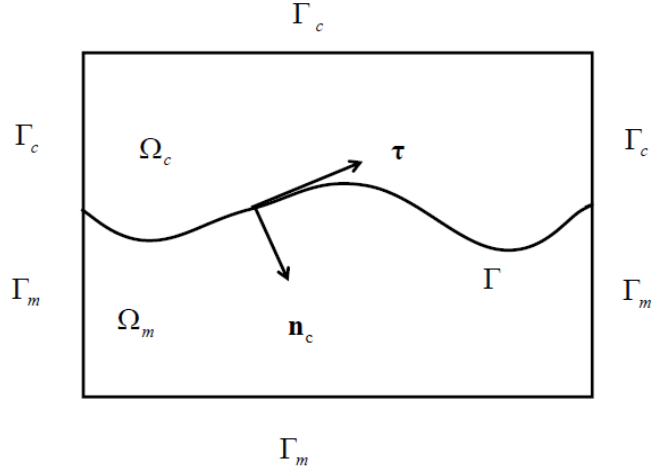}
\caption{A sketch of the porous median domain $\Omega_m$, fluid
domain $\Omega_c$, and the interface $\Gamma$.}
\label{Stokes_marcy_domain_illustration}
\end{figure}

We consider the coupled CHNSD system on  a bounded connected domain
$\Omega=\Omega_c \bigcup \Omega_m \subset {\mathbb
R}^{\mbox{\textbf{d}}}, \ (\textbf{d} =2, 3)$ consisting of a free-flow region $\Omega_c$ and a porous media region $\Omega_m$.
Let $\partial{\Omega}_c$ and $\partial{\Omega}_m$ denote the   Lipschitz continuous boundaries of $\Omega_c$ and $\Omega_m$ with the outward unit normal vectors $\Anc$ and $\Anm$ to the fluid and the porous media regions, respectively. The interface  between the two  parts  is denoted by $\Gamma$, i.e $\Gamma:=\partial{\Omega}_m\cap\partial{\Omega}_c$.
A typical two-dimensional geometry is illustrated in Figure \ref{Stokes_marcy_domain_illustration}.

Let $w_j~(j=c,m)$ denote the chemical potential  and $M_j~(j=c,m)$ denote a mobility constant related to the relaxation time scale. {Let $f(\phi)$ be a polynomial of $\phi$ such that $f(\phi)=F'(\phi)$, where $F(\phi)$ represents the Helmholtz free energy and is commonly taken to be a non-convex function of $\phi$ for two  immiscible flows. In this article, we consider the Ginzburg-Landau double-well potential $F(\phi)=\frac{1}{4\epsilon}(\phi^2-1)^2$ with the width of mixing layer $\epsilon$.} $\rho$ and $\eta$ are the density and
viscosity of the mixture, denoted by
\begin{eqnarray}\label{mix_define}
  \rho=\frac{\rho_1-\rho_2}{2}\phi+\frac{\rho_1+\rho_2}{2},\quad \nu=\frac{\nu_1-\nu_2}{2}\phi+\frac{\nu_1+\nu_2}{2}.
\end{eqnarray}
The gravity vector is $\Ag=g\Aj$ with the gravity constant $g$ and the  unit upward vector $\Aj$.
$\rho\Ag$ denotes the external gravitational forces.
Furthermore, $\gamma$ and $\epsilon$ denote the elastic relaxation time and the capillary width, respectively, of the thin interfacial region.
The order parameters (phase functions) are denoted by $\phi_j~(j=c,m)$ in $\Omega_j~(j=c,m)$ which assume distinct values $\pm 1$ respectively in the bulk phases away from the diffuse interface and varies smoothly inside it.

In the porous media region
$\Omega_m$, consider the porous media flow governed by
the following Cahn-Hilliard-Darcy (CHD) system:
\begin{eqnarray}
 \mK^{-1}\Aum+\nabla p_m+\phi_m\nabla w_m&=&\rho \Ag,\label{Darcy_law_time_BJ}\\
  \nabla \cdot \Aum&=&0,\label{Darcy_divergence_free_time_BJ}\\
  \frac{\partial \phi_m}{\partial t}+\nabla\cdot(\Aum\phi_m)-\nabla \cdot \left(M_m \nabla w_m\right)&=&0,\label{equation_for_marcy_time_BJ}\\
  w_m+\gamma\epsilon \triangle
\phi_m-\gamma f(\phi_m)&=& 0,\label{Darcy_chemical_potential_time_BJ}
\end{eqnarray}
where $\Aum$ is the fluid discharge rate in the porous
media, $\mathbb K$ is the hydraulic conductivity tensor, $p_m$ is the hydraulic head, and the term $w_m\nabla\phi_m$ is the induced extra stress from the free energy. Assuming external forces  to be zero, and inserting Darcy's law \eqref{Darcy_law_time_BJ} into the mass conservation equation \eqref{Darcy_divergence_free_time_BJ}, we will consider the second order formulation as follows:
\begin{eqnarray}
 -\nabla \cdot ({\mathbb K}\nabla p_m+{\mathbb K}\phi_m\nabla w_m)&=&0.\label{second_order_darcy}
\end{eqnarray}
After solving this equation, one can recover the Darcy velocity   via \eqref{Darcy_law_time_BJ}.

In the fluid region $\Omega_c$,  consider the two phase fluid flows governed by a coupled  Cahn-Hilliard-Navier-Stokes (CHNS) system with different densities and viscosities:
\begin{eqnarray}
  \rho\left(\frac{\partial \Auc}{\partial t}+(\Auc \cdot \nabla) \Auc\right)-\nabla\cdot\mathbb{T}(\Auc,p_c)+\phi_c\nabla w_c&=&\rho \Ag,\label{time_ctokes_equation_BJ}\\
  \nabla \cdot \Auc&=&0,\label{Stokes_divergence_free_time_BJ}\\
  \frac{\partial \phi_c}{\partial t}+\nabla\cdot(\Auc\phi_c)-\nabla \cdot \left(M_c \nabla w_c\right)&=&0,\label{Stokes_for_marcy_time_BJ}\\
  w_c+\gamma\epsilon \triangle
\phi_c-\gamma f(\phi_c)&=& 0,\label{Stokes_chemical_potential_time_BJ}
\end{eqnarray}
where $\Auc$ is the fluid
velocity, $p_c$ is the kinematic pressure, $\nu$ is the kinematic viscosity of the
fluid, $\mathbb{T}(\Auc,p_c)=2\nu \mathbb{D}(u_c)-p_c\mathbb{I}$ is the stress tensor, $\mathbb{D}(u_c)=(\nabla \Auc+\nabla^T\Auc)/2$ is the deformation tensor, and  ${\mathbb I}$ is the identity matrix.

We now introduce domain interface boundary conditions in order to couple the CHD system \eqref{Darcy_divergence_free_time_BJ}-\eqref{Darcy_chemical_potential_time_BJ} and CHNS system \eqref{time_ctokes_equation_BJ}-\eqref{Stokes_chemical_potential_time_BJ}.
The continuity of normal component of velocity  is assumed  across the interface
\begin{eqnarray}
   \Auc\cdot \Anc= -\Aum\cdot \Anm. \label{BJS-1_time}
\end{eqnarray}
The balance of  normal force over the interface is satisfied by~\cite{DHan_DSun_XWang_1,MCai_MMu_JXu_2,VGirault_BRiviere_1}
\begin{eqnarray}
    -\Anc\cdot (\mathbb{T}(\Auc,p_c)\cdot \Anc)+\frac{\rho}{2}|\Auc|^2= p_m. \label{BJS-3_time}
\end{eqnarray}
The   Beavers-Joseph-Saffman-Jones (BJS) interface condition \cite{GBeavers_DJoseph_1} holds as follows
\begin{eqnarray}
   - {\mathbf{\btau}}_j\cdot(\mathbb{T}(\Auc,p_c)\cdot \Anc)
   = \frac{\alpha\nu \sqrt{\mbox{d}}}{\sqrt{\mbox{trace$(\prod)$}}}\mathbf{\btau}_j \cdot \Auc, \label{BJS-2_time}
\end{eqnarray}
where $\mathbf{\btau}_j$ $(j=1,\cdots,d-1)$ are mutually
orthogonal unit tangential vectors through the interface $\Gamma$, and
$\prod$ is the permeability of the
porous media.

Moreover, we assume the continuity conditions for the phase field function,
the chemical potential, and their normal derivatives  on the interface $\Gamma$ \cite{DHan_DSun_XWang_1,DHan_XWang_HWu_1},
\begin{eqnarray}
\phi_c&=&\phi_m,\label{phase_interface_con_1}\\
 w_c&=&w_m,\label{phase_interface_con_2}\\
\nabla \phi_c\cdot \Anc &=&-\nabla \phi_m \cdot \Anm,\label{phase_interface_con_3}\\
M_c \nabla w_c \cdot \Anc&=&-M_m \nabla w_m \cdot \Anm .\label{phase_interface_con_4}
\end{eqnarray}

For the boundary conditions and initial conditions, we consider
\begin{eqnarray*}
\Aum\cdot\Anm|_{\Gamma_m} = 0,\label{DCH_bound_con_velocity}~~
\nabla \phi_m \cdot \Anm|_{\Gamma_m}= 0,\label{DCH_bound_con_phase}~~
M_m \nabla w_m \cdot \Anm|_{\Gamma_m} = 0, \label{DCH_bound_con_potential}
\end{eqnarray*}
on $\Gamma_m=\partial\Omega_m\backslash \Gamma$, and
\begin{eqnarray*}
\Auc|_{\Gamma_c} = 0,\label{NSCH_bound_con_velocity}~~
\nabla \phi_c \cdot \Anc|_{\Gamma_c} = 0,\label{NSCH_bound_con_phase}~~
M_c \nabla w_c \cdot \Anc|_{\Gamma_c} = 0,\label{NSCH_bound_con_potential}
\end{eqnarray*}
on $\Gamma_c=\partial\Omega_c\backslash \Gamma$. The initial conditions can be simply given as
\begin{eqnarray*}
   \phi_j(0,x,y) = \phi_j^0(x,y),\,\, j=c,m,~~
   \Auc(0,x,y) = \Auc^0(x,y).\label{initial condition for u_NSCH}
\end{eqnarray*}
For the ease of presentation, we assume external forces $\rho\Ag$ on the right side of equations  \eqref{Darcy_law_time_BJ} and \eqref{time_ctokes_equation_BJ} to be zero as these forces are given quantities which enter the system linearly.
Hence, they do not have a qualitative effect on estimates or results.
Without loss of generality,  we also assume that $\mathbb{K}$ is a bounded, symmetric and uniformly positive definite matrix.


\subsection{The weak formulation}
 In this subsection, we present the weak formulation of the CHNSD model system \eqref{Darcy_law_time_BJ}-\eqref{phase_interface_con_4}.
 Let $H^{m}\left(\Omega\right)$ be the classical Sobolev space with the norm $\|\cdot\|_m$,
where $m$ is a nonnegative integer. The norm $\|\cdot\|_\infty$ denotes the essential supremum. For the sake of simplicity,
we denote $L^2$ norm $\|\cdot\|_0$ by $\|\cdot\|$.
 Furthermore, we set $\BV=[H_0^1(\Omega)]^d=\{\Av \in [H^1(\Omega)]^d:\Av|_{\partial\Omega}=0\}$.

Given $v\in L^1(\Omega_j)~(j=c,m)$, we denote its mean value by $\hv=|\Omega_j|^{-1}\int_{\Omega_j}\,v(x)\,\mbox{d}x$.
Then we define the space
\begin{eqnarray}
\dot{L}^2(\Omega_j):=\{v\in L^2(\Omega_j):\int_{\Omega_j} v\,\mbox{d}\bx=0\}.
\end{eqnarray}
 Let  $\dot{H}^1(\Omega_j)=H^1(\Omega_j)\cap \dot{L}^2(\Omega_j)$
be a Hilbert space with inner product $(u,v)_{H^1}=\int_{\Omega_j}\,\nabla u\cdot \nabla v\,d\bx$ due to the classical Poincar\'{e} inequality for functions with zero mean. We  denote its dual space by $(\dot{H}^1(\Omega_j))'$. For the  coupled  CHNSD system, we introduce the following  spaces utilized throughout this paper
\begin{eqnarray}
   &&\BX_c = \{\Av \in [H^1(\Omega_c)]^d\,\,\,|\,\,\, \Av=0\ \mbox{on $ \Gamma_c$}\},\nonumber\\
   &&\BX_m = \{\Av \in [H^1(\Omega_m)]^d\,\,\,|\,\,\, \Av\cdot\Anm=0\ \mbox{on $ \Gamma_m$}\},\nonumber\\
 && \BX_{j,div}= \{\Av \in X_j\,\,\,|\,\,\, \nabla \cdot \Av=0\},\nonumber\\
 &&   Q_c=L^2(\Omega_c),\quad  Q_m= \dot{H}^1(\Omega_m),\nonumber\\
  && \quad Y_j= H^1(\Omega_j),\quad Y= H^1(\Omega),\quad j=c,m.\nonumber
\end{eqnarray}
Define $P_\tau$ to be the projection onto the tangent space on $\Gamma$, i.e.
$P_\tau\Au = \sum_{j=1}^{d-1}(\Au\cdot\btau_j)\btau_j$.
For the domain $\Omega_j~(j=c,m)$, $(\cdot,\cdot)$
denotes the $L^2$ inner product on the domain $\Omega_j$
decided by the subscript of integrated functions, and
$\langle\cdot,\cdot\rangle$ denotes the $L^2$ inner product on the
interface $\Gamma$. Then it is clear that
\begin{eqnarray*}
&&(u_m,v_m)=\int_{\Omega_m}u_mv_md\bx,\quad (u_c,v_c)=\int_{\Omega_c}u_cv_cd\bx,\quad (u,v)=\int_{\Omega_m}u_mv_md\bx+\int_{\Omega_c}u_cv_cd\bx,\\
&& \|u_m\|:=\left(\int_{\Omega_m}|u_m|^2d\bx\right)^{\frac{1}{2}},\quad \|u_c\|:=\left(\int_{\Omega_c}|u_c|^2d\bx\right)^{\frac{1}{2}},\quad \|u\|^2=\int_{\Omega_m}|u_m|^2d\bx+\int_{\Omega_c}|u_c|^2d\bx,
\end{eqnarray*}
where $u_m:=u|_{\Omega_m}$ and $u_c:=u|_{\Omega_c}$.
We also denote $H'$ the dual space of $H$ with the duality induced by the $L^2$ inner product.

Different from the equal density case \cite{MGao_XWang_1}, it is difficult to eliminate the nonlinear convective term of Navier-Stokes equation in the proof of the energy law. Hence a new variable $\sigma=\sqrt{\rho}$ is introduced to replace $\rho$ \cite{JLGuermond_LQuartapelle_2000}. Using the mass conservation
\begin{eqnarray}\label{mass_conservation}
\frac{\partial \rho}{\partial t}+\nabla\cdot(\rho\Auc)=0,
\end{eqnarray}
one can derive
\begin{eqnarray*}
\sigma\frac{\partial (\sigma\Auc)}{\partial t}=\rho\frac{\partial\Auc}{\partial t}+\frac{1}{2}\frac{\partial \rho}{\partial t}\Auc=\rho\frac{\partial\Auc}{\partial t}-\frac{1}{2}\nabla\cdot(\rho\Auc)\Auc.
\end{eqnarray*}
Therefore, \eqref{time_ctokes_equation_BJ} can be rewritten 
by replacing $\rho\frac{\partial\Auc}{\partial t}$ with
$\sigma\frac{\partial (\sigma\Auc)}{\partial t}+\frac{1}{2}\nabla\cdot(\rho\Auc)\Auc$  as
\begin{eqnarray}
&&\sigma\frac{\partial (\sigma\Auc)}{\partial t}+\rho\left(\Auc \cdot \nabla\right)\Auc-\nabla\cdot \mathbb{T}(\Auc,p_c)+\phi_c \nabla w_c
 +\frac{1}{2}\nabla\cdot(\rho\Auc)\Auc=0.\quad\label{time_ctokes_equation_BJ_1_an_time}
\end{eqnarray}
The application of this technique and the resulting   \eqref{time_ctokes_equation_BJ_1_an_time} in the context of multiphase flows first appears in \cite{ShenYang2010} by Shen and Yang. It is noted by Lowengrub and Truskinovsky in \cite{JLowengrub_LTruskinovsky_1998} that the mass-averaged velocity which maintains the continuity equation \eqref{mass_conservation} is quasi-incompressible, that is,  the mixture of two incompressible fluids is slightly compressible inside the diffusive interface. Hence the divergence-free condition in our model amounts to an approximation to the quasi-incompressibility of the mass-averaged velocity. This adoption is for the convenience of numerical modeling so that classical numerical methods  for incompressible fluid such as pressure-correction can be employed. The approximation can  be justified from the point-of-view of sharp interface limit, in the sense that our model will recover the sharp interface model as the interfacial width goes to zero in the case of single domains (the sharp interface model for two-phase flow in the coupled setting remains open). It is also a common practice to adopt the simplification of incompressibility when the Mach number is small. We point out that one could use the divergence-free (solenoidal) volume-averaged velocity as is proposed in \cite{HAbels_HGarcke_GGrun_2012} by Abels et al. In the formalism of volume-averaged velocity, the continuity equation \eqref{mass_conservation} is no longer valid, and there is an extra advection term from the chemical flux in the momentum equation. The numerical modeling utilizing the volume-averaged velocity is deferred to a future work.

By applying the interface conditions \eqref{BJS-1_time}-\eqref{phase_interface_con_4}, the weak formulation of the proposed Cahn-Hilliard-Navier-Stokes-Darcy is given as follows: find
\begin{eqnarray*}
(p_m,\Auc,p_c,\phi,w)\in (Q_m,\BX_c,Q_c,Y,Y)
\end{eqnarray*}
such that
\begin{eqnarray}
    &&({\mathbb K}\nabla p_m,\nabla q)+({\mathbb K}\phi_m \nabla w_m,\nabla q)-\langle\Auc\cdot\Anc,q\rangle = 0,~\forall~q \in Q_m,\quad\label{DCH_weak_formulation_BJS_1_time}\\
  &&(\sigma\frac{\partial (\sigma\Auc)}{\partial t},\Av)+\left(\rho\left(\Auc \cdot \nabla\right)\Auc,\Av\right)+(2\nu\mathbb{D}(\Auc),\mathbb{D}(\Av)) -(p_c,\nabla\cdot\Av)+(\phi_c \nabla w_c,\Av)+\frac{1}{2}(\nabla\cdot(\rho\Auc)\Auc,\Av)
  \nonumber\\
  &&\qquad\qquad\qquad
  +\langle p_m-\frac{\rho}{2}|\Auc|^2, \Av\cdot \Anc \rangle
       +\frac{\alpha \sqrt{\mbox{d}}}{\sqrt{\mbox{trace$(\prod)$}}}\langle \nu P_\tau\Auc,P_\tau\Av\rangle=0,~\forall~ \Av\in \BX_c,\quad\label{NSCH_weak_formulation_BJS_1_time}\\
  && (\nabla\cdot\Auc,q) =0,~\forall~\, q \in Q_c,\quad\label{NSCH_weak_formulation_BJS_2_time}\\
  && (\frac{\partial \phi}{\partial t},\psi)-(\Au\phi,\nabla\psi)+(M\nabla w,\nabla\psi)=0,~\forall~\psi\in Y,\quad\label{NSCH_weak_formulation_BJ_3_time}\\
   && (w,\omega)-\gamma\epsilon(\nabla\phi,\nabla \omega)-\gamma(f(\phi),\omega)=0,~\forall~\omega\in Y,
   \quad\label{NSCH_weak_formulation_BJ_4_time}
\end{eqnarray}
where $t\in [0,T]$, $T$ is a finite time, $\Aum \in L^\infty(0,T;[L^2(\Omega_m)]^d) \cap L^2(0,T;\BX_{m})$, $\Auc \in L^\infty(0,T;[L^2(\Omega_c)]^d) \cap L^2(0,T;\BX_{c,div})$, $\frac{\partial \Auc}{\partial t}\in L^2(0,T;\BX'_{c,div})$, $p_j \in L^2(0,T;Q_j)$, $\phi_j \in L^\infty(0,T;Y_j) \cap L^2(0,T;H^3(\Omega_j))$, $\frac{\partial \phi_j}{\partial t} \in  L^2(0,T;Y'_j)$, $w_j \in L^2(0,T;Y_j)$, and $j=\{c,m\}$,
~$\Aum$~is defined by
\begin{eqnarray}\label{newdefum}
\Aum=-{\mathbb K}\nabla p_m-{\mathbb K}\phi_m\nabla w_m
\end{eqnarray}
based on \eqref{Darcy_law_time_BJ}.
 We inherit the idea from~\cite{WChen_DHan_XWang_2017}~that we can solve a Cahn-Hilliard equations on the whole domain $\Omega$.
 This is an alternative to~\cite{YGao_XHe_LMei_XYang_2018},   where  two Cahn-Hilliard equations are solved on $\Omega_m$
and $\Omega_c$ separately.
 The other three interface conditions~\eqref{BJS-1_time}-\eqref{BJS-2_time}~are utilized in the traditional way for the single-phase Navier-Stokes-Darcy model in the literature~\cite{XMHe_JLi_YLin_JMing_1,VGirault_BRiviere_1,MLHadji_AAssala_FZNouri_1}.

\subsection{A dissipative energy law}
In order to show that the above weak formulation obeys a dissipative energy law, we first note that the total energy of the coupled system is given by
 \begin{eqnarray}
   E(t)=\frac{1}{2}\|\sigma\Auc\|^2
   +\gamma[\frac{\epsilon}{2}\|\nabla \phi\|^2+(F(\phi),1)].\quad\label{Energy_equation}
 \end{eqnarray}

\begin{theorem}\label{thm_energy}
Assume $(\Aum,\Auc,\phi)$ is a smooth solution to the initial boundary
value problem \eqref{Darcy_law_time_BJ}-\eqref{phase_interface_con_4}. Then $(\Aum,\Auc,\phi)$ satisfies the
basic energy law
 \begin{eqnarray}
  \frac{d}{dt}E(t)=-\mathcal{D}(t),\quad\label{Energy_equation_time}
 \end{eqnarray}
 where the energy dissipation $\mathcal{D}$ is given by
 \begin{eqnarray}\label{Energy_equation_dissipation}
   &&\mathcal{D}(t)=\|\sqrt{2\nu}\mathbb{D}(\Auc)\|^2+ M\|\nabla w\|^2+\|\sqrt{\mK^{-1}}\Aum\|^2+\frac{\alpha \sqrt{\mbox{d}}}{\sqrt{\mbox{trace$(\prod)$}}}\langle\nu P_\tau\Auc,P_\tau\Auc\rangle.\qquad
 \end{eqnarray}
\end{theorem}

\noindent{\bf Proof.}
First, choose the test functions $\Av=\Auc$  and $q=p_c$ in~\eqref{NSCH_weak_formulation_BJS_1_time}-\eqref{NSCH_weak_formulation_BJS_2_time}.
Adding the resultants together, we get
\begin{eqnarray}\label{weak_NSCH_Energy_1}
&&\frac{1}{2}\frac{d}{dt}\|\sigma\Auc\|^2+(\rho(\Auc\cdot\nabla) \Auc ,\Auc)+\frac{1}{2}(\nabla\cdot(\rho\Auc)\Auc,\Auc)+(\phi_c \nabla w_c,\Auc)
\nonumber\\
&&\hskip 1.9cm+\|\sqrt{2\nu}\mathbb{D}(\Auc)\|^2+\frac{\alpha \sqrt{\mbox{d}}}{\sqrt{\mbox{trace$(\prod)$}}}\langle\nu P_\tau\Auc,P_\tau\Auc\rangle+\langle\Auc\cdot\Anc,p_m-\frac{\rho}{2}|\Auc|^2\rangle=0.\quad\quad
 \end{eqnarray}
Using integration by parts, we can show that
\begin{eqnarray}\label{trilinear_property_NS_22}
 ((\Auc\cdot\nabla) \Av ,\Av)+\frac{1}{2}((\nabla\cdot \Auc)\Av,\Av)=\frac{1}{2}\langle\Auc\cdot\Anc,\Av\cdot\Av\rangle, ~\forall \Av\in V.
 \end{eqnarray}
Thanks to \eqref{trilinear_property_NS_22}, we have
\begin{eqnarray}
 &&((\rho\Auc\cdot\nabla) \Auc ,\Auc)+\frac{1}{2}(\nabla\cdot (\rho\Auc)\Auc,\Auc)=\frac{1}{2}\langle\rho\Auc\cdot\Anc,\Auc\cdot\Auc\rangle.\label{trilinear_property_NS_2}
\end{eqnarray}
Thus, applying \eqref{trilinear_property_NS_2}  in \eqref{weak_NSCH_Energy_1}, we obtain
\begin{eqnarray}\label{weak_NSCH_Energy_2}
&&\frac{1}{2}\frac{d}{dt}\|\sigma\Auc\|^2+(\phi_c \nabla w_c,\Auc)+\|\sqrt{2\nu}\mathbb{D}(\Auc)\|^2+\langle\Auc\cdot\Anc,p_m\rangle
+\frac{\alpha \sqrt{\mbox{d}}}{\sqrt{\mbox{trace$(\prod)$}}}\langle\nu P_\tau\Auc,P_\tau\Auc\rangle=0.\quad
 \end{eqnarray}

Second, by taking $q=p_m$ in \eqref{DCH_weak_formulation_BJS_1_time}, and applying \eqref{newdefum}, we obtain
\begin{eqnarray}\label{weak_DCH_Energy_1}
\begin{aligned}
-(\Aum,\nabla p_m)-\langle\Auc\cdot\Anc,p_m\rangle=0.
\end{aligned}
\end{eqnarray}
Taking the inner product of \eqref{newdefum} with $\Aum$, we have
  \begin{eqnarray}\label{weak_DCH_Energy_2}
  \|\sqrt{\mK^{-1}}\Aum\|^2=-(\nabla p_m,\Aum)-(\phi_m \nabla w_m,\Aum).
\end{eqnarray}
Adding \eqref{weak_DCH_Energy_1} and \eqref{weak_DCH_Energy_2}, we obtain
  \begin{eqnarray}\label{weak_DCH_Energy_3}
  \|\sqrt{\mK^{-1}}\Aum\|^2+(\phi_m \nabla w_m,\Aum)-\langle\Auc\cdot\Anc,p_m\rangle=0.\quad
  \end{eqnarray}
By taking $\psi=w$ and $\omega=-\frac{\partial \phi}{\partial t}$ in \eqref{NSCH_weak_formulation_BJ_3_time} and  \eqref{NSCH_weak_formulation_BJ_4_time}, respectively, and adding these two equations, we derive
 \begin{eqnarray}\label{weak_DCH_Energy_4}
\begin{aligned}
&\gamma[\frac{\epsilon}{2}\frac{d}{dt}\|\nabla \phi\|^2+\frac{d}{dt}(F(\phi),1)]+M\|\nabla w\|^2-(\Au\phi,\nabla w)=0.
\end{aligned}
 \end{eqnarray}

Summing  the above   resultants  \eqref{weak_NSCH_Energy_2}, \eqref{weak_DCH_Energy_3} and \eqref{weak_DCH_Energy_4} together,   we obtain~\eqref{Energy_equation_time}.
This completes the proof of Theorem \ref{thm_energy}.
\mbox{}\hfill$\Box$

\subsection{An unconditionally stable coupled time-stepping method}\label{sec(3)}
Unconditionally stable but coupled time-stepping methods can be readily constructed for solving the CHNSD system \eqref{DCH_weak_formulation_BJS_1_time}-\eqref{NSCH_weak_formulation_BJ_4_time}. Here we present such an example and discuss its   energy stability.

We shall follow the stabilization technique~\cite{ShenYang2010,ZhuCST1999,Xu06}   to handle the non-convex double-well potential $F(\phi)$. In order to ensure the stability of this  approach, we assume that the potential function $F(\phi)$ satisfying the following condition:
there exists a constant $L$ such that
\begin{eqnarray}\label{max_priciple}
\max_{\phi \in \mathbb{R}} | F''(\phi) |\leq L.
\end{eqnarray}
It is clear that the common Ginzburg-Landau double well potential $F(\phi)$  does not satisfy~\eqref{max_priciple}.
Following \cite{NCondette_CMMelcher_ESuli_2011,ShenYang2010},  one truncates $F(\phi)$, still denoted by $F(\phi)$, such that  \eqref{max_priciple} holds with  $L=\frac{2}{\epsilon}$ in \eqref{max_priciple}. We point out that  both the IEQ method and the SAV approach will lead to linear schemes with energy laws reformulated in terms of Lagrange multipliers.

Let $t_n, n=0,1 \cdots M$ be a uniform partition of $[0, T]$ with  $\Delta t=t_{n+1}-t_n=\frac{T}{M}$ being the time step size.  Then, we construct the following discrete time, and continuous space scheme in the weak form \eqref{DCH_weak_formulation_BJS_1_time}-\eqref{NSCH_weak_formulation_BJ_4_time}:
 Find
 \begin{eqnarray*}
 (p_m^{n+1},\Auc^{n+1},p_c^{n+1},\phi^{n+1},w^{n+1})\in (Q_m, \BX_c, Q_c, Y, Y)
 \end{eqnarray*}
such that for all $(q,\Av,q,\psi,\omega)\in (Q_m,\BX_c,Q_c,Y,Y)$
\begin{align}
 &({\mathbb K}\nabla p_m^{n+1},\nabla q)+ ({\mathbb K} \phi_m^n\nabla w_m^{n+1},\nabla q)-\langle\Auc^{n+1}\cdot\Anc,q\rangle = 0,~\forall~\, q \in Q_m,
 \label{DCH_semi_disretization_BJS_1_time}\\
  &(\sigma^{n+1}\frac{\sigma^{n+1}\Auc^{n+1}-\sigma^n\Auc^n}{\Delta t},\Av)
     +\left(\rho^n\left(\Auc^{n} \cdot \nabla\right)\Auc^{n+1},\Av\right)
     +(2\nu^{n}\mathbb{D}(\Auc^{n+1}),\mathbb{D}(\Av))
     \nonumber\\
    &\hskip 1.0cm-(p_c^{n+1},\nabla\cdot\Av)+(\phi_c^{n} \nabla w_c^{n+1} ,\Av)+\frac{1}{2}(\nabla\cdot(\rho^n\Auc^n)\Auc^{n+1},\Av)
    +\langle p_m^{n+1}, \Av\cdot \Anc \rangle \nonumber\\
    &\hskip 1.0cm-\frac{1}{2}\langle \rho^n \Auc^n\cdot\Auc^{n+1}, \Av\cdot \Anc \rangle
   +\frac{\alpha \sqrt{\mbox{d}}}{\sqrt{\mbox{trace$(\prod)$}}}\langle \nu^n P_\tau\Auc^{n+1},P_\tau\Av\rangle=0,
   ~\forall~ \Av\in \BX_c,\quad
  \label{NSCH_semi_disretization_BJS_1_time}\\
   & (\nabla\cdot\Auc^{n+1},q) =0,~\forall~\, q \in Q_c,\quad
   \label{NSCH_semi_disretization_BJS_2_time}\\
  & (\frac{\phi^{n+1}-\phi^n}{\Delta t},\psi)-(\Au^{n+1}\phi^n,\nabla\psi)+(M \nabla w^{n+1},\nabla\psi)=0,~\forall~\psi\in Y,\quad
  \label{NSCH_semi_disretization_BJ_3_time}\\
   & (w^{n+1},\omega)-\gamma\epsilon (\nabla\phi^{n+1},\nabla \omega)-\frac{\gamma}{\epsilon}(\phi^{n+1}-\phi^{n},\omega)
   -\gamma(f(\phi^n),\omega)=0,
~\forall~\omega\in Y,\quad
\label{NSCH_semi_disretization_BJ_4_time}
\end{align}
where
\begin{eqnarray}\label{DCH_velocity_semi}
\Aum^{n+1}=-{\mathbb K}\nabla p_m^{n+1}-{\mathbb K}\phi_m^n\nabla w_m^{n+1}.
\end{eqnarray}

We now proceed to prove the energy stability theorem as follows.
\begin{theorem}\label{thm_energy_semi}
The  scheme~\eqref{DCH_semi_disretization_BJS_1_time}-\eqref{NSCH_semi_disretization_BJ_4_time} is unconditionally energy stable, in the sense that its  approximation  $(\Auc^{n+1},\phi_m^{n+1},\phi_c^{n+1})$ satisfies the following discrete energy law:
 \begin{eqnarray}
 &&E^{n+1}-E^{n}\leq-\mathcal{D}^{n+1},\quad\label{Decouple_Energy_equation_time}
 \end{eqnarray}
 where the discrete energy
 $E$ is defined as
 \begin{eqnarray}
   E^{n}=\frac{1}{2}\|\sigma^n\Auch^n\|^2
+\gamma[\frac{\epsilon}{2}\|\nabla\phi^n\|^2+(F(\phi^n),1)],\quad\label{dis_Energy_equation}
 \end{eqnarray}
and  the energy dissipation $\mathcal{D}^{n+1}$ is given by
\begin{eqnarray}\label{Decouple_Energy_equation_dissipation}
&&\mathcal{D}^{n+1}=\frac{1}{2}\|\sigma^{n+1}\Auc^{n+1}-\sigma^n\Auc^n\|^2
+\Delta t\|\sqrt{2\nu^{n}}\mathbb{D}(\Auc^{n+1})\|^2+\Delta t\|\sqrt{\mK^{-1}}\Aum^{n+1}\|^2
+\frac{ \gamma\epsilon}{2}\|\nabla(\phi^{n+1}-\phi^n)\|^2
\nonumber\\
&&\hskip 1.0cm+\Delta t M\|\nabla w^{n+1}\|^2
+\Delta t\frac{\alpha \sqrt{\mbox{d}}}{\sqrt{\mbox{trace$(\prod)$}}}\langle \nu^n P_\tau\Auc^{n+1},P_\tau\Auc^{n+1}\rangle.
 \end{eqnarray}
\end{theorem}

\noindent{\bf Proof.} We first consider the Cahn-Hilliard part.
Taking $\psi=\Delta t w^{n+1}$ in~\eqref{NSCH_semi_disretization_BJ_3_time}, using
the identity
 \begin{eqnarray}\label{identity}
  2a(a-b)=a^2-b^2+(a-b)^2,
 \end{eqnarray}
 we get
 \begin{eqnarray}\label{Decouple_NSCH_Energy_8}
&&(\phi^{n+1}-\phi^n,w^{n+1})-\Delta t (\Au^{n+1}\phi^n,\nabla w^{n+1})+\Delta tM\|\nabla w^{n+1}\|^2=0.
 \end{eqnarray}

We take $\omega=-(\phi_c^{n+1}-\phi_c^n)$ in~\eqref{NSCH_semi_disretization_BJ_4_time},
 use~\eqref{identity} and the Taylor expansion
\begin{eqnarray}\label{Taylor_expansion}
F(\phi^{n+1})-F(\phi^n)=f(\phi^n)(\phi^{n+1}-\phi^n)+\frac{F''(\xi^n)}{2}(\phi^{n+1}-\phi^n)^2,
 \end{eqnarray}
to get
  \begin{eqnarray}\label{Decouple_NSCH_Energy_90}
&&-(w^{n+1},\phi^{n+1}-\phi^n)
+\frac{\gamma\epsilon}{2}[\|\nabla\phi^{n+1}\|^2-\|\nabla\phi^n\|^2+\|\nabla(\phi^{n+1}-\phi^n)\|^2]
+\frac{\gamma}{\epsilon}\|\phi^{n+1}-\phi^{n}\|^2\nonumber\\
&&\hskip 3.2cm+\gamma(F(\phi^{n+1})-F(\phi^n),1)
\leq \frac{\gamma}{2}|F''(\xi^n)|\|\phi^{n+1}-\phi^{n}\|^2 .
 \end{eqnarray}
Then, combining~\eqref{max_priciple}, we derive
 \begin{eqnarray}\label{Decouple_NSCH_Energy_9}
-(w^{n+1},\phi^{n+1}-\phi^n)
+\frac{\gamma\epsilon}{2}[\|\nabla\phi^{n+1}\|^2-\|\nabla\phi^n\|^2+\|\nabla(\phi^{n+1}-\phi^n)\|^2]
+\gamma(F(\phi^{n+1})-F(\phi^n),1)\leq0.\quad
 \end{eqnarray}
 Adding~\eqref{Decouple_NSCH_Energy_8} and \eqref{Decouple_NSCH_Energy_9} together,
we get
 \begin{eqnarray}\label{Decouple_NSCH_Energy S1_9}
&&\frac{\gamma\epsilon}{2}[\|\nabla\phi^{n+1}\|^2-\|\nabla\phi^n\|^2]
+\gamma(F(\phi^{n+1})-F(\phi^n),1)+\Delta tM\|\nabla w^{n+1}\|^2\nonumber\\
&&\hskip 2.5cm
+\frac{\gamma\epsilon}{2}\|\nabla(\phi^{n+1}-\phi^n)\|^2-\Delta t (\Au^{n+1}\phi^n,\nabla w^{n+1})\leq0.
\end{eqnarray}

Then, we consider conduit part.
Thanks to \eqref{trilinear_property_NS_22}, we have
 \begin{eqnarray}
 &&((\rho^n\Auc^n\cdot\nabla) \Auc^{n+1} ,\Auc^{n+1})+\frac{1}{2}(\nabla\cdot (\rho^n\Auc^n)\Auc^{n+1},\Auc^{n+1})=\frac{1}{2}\langle\rho^n\Auc^n\cdot\Auc^{n+1},\Auc^{n+1}\cdot\Anc\rangle.\label{trilinear_property_NS_32}
\end{eqnarray}
By taking the test function $\Av=\Delta t\Auc^{n+1}$ in~\eqref{NSCH_semi_disretization_BJS_1_time}, $q=\Delta tp_c^{n+1}$ in \eqref{NSCH_semi_disretization_BJS_2_time}, summing the resultants,
applying~\eqref{trilinear_property_NS_32}  and \eqref{identity},
we obtain
\begin{eqnarray}\label{Decouple_NSCH_Energy S1_1}
&&\frac{1}{2}[\|\sigma^{n+1}\Auc^{n+1}\|^2-\|\sigma^{n}\Auc^n\|^2+\|\sigma^{n+1}\Auc^{n+1}-\sigma^{n}\Auc^n\|^2]
+\Delta t\|\sqrt{2\nu^{n}}\mathbb{D}(\Auc^{n+1})\|^2+\Delta t( \phi_c^{n}\nabla w_c^{n+1},\Auc^{n+1})\nonumber\\
&&\hskip 2.5cm+\Delta t \langle\Auc^{n+1}\cdot\Anc,p_m^{n+1}\rangle+\Delta t\frac{\alpha \sqrt{\mbox{d}}}{\sqrt{\mbox{trace$(\prod)$}}}\langle \nu^{n} P_\tau\Auc^{n+1},P_\tau\Auc^{n+1}\rangle=0.
\end{eqnarray}

Next, we consider the matrix part. Choosing  $q=\Delta tp_h^{n+1}$ in \eqref{DCH_semi_disretization_BJS_1_time} and taking the inner product of~\eqref{DCH_velocity_semi} with $v=\Delta t \Aum^{n+1}$,
 then adding the resultants together, we derive
 \begin{eqnarray}\label{Decouple_DCH_Energy_S1_3}
 &&\Delta t\|\sqrt{\mK^{-1}}\Aum^{n+1}\|^2+\Delta t(\phi_m^n\nabla w_m^{n+1},\Aum^{n+1})-\Delta t\langle\Auc^{n+1}\cdot\Anc,p_m^{n+1}\rangle=0 .
\end{eqnarray}

Summing~\eqref{Decouple_NSCH_Energy S1_9}, \eqref{Decouple_NSCH_Energy S1_1} and \eqref{Decouple_DCH_Energy_S1_3} together, we have
\begin{eqnarray}\label{Decouple_DCH_Energy_9}
\begin{aligned}
E^{n+1}-E^n
&\leq-\frac{1}{2}\|\sigma^{n+1}\Auc^{n+1}-\sigma^{n}\Auc^n\|^2
-\Delta t\|\sqrt{2\nu^{n}}\mathbb{D}(\Auc^{n+1})\|^2-\Delta t \|\sqrt{\mK^{-1}}\Aum^{n+1}\|^2
\\
&-\frac{ \gamma\epsilon}{2}\|\nabla(\phi^{n+1}-\phi^n)\|^2
-\Delta t M\|\nabla w^{n+1}\|^2
-\Delta t\frac{\alpha \sqrt{\mbox{d}}}{\sqrt{\mbox{trace$(\prod)$}}}\langle \nu^n P_\tau\Auc^{n+1},P_\tau\Auc^{n+1}\rangle,
\end{aligned}
\end{eqnarray}
namely, we obtain \eqref{Decouple_Energy_equation_time}.  Therefore, the conclusion of Theorem \ref{thm_energy_semi} follows. \mbox{}\hfill$\Box$

\section{An unconditionally stable decoupled numerical scheme}\label{sec(4)}
In this section, we present  an unconditionally stable decoupled numerical scheme  for solving the CHNSD model.
Finite elements are used for the spatial discretization.
 Let $\Im_h$ be a quasi-uniform triangulation  of domain $\Omega$ under mesh size $h$.
We introduce the  finite element spaces $Y_h\subset Y$,  $Y_{jh}\subset Y_j$,
$\BX_{ch}\subset \BX_c$ and $Q_{jh}\subset Q_j$ with $j=c,m$.
Here we assume $\BX_{ch}\subset \BX_c$ and $Q_{ch}\subset Q_c$ satisfy an
inf-sup condition for the divergence operator in the following form:
There exists a constant $C>0$ independent of $h$ such that the  LBB condition
\begin{eqnarray*}
\inf_{0\neq q_h}\sup_{0\neq \Avh}
\frac{(\nabla\cdot\Avh,q_h)}{\|\Avh\|_1}> C\|q_h\|, \; \forall~ q_h\in Q_{ch}, \Avh\in \BX_{ch}
\label{inf_sup_condition}
\end{eqnarray*}
holds.

We first recall the following lemma for the estimate of the interface term from \cite{WChen_DHan_XWang_2017,MMoraiti_1}:
\begin{lemma}\label{lemma31}
There exists a constant $C$ such that, for $\Av\in\BX_{c}$, $q_{mh}\in Q_{mh}$
\begin{eqnarray}\label{bound_interface1}
|\langle\Av\cdot\Anc,q_{mh}\rangle|\leq C \|\Av\|_{\BX_{div}}\|\nabla q_{mh}\|,
\end{eqnarray}
where $\|\Av\|_{\BX_{div}}^2=\|\Av\|^2+\|\nabla\cdot\Av\|^2$.
\end{lemma}

In order to decouple the velocity and pressure in the Navier-Stokes equations, we follow the idea of artificial compressibility method \cite{ AJChorin_2,VDeCaria_WLayton_MMcLaughlin_CMAME_2017,XHe_NJiang_CQiu_IJNME_2019,RTemam_BSMF_1968, NNYanenko_1971} and
 replace the divergence-free condition by
\begin{eqnarray*}
 \nabla\cdot\Av-\delta p_t=0,
\end{eqnarray*}
where $\delta$ is an  artificial compression parameter such that the pressure can be solved explicitly.  We propose the following  decoupled, unconditionally stable,  linear scheme:

 {\bf Step 1.} Find $(\phi_h^{n+1},w_h^{n+1})\in Y_{h}\times Y_{h}$, such that
\begin{eqnarray}
  && (\frac{\phi_h^{n+1}-\phi_h^n}{\Delta t},\psi_h)-(\bar{u}_h^{n+1} \phi_h^n,\nabla\psi_h)+(M\nabla w_h^{n+1},\nabla \psi_h)=0,~\forall~\psi_h\in Y_h,
  \quad\label{2Decouple_DCH_full_disretization_BJS_2_time}\\
  && (w_h^{n+1},\wh)-\gamma\epsilon(\nabla\phi_h^{n+1},\nabla \wh)-\frac{\gamma}{\epsilon}(\phi_h^{n+1}-\phi_h^{n},\wh)
  -\gamma(f(\phi_h^n),\wh)=0,~\forall~\wh\in Y_h,\quad\label{2Decouple_DCH_full_disretization_BJS_3_time}
  \end{eqnarray}
  where
  \begin{eqnarray*}
\bar{u}_h^{n+1}= \left\{
  	\begin{array}{ll}
      \Austar^n,\quad \bx\in\Omega_c,\\
      \Aumh^{n+1},\quad \bx\in\Omega_m,
    \end{array}
   \right.
\end{eqnarray*}
and  $\Austar^n$ and $\Aumh^{n+1}$ are defined as following
 \begin{eqnarray}
\Austar^n&=&\Auch^n-\frac{1}{\rho^n}\Delta t \phi_{ch}^n\nabla w_{ch}^{n+1},\label{CS2_flux_identity}\\
\Aumh^{n+1}&=&-{\mathbb K}\nabla p_{mh}^n-{\mathbb K}\phi_{mh}^n \nabla w_{mh}^{n+1}.\label{DCH_velocity_time}
 \end{eqnarray}
  {\bf Step 2.} Find $p_{mh}^{n+1}\in  Q_{mh}$, such that
\begin{eqnarray}
 &&({\mathbb K}\nabla p_{mh}^{n+1},\nabla q_h)+({\mathbb K}\phi_{mh}^n\nabla w_{mh}^{n+1},\nabla q_h)+\beta\Delta t(\nabla p_{mh}^{n+1},\nabla q_h)-\langle\Auch^{n}\cdot\Anc,q_h\rangle = 0,~\forall~\, q_h \in Q_{mh}.\label{2Decouple_DCH_full_disretization_BJS_1_time}
  \end{eqnarray}
    {\bf Step 3.}  Find $\Auch^{n+1}\in  \BX_{ch}$,
     such that
  \begin{eqnarray}\label{2Decouple_NSCH_full_disretization_BJS_1_time}
     &&(\frac{\bar{\rho}^{n+1}\Auch^{n+1}-\rho^{n}\Auch^n}{\Delta t},\Avh)
     +\left(\rho^{n}\left(\Auch^{n} \cdot \nabla\right)\Auch^{n+1},\Avh\right)
     +(2\nu^{n}\mathbb{D}(\Auch^{n+1}),\mathbb{D}(\Avh))
    +( \phi_{ch}^{n}\nabla w_{ch}^{n+1},\Avh) \nonumber\\
    &&\hskip 2.0cm
    -(2p_{ch}^{n}-p_{ch}^{n-1},\nabla\cdot\Avh)+\frac{1}{2}(\nabla\cdot(\rho^{n}\Auch^n)\Auch^{n+1},\Avh)
+\frac{\xi}{\Delta t} (\nabla\cdot (\Auch^{n+1}-\Auch^{n}),\nabla\cdot\Avh)\quad \\
    &&\hskip 2.0cm+\langle p_{mh}^{n+1}, \Avh\cdot \Anc \rangle-\frac{1}{2}\langle \rho^n\Auch^n\cdot\Auch^{n+1}, \Avh\cdot \Anc \rangle+\frac{\alpha \sqrt{\mbox{d}}}{\sqrt{\mbox{trace$(\prod)$}}}\langle \nu^n P_\tau\Auch^{n+1},P_\tau\Avh\rangle=0,~\forall~\Avh\in \BX_{ch},\nonumber
    \end{eqnarray}
    with $\bar{\rho}^{n+1}=\dfrac{{\rho}^{n+1}+{\rho}^{n}}{2}$.\\
     {\bf Step 4:}  Find $p_{ch}^{n+1}\in   Q_{ch}$,
     such that
  \begin{eqnarray}
    &&(p_{ch}^{n+1}-p_{ch}^{n},q_h)=-\frac{\zeta}{\Delta t}(\nabla\cdot\Auch^{n+1},q_h),~\forall~\, q_h \in Q_{ch},\quad\label{2Decouple_NSCH_full_disretization_BJS_5_time}
   \end{eqnarray}
with  $\zeta=\frac{1}{4}\min\{\rho_1,\rho_2\}$.

\begin{remark}\label{re3.7}
The term $\beta \Delta t (\nabla p_{mh}^{n+1},\nabla q_h)$  in \eqref{2Decouple_DCH_full_disretization_BJS_1_time} is a  stabilization  term in order to deduce the unconditional stability for the linearized  scheme \eqref{2Decouple_DCH_full_disretization_BJS_1_time}.
 The parameter $\alpha$ depends only on the geometry of $\Omega$.
\end{remark}

\begin{remark}\label{re3.8}
The term $\frac{\xi}{\Delta t} \nabla (\nabla\cdot (\Auch^{n+1}-\Auch^{n}))$ in \eqref{2Decouple_NSCH_full_disretization_BJS_1_time} is a   term to ensure the energy stability of continuity equation~\cite{VDeCaria_WLayton_MMcLaughlin_CMAME_2017,JAFiordilino_WLayton_YRong_CMAME_2018}. Thus, one can derive the stability of numerical method under some approximate constant for $\xi$.
\end{remark}

\begin{remark}\label{re3.6}
 The scheme \eqref{2Decouple_DCH_full_disretization_BJS_2_time}-\eqref{2Decouple_NSCH_full_disretization_BJS_5_time}
is a decoupled, linear scheme. Indeed, \eqref{2Decouple_DCH_full_disretization_BJS_2_time}-\eqref{2Decouple_DCH_full_disretization_BJS_3_time}, \eqref{2Decouple_DCH_full_disretization_BJS_1_time}, \eqref{2Decouple_NSCH_full_disretization_BJS_1_time} and \eqref{2Decouple_NSCH_full_disretization_BJS_5_time}  are decoupled linear elliptic equations for $\phi_h^{n+1}$, $w_h^{n+1}$, $p_{mh}^{n+1}$, $\Auch^{n+1}$ and $p_{ch}^{n+1}$. Therefore, at each time step, one only needs to solve a sequence of linear equations which can be solved very efficiently.
\end{remark}

We now prove the energy stability theorem as follows.
\begin{theorem}\label{S2_thm_energy_Decouple}
Let $(\Auch^{n+1},p_{mh}^{n+1},p_{ch}^{n+1},\phi_h^{n+1})$ be a smooth solution to the initial boundary value problem \eqref{2Decouple_DCH_full_disretization_BJS_2_time}-\eqref{2Decouple_NSCH_full_disretization_BJS_5_time}. Then the approximation  $(\Auch^{n+1},p_{mh}^{n+1},p_{ch}^{n+1},\phi_h^{n+1})$  satisfies the following modified discrete energy law:
 \begin{eqnarray}
  &&\mathrm{\mathcal{E}}^{n+1}-\mathrm{\mathcal{E}}^{n}\leq-\mathcal{D}^{n+1},\quad\label{S2_Decouple_Energy_equation_time}
 \end{eqnarray}
 where the modified discrete energy $\mathrm{\mathcal{E}}^{n+1}$ is defined as
 \begin{eqnarray}\label{S2_dis_Energy_equation}
&&\mathrm{\mathcal{E}}^{n}=E^n+\frac{\xi}{2} \|\nabla\cdot \Auch^{n}\|^2+\frac{\Delta t^2}{2\zeta} \| p_{ch}^{n}\|^2+\frac{1}{2}\Delta t \|\sqrt{\mK}\nabla p_{mh}^n\|^2,
 \end{eqnarray}
 with
 \begin{eqnarray*}
   E^{n}=\int_{\Omega_c}\frac{1}{2}|\sigma^n\Auch^{n}|^2d\bx+\gamma\int_{\Omega}[\frac{\epsilon}{2}|\nabla \phi_h^{n}|^2+F(\phi_h^{n})]d\bx,
 \end{eqnarray*}
and  the energy dissipation $\mathcal{D}^{n+1}$ is given by
\begin{eqnarray}\label{S2_Decouple_Energy_equation_dissipation}
&&\mathcal{D}^{n+1}=\Delta t\|\sqrt{2\nu^{n}}\mathbb{D}(\Auch^{n+1})\|^2+\Delta tM\|\nabla w_h^{n+1}\|^2
+\frac{\gamma\epsilon}{2}\|\nabla\phi_h^{n+1}-\nabla\phi_h^n\|^2+\frac{\Delta t^2}{2\zeta}\|p_{ch}^{n}-p_{ch}^{n-1}\|^2
\nonumber\\
&&\hskip 1.3cm
+\frac{1}{4}\Delta t\|\sqrt{\mK}\nabla (p_{mh}^{n+1}-p_{mh}^{n})\|^2+\Delta t\frac{\alpha \sqrt{\mbox{d}}}{\sqrt{\mbox{trace$(\prod)$}}}\langle \nu^{n} P_\tau\Auch^{n+1},P_\tau\Auch^{n+1}\rangle.
 \end{eqnarray}

\end{theorem}

\noindent{\bf Proof.} We firstly consider the full discretization~\eqref{2Decouple_DCH_full_disretization_BJS_2_time} and~\eqref{2Decouple_DCH_full_disretization_BJS_3_time} for Cahn-Hilliard euqation on whole domain $\Omega$.
Taking $\psi_h=\Delta t w_h^{n+1}$ in~\eqref{2Decouple_DCH_full_disretization_BJS_2_time}, we get
 \begin{eqnarray}\label{2Decouple_NSCH_Energy_8}
&&(\phi_h^{n+1}-\phi_h^n,w_h^{n+1})-\Delta t (\Aumbhn\phi_h^n,\nabla w_h^{n+1})+\Delta tM\|\nabla w_h^{n+1}\|^2=0.
 \end{eqnarray}
We take $\wh=-(\phi_h^{n+1}-\phi_h^n)$ in~\eqref{2Decouple_DCH_full_disretization_BJS_3_time}, the equality \eqref{identity}
  and the Taylor expansion \eqref{Taylor_expansion}
to get
  \begin{eqnarray}\label{2Decouple_NSCH_Energy_90}
&&-(w_h^{n+1},\phi_h^{n+1}-\phi_h^n)
+\frac{\gamma\epsilon}{2}[\|\nabla\phi_h^{n+1}\|^2-\|\nabla\phi_h^n\|^2]
+\gamma(F(\phi_h^{n+1})-F(\phi_h^n),1)+\frac{\gamma}{\epsilon}\|\phi_h^{n+1}-\phi_h^{n}\|^2\nonumber\\
&&\hskip 3.3cm+\frac{\gamma\epsilon}{2}\|\nabla\phi_{ch}^{n+1}-\nabla\phi_{ch}^n\|^2
\leq \frac{\gamma}{2}|F''(\xi^n)|\|\phi_h^{n+1}-\phi_h^{n}\|^2.\quad
 \end{eqnarray}
Then, combining~\eqref{max_priciple} and \eqref{2Decouple_NSCH_Energy_90}, we derive
 \begin{eqnarray}\label{2Decouple_NSCH_Energy_9}
&&-(w_h^{n+1},\phi_h^{n+1}-\phi_h^n)
+\frac{\gamma\epsilon}{2}[\|\nabla\phi_h^{n+1}\|^2-\|\nabla\phi_h^n\|^2]
+\gamma(F(\phi_h^{n+1})-F(\phi_h^n),1)
\leq-\frac{\gamma\epsilon}{2}\|\nabla\phi_h^{n+1}-\nabla\phi_h^n\|^2.\quad
 \end{eqnarray}
Adding \eqref{2Decouple_NSCH_Energy_8} and \eqref{2Decouple_NSCH_Energy_9}, we obtain
 \begin{eqnarray}\label{2Decouple_CH_Energy_3}
&&\frac{\gamma\epsilon}{2}[\|\nabla\phi_h^{n+1}\|^2-\|\nabla\phi_h^n\|^2]
+\gamma(F(\phi_h^{n+1})-F(\phi_h^n),1)-\Delta t (\Aumbhn\phi_h^n,\nabla w_h^{n+1})\nonumber\\
&&\hskip 2.5cm\leq-\frac{\gamma\epsilon}{2}\|\nabla\phi_h^{n+1}-\nabla\phi_h^n\|^2-\Delta tM\|\nabla w_h^{n+1}\|^2.\quad
 \end{eqnarray}

Next, we discuss the conduit part. Taking the test function $\Avh=\Delta t\Auch^{n+1}$ in~\eqref{2Decouple_NSCH_full_disretization_BJS_1_time},
combining~\eqref{trilinear_property_NS_32}, \eqref{CS2_flux_identity}, and
the identity \eqref{identity}, we obtain
 \begin{eqnarray}\label{Decouple_NSCH_Energy S2_1}
&&\frac{1}{2}[\|\sigma^{n+1}\Auch^{n+1}\|^2-\|\sigma^{n}\Austar^n\|^2
+\|\sigma^{n}\left(\Auch^{n+1}-\Austar^n\right)\|^2]+\Delta t\|\sqrt{2\nu^{n}}\mathbb{D}(\Auch^{n+1})\|^2
\nonumber\\
&&\hskip 2.1cm+\frac{\xi}{2} [\|\nabla\cdot \Auch^{n+1}\|^2-\|\nabla\cdot \Auch^{n}\|^2+\|\nabla\cdot(\Auch^{n+1}-\Auch^{n})\|^2]\nonumber\\
&&\hskip 2.1cm+\Delta t(p_{ch}^{n+1}-2p_{ch}^{n}+p_{ch}^{n-1},\nabla\cdot\Auch^{n+1})
-\Delta t(p_{ch}^{n+1},\nabla\cdot\Auch^{n+1})
\nonumber\\
&&
\hskip 2.1cm+\Delta t \langle\Auch^{n+1}\cdot\Anc,p_{mh}^{n+1}\rangle+\Delta t\frac{\alpha \sqrt{\mbox{d}}}{\sqrt{\mbox{trace$(\prod)$}}}\langle \nu^{n} P_\tau\Auch^{n+1},P_\tau\Auch^{n+1}\rangle=0.
\end{eqnarray}
Taking $q_h=\dfrac{\Delta t^2}{\zeta}(p_{ch}^{n+1}-2p_{ch}^{n}+p_{ch}^{n-1})$ in
\eqref{2Decouple_NSCH_full_disretization_BJS_5_time}, and using~\eqref{identity},
we have
\begin{eqnarray}\label{Decouple_NSCH_Energy S2_2}
&&\frac{\Delta t^2}{2\zeta}[\|p_{ch}^{n+1}-p_{ch}^{n}\|^2-\|p_{ch}^{n}-p_{ch}^{n-1}\|^2+\|p_{ch}^{n+1}-2p_{ch}^n+p_{ch}^{n-1}\|^2]
    =\Delta t(\nabla\cdot\Auch^{n+1},p_{ch}^{n+1}-2p_{ch}^n+p_{ch}^{n-1}).\quad
\end{eqnarray}

Taking $q_h=-\dfrac{\Delta t^2}{\zeta}p_{ch}^{n+1}$ in
\eqref{2Decouple_NSCH_full_disretization_BJS_5_time}, and using~\eqref{identity},
we obtain
\begin{eqnarray}\label{Decouple_NSCH_Energy S2_3}
&&\frac{\Delta t^2}{2\zeta}[\| p_{ch}^{n+1}\|^2-\|p_{ch}^{n}\|^2+\|p_{ch}^{n+1}-p_{ch}^{n}\|^2]
    =-\Delta t(\nabla\cdot\Auch^{n+1},p_{ch}^{n+1}).\qquad
\end{eqnarray}
Adding \eqref{Decouple_NSCH_Energy S2_2} and \eqref{Decouple_NSCH_Energy S2_3} to get
\begin{eqnarray}\label{Decouple_NSCH_Energy S2_4}
&&\frac{\Delta t^2}{2\zeta}[\| p_{ch}^{n+1}\|^2-\| p_{ch}^{n}\|^2+\|p_{ch}^{n}-p_{ch}^{n-1}\|^2
-\|p_{ch}^{n+1}-2p_{ch}^{n}+p_{ch}^{n-1}\|^2]\nonumber\\
&&\hskip 2.4cm=\Delta t(\nabla\cdot\Auch^{n+1},p_{ch}^{n+1}-2p_{ch}^{n}+p_{ch}^{n-1})-\Delta t(\nabla\cdot\Auch^{n+1},p_{ch}^{n+1}).
\end{eqnarray}
Now, we estimate the term $\|p_{ch}^{n+1}-2p_{ch}^{n}+p_{ch}^{n-1}\|^2$ on the right hand side of \eqref{Decouple_NSCH_Energy S2_4}. Taking the difference of \eqref{2Decouple_NSCH_full_disretization_BJS_5_time} at step $t^{n+1}$ and step $t^n$ to derive,
\begin{eqnarray}\label{Decouple_NSCH_Energy S2_5}
&&p_{ch}^{n+1}-2p_{ch}^{n}+p_{ch}^{n-1} =-\frac{\zeta}{\Delta t}\nabla\cdot(\Auch^{n+1}-\Auch^n),\qquad
\end{eqnarray}
which implies
\begin{eqnarray}\label{Decouple_NSCH_Energy S2_7}
    &&\frac{\Delta t^2}{2\zeta}\|p_{ch}^{n+1}-2p_{ch}^{n}+p_{ch}^{n-1}\|^2\leq\frac{\zeta}{2}\|\nabla\cdot(\Auch^{n+1}-\Auch^n)\|^2.
\end{eqnarray}

Adding \eqref{Decouple_NSCH_Energy S2_1}, \eqref{Decouple_NSCH_Energy S2_4} and \eqref{Decouple_NSCH_Energy S2_7}, we obtain
\begin{eqnarray}\label{Decouple_NSCH_Energy S2_8}
&&\frac{1}{2}[\|\sigma^{n+1}\Auch^{n+1}\|^2-\|\sigma^{n}\Austar^n\|^2
+\|\sigma^{n}\left(\Auch^{n+1}-\Austar^n\right)\|^2]
+\Delta t\|\sqrt{2\nu^{n}}\mathbb{D}(\Auch^{n+1})\|^2+\frac{\Delta t^2}{2\zeta}[\| p_{ch}^{n+1}\|^2-\| p_{ch}^{n}\|^2]
\nonumber\\
&&\hskip 1.7cm+\frac{\xi}{2} [\|\nabla\cdot \Auch^{n+1}\|^2-\|\nabla\cdot \Auch^{n}\|^2+\|\nabla\cdot(\Auch^{n+1}-\Auch^{n})\|^2]
+\frac{\Delta t^2}{2\zeta}\|p_{ch}^n-p_{ch}^{n-1}\|^2\nonumber\\
&&\hskip 2.4cm+\Delta t \langle\Auch^{n+1}\cdot\Anc,p_{mh}^{n+1}\rangle+\Delta t\frac{\alpha \sqrt{\mbox{d}}}{\sqrt{\mbox{trace$(\prod)$}}}\langle \nu^{n} P_\tau\Auch^{n+1},P_\tau\Auch^{n+1}\rangle
\nonumber\\
&&\hskip 2.4cm\leq \frac{\zeta}{2}\|\nabla\cdot(\Auch^{n+1}-\Auch^n)\|^2.
\end{eqnarray}
We rewrite \eqref{CS2_flux_identity} as
\begin{eqnarray}\label{CS2_flux_identity_2}
\frac{\rho^n(\Austar^n-\Auch^n)}{\Delta t}&=& -\phi_{ch}^n\nabla w_{ch}^{n+1},
 \end{eqnarray}
and take the inner product of \eqref{CS2_flux_identity_2} with $\Delta t \Austar^n$ to obtain
 by using the identity \eqref{identity}
\begin{eqnarray}\label{Decouple_NSCH_Energy S2_9}
&&\frac{1}{2}[\|\sigma^n\Austar^n\|^2-\|\sigma^n\Auch^n\|^2+\|\sigma^n(\Austar^n-\Auch^n)\|^2]=- \Delta t (\phi_{ch}^n\nabla w_{ch}^{n+1},\Austar^n).
 \end{eqnarray}
 Adding \eqref{Decouple_NSCH_Energy S2_8} and \eqref{Decouple_NSCH_Energy S2_9}, we obtain
 \begin{eqnarray}\label{Decouple_NSCH_Energy S2_11}
&&\frac{1}{2}[\|\sigma^{n+1}\Auch^{n+1}\|^2-\|\sigma^n\Auch^n\|^2+\|\sigma^n(\Austar^n-\Auch^n)\|^2+\|\sigma^n(\Auch^{n+1}-\Austar^n)\|^2]
+\frac{\Delta t^2}{2\zeta}[\| p_{ch}^{n+1}\|^2-\| p_{ch}^{n}\|^2]
\nonumber\\
&&\hskip 2.4cm+\frac{\xi}{2} [\|\nabla\cdot \Auch^{n+1}\|^2-\|\nabla\cdot \Auch^{n}\|^2]+\frac{\xi}{2}\|\nabla\cdot(\Auch^{n+1}-\Auch^{n})\|^2
+\frac{\Delta t^2}{2\zeta}\|p_{ch}^{n}-p_{ch}^{n-1}\|^2
\nonumber\\
&&\hskip 2.4cm+\Delta t\|\sqrt{2\nu^{n}}\mathbb{D}(\Auch^{n+1})\|^2+\Delta t\frac{\alpha \sqrt{\mbox{d}}}{\sqrt{\mbox{trace$(\prod)$}}}\langle \nu^{n} P_\tau\Auch^{n+1},P_\tau\Auch^{n+1}\rangle+\Delta t \langle\Auch^{n+1}\cdot\Anc,p_{mh}^{n+1}\rangle\nonumber\\
&&\hskip 2.4cm
\leq
-\Delta t (\phi_{ch}^n\nabla  w_{ch}^{n+1},\Austar^n)+\frac{\zeta}{2}\|\nabla\cdot(\Auch^{n+1}-\Auch^n)\|^2.
\end{eqnarray}

Then, we study the matrix part.
We take the inner product of \eqref{DCH_velocity_time} with $\Delta t \Aumh^{n+1}$ to get
  \begin{eqnarray}\label{Decouple_DCH_Energy_3}
  \Delta t \|\sqrt{\mK^{-1}}\Aumh^{n+1}\|^2=\Delta t(\nabla (p_{mh}^{n+1}-p_{mh}^n),\Aumh^{n+1})-\Delta t(\nabla p_{mh}^{n+1},\Aumh^{n+1})-\Delta t(\phi_{mh}^n\nabla w_{mh}^{n+1},\Aumh^{n+1}).
 \end{eqnarray}
From  \eqref{DCH_velocity_time},  \eqref{2Decouple_DCH_full_disretization_BJS_1_time} can be written as
 \begin{eqnarray}\label{Decouple_DCH_full_disretization_BJS_11_time}
 -(\Aumh^{n+1},\nabla q_h)+({\mathbb K}\nabla (p_{mh}^{n+1}-p_{mh}^n),\nabla q_h)+\beta\Delta t (\nabla p_{mh}^{n+1},\nabla q_h)-\langle\Auch^{n}\cdot\Anc,q_h\rangle = 0.
  \end{eqnarray}
We take $q_h=\Delta t p_{mh}^{n+1}$ in \eqref{Decouple_DCH_full_disretization_BJS_11_time} and utilize the identity \eqref{identity} to obtain
 \begin{eqnarray}\label{Decouple_DCH_Energy_4}
 &&-\Delta t (\Aumh^{n+1},\nabla p_{mh}^{n+1})+\frac{1}{2}\Delta t [\|\sqrt{\mK}\nabla p_{mh}^{n+1}\|^2-\|\sqrt{\mK}\nabla p_{mh}^n\|^2+\|\sqrt{\mK}\nabla (p_{mh}^{n+1}-p_{mh}^{n})\|^2]+\beta\Delta t^2 \|\nabla p_{mh}^{n+1}\|^2
 \nonumber\\
&&\hskip 3cm=\Delta t \langle\Auch^{n}\cdot \Anc, p_{mh}^{n+1} \rangle.\quad\quad
 \end{eqnarray}
Taking the sum of \eqref{Decouple_DCH_Energy_3} and \eqref{Decouple_DCH_Energy_4}, we get
  \begin{eqnarray}\label{2Decouple_DCH_Energy_5}
 &&\Delta t\|\sqrt{\mK^{-1}}\Aumh^{n+1}\|^2
+\frac{1}{2}\Delta t [\|\sqrt{\mK}\nabla p_{mh}^{n+1}\|^2-\|\sqrt{\mK}\nabla p_{mh}^n\|^2+\|\sqrt{\mK}\nabla (p_{mh}^{n+1}-p_{mh}^{n})\|^2]+\beta\Delta t^2 \|\nabla p_{mh}^{n+1}\|^2\nonumber\\
 &&\hskip  2.5cm=\Delta t(\nabla (p_{mh}^{n+1}-p_{mh}^n),\Aumh^{n+1})-\Delta t(\phi_{mh}^n\nabla w_{mh}^{n+1},\Aumh^{n+1})+\Delta t\langle\Auch^{n}\cdot \Anc, p_{mh}^{n+1} \rangle.
 \end{eqnarray}
Now, we estimate the term $(\nabla (p_{mh}^{n+1}-p_{mh}^n),\Aumh^{n+1})$. Combining
\begin{eqnarray}
\Delta t|(\nabla (p_{mh}^{n+1}-p_{mh}^n),\Aumh^{n+1})|\leq \Delta t \|\sqrt{\mK^{-1}}\Aumh^{n+1}\|^2+\frac{1}{4}\Delta t \|\sqrt{\mK}\nabla (p_{mh}^{n+1}-p_{mh}^{n})\|^2,\label{Decouple_DCH_Energy_6}
\end{eqnarray}
then, we obtain from \eqref{2Decouple_DCH_Energy_5}
\begin{eqnarray}
&&\frac{1}{2}\Delta t [\|\sqrt{\mK}\nabla p_{mh}^{n+1}\|^2-\|\sqrt{\mK}\nabla p_{mh}^n\|^2]+\frac{1}{4}\Delta t\|\sqrt{\mK}\nabla (p_{mh}^{n+1}-p_{mh}^{n})\|^2+\beta\Delta t^2 \|\nabla p_{mh}^{n+1}\|^2\nonumber\\
 &&\hskip  2.5cm\leq-\Delta t(\phi_{mh}^n\nabla w_{mh}^{n+1},\Aumh^{n+1})+\Delta t\langle\Auch^{n}\cdot \Anc, p_{mh}^{n+1} \rangle.\label{2Decouple_DCH_Energy_6}
\end{eqnarray}

Adding~\eqref{2Decouple_CH_Energy_3}, ~\eqref{Decouple_NSCH_Energy S2_11} and \eqref{2Decouple_DCH_Energy_6} together, we obtain
\begin{eqnarray}\label{Decouple_NSCHDCH_Energy S2}
&&\frac{1}{2}[\|\sigma^{n+1}\Auch^{n+1}\|^2-\|\sigma^n\Auch^n\|^2]
+\frac{\gamma\epsilon}{2}[\|\nabla\phi_h^{n+1}\|^2-\|\nabla\phi_h^n\|^2]
+\gamma(F(\phi_h^{n+1})-F(\phi_h^n),1)
\nonumber\\
&&\hskip 2.4cm
+\frac{\Delta t^2}{2\zeta}[\| p_{ch}^{n+1}\|^2-\| p_{ch}^{n}\|^2]+\frac{\xi}{2} [\|\nabla\cdot \Auch^{n+1}\|^2-\|\nabla\cdot \Auch^{n}\|^2]+\frac{\xi}{2}\|\nabla\cdot(\Auch^{n+1}-\Auch^{n})\|^2
\nonumber\\
&&\hskip 2.4cm+\frac{1}{2}\Delta t [\|\sqrt{\mK}\nabla p_{mh}^{n+1}\|^2-\|\sqrt{\mK}\nabla p_{mh}^n\|^2]
+\frac{\Delta t^2}{2\zeta}\|p_{ch}^{n}-p_{ch}^{n-1}\|^2
+\frac{1}{4}\Delta t\|\sqrt{\mK}\nabla (p_{mh}^{n+1}-p_{mh}^{n})\|^2\nonumber\\
&&\hskip 2.4cm
 +\Delta t\|\sqrt{2\nu^{n}}\mathbb{D}(\Auch^{n+1})\|^2
+\Delta tM\|\nabla w_h^{n+1}\|^2+\beta\Delta t^2 \|\nabla p_{mh}^{n+1}\|^2
+\frac{\gamma\epsilon}{2}\|\nabla\phi_h^{n+1}-\nabla\phi_h^n\|^2
\nonumber\\
&&
\hskip 2.4cm
+\frac{1}{2}[\|\sigma^n(\Austar^n-\Auch^n)\|^2+\|\sigma^n(\Auch^{n+1}-\Austar^n)\|^2]
+\Delta t\frac{\alpha \sqrt{\mbox{d}}}{\sqrt{\mbox{trace$(\prod)$}}}\langle \nu^{n} P_\tau\Auch^{n+1},P_\tau\Auch^{n+1}\rangle
\nonumber\\
&&\hskip 2.4cm \leq\frac{\zeta}{2}\|\nabla\cdot(\Auch^{n+1}-\Auch^n)\|^2
+\Delta t\langle(\Auch^n-\Auch^{n+1})\cdot\Anc,p_{mh}^{n+1}\rangle.
\end{eqnarray}

Now, we estimate the last interface term in the above equation.
Using using Lemma~\ref{lemma31} and the triangle inequality,
\begin{eqnarray}\label{S2_estimate_interface}
&&\Delta t|\langle\left(\Auch^{n}-\Auch^{n+1}\right)\cdot \Anc,p_{mh}^{n+1}\rangle|
\leq
C\Delta t\|\Auch^{n}-\Auch^{n+1}\|_{\BX_{div}}\|\nabla p_{mh}^{n+1}\|\nonumber\\
&&\hskip 4.5cm\leq \frac{1}{4}\min\{\rho_1,\rho_2\}\|\Auch^{n}-\Auch^{n+1}\|_{\BX_{div}}^2+\tilde{C}\Delta t^2 \|\nabla p_{mh}^{n+1}\|^2\nonumber\\
&&\hskip 4.5cm= \frac{1}{4}\min\{\rho_1,\rho_2\}\|\Auch^{n}-\Auch^{n+1}\|^2+\frac{1}{4}\min\{\rho_1,\rho_2\}\|\nabla\cdot(\Auch^{n}-\Auch^{n+1})\|^2+\tilde{C}\Delta t^2 \|\nabla p_{mh}^{n+1}\|^2\nonumber\\
&&\hskip 4.5cm\leq
\frac{1}{4}\|\sigma^n(\Auch^{n+1}-\Auch^{n})\|^2+\frac{1}{4}\min\{\rho_1,\rho_2\}\|\nabla\cdot(\Auch^{n}-\Auch^{n+1})\|^2+\tilde{C}\Delta t^2 \|\nabla p_{mh}^{n+1}\|^2.\quad
\end{eqnarray}
On the other hand, we derive from the triangle inequality that
 \begin{eqnarray}\label{Decouple_NSCH_Energy S2_10}
&&-\frac{1}{2}[\|\sigma^n(\Austar^n-\Auch^n)\|^2+\|\sigma^n(\Auch^{n+1}-\Austar^n)\|^2]\leq -\frac{1}{4}\|\sigma^n(\Auch^{n+1}-\Auch^n)\|^2.
 \end{eqnarray}

Adding \eqref{Decouple_NSCHDCH_Energy S2}, \eqref{S2_estimate_interface} and \eqref{Decouple_NSCH_Energy S2_10}, we obtain
\begin{eqnarray}\label{Decouple_NSCHDCH_Energy_all S2}
&&\mathrm{\mathcal{E}}^{n+1}-\mathrm{\mathcal{E}}^{n}
 \leq-\Delta t\|\sqrt{2\nu^{n}}\mathbb{D}(\Auch^{n+1})\|^2
-\Delta tM\|\nabla w_h^{n+1}\|^2-\frac{\Delta t^2}{2\zeta}\|p_{ch}^{n}-p_{ch}^{n-1}\|^2
-\frac{1}{4}\Delta t\|\sqrt{\mK}\nabla (p_{mh}^{n+1}-p_{mh}^{n})\|^2\nonumber\\
&&\hskip 2.4cm
-\frac{\gamma\epsilon}{2}\|\nabla\phi_h^{n+1}-\nabla\phi_h^n\|^2-\Delta t\frac{\alpha \sqrt{\mbox{d}}}{\sqrt{\mbox{trace$(\prod)$}}}\langle \nu^{n} P_\tau\Auch^{n+1},P_\tau\Auch^{n+1}\rangle
\nonumber \\
&&\hskip 2.5cm-(\beta-\tilde{C})\Delta t^2 \|\nabla p_{mh}^{n+1}\|^2-\frac{1}{2}(\xi-\zeta-\frac{1}{2}\min\{\rho_1,\rho_2\})\|\nabla\cdot(\Auch^{n+1}-\Auch^{n})\|^2.
\quad
\end{eqnarray}

If we now impose $\xi\geq \zeta+\frac{1}{2}\min\{\rho_1,\rho_2\}$ and $\beta\geq 2\tilde{C}$ which only depends on
the geometry of $\Omega_m$, $\Omega_c$, $\rho_1$ and $\rho_2$, then one leads to the energy stability and  complete the proof
of Theorem \ref{S2_thm_energy_Decouple}. \mbox{}\hfill$\Box$

\section{Numerical results}\label{sec(5)}
In this section, we will use three numerical examples to
illustrate the features of proposed model and numerical methods.
The first example is provided to illustrate  the convergence and accuracy.
The second test is designed to verify that the proposed algorithm \eqref{2Decouple_DCH_full_disretization_BJS_2_time}-\eqref{2Decouple_NSCH_full_disretization_BJS_5_time}   obeys the energy dissipation of the CHNSD model \eqref{Darcy_law_time_BJ}-\eqref{phase_interface_con_4}.
The last experiment  presents the simulation of
a lighter bubble rising through the interface driven by buoyancy forces.
For all examples,
we  employ the  celebrated  Taylor-Hood elements for the Navier-Stokes equation and  linear elements for the Darcy equation.
For the single Cahn-Hilliard equation in the coupling free flow and porous media,
we consider the quadratic elements.\\

\noindent
\textbf{Example 1: Convergence and accuracy.} Consider the  CHNSD model  on~$\Omega=[0,1]\times[0,2]$~where~$\Omega_m=[0,1]\times[0,1]$~and $\Omega_c=[0,1]\times[1,2]$.
Set~$\nu=1$,~$\rho_1=1$,~$\rho_2=3$,~$M_m=1$,~$\gamma=1$,~$\epsilon=1$,~$M_c=1$~,~$\mathbb{K}=\mathbb{I}$,$\beta=5$, and $\xi=5$.
The simulation is performed out at terminational time $T=0.2$.
The exact solutions are chosen as:
\begin{eqnarray}\label{accuray_example}
\left\{ \begin{array}{l}
\phi=g(x)g(y)\cos(\pi t),\\
p_m=g(x)g_m(y)\cos(\pi t), \\
\Auc=[x^2(y-1)^2,~-\frac{2}{3}x(y-1)^3]^T\cos(\pi t),  \\
p_c=\cos(\pi t)g(x)g_c(y),
\end{array} \right.
\label{exact_solution_example_2}
\end{eqnarray}
where~$g(x)=16x^2(x-1)^2, g(y)=16y^2(y-2)^2,g_m(y)=16y^2(y-1)^2,
g_c(y)=16(y-1)^2(y-2)^2$. The boundary condition functions and the
source terms can be computed based on the exact solutions.

To examine the accuracy of proposed scheme,  we compute the pointwise convergence rate and  define the rate of convergence in space as follows
\begin{eqnarray*}
\mbox{order}_h=\frac{\log(|e_{v,h_j}|/|e_{v,h_{j+1}}|)}{\log(h_j/h_{j+1})}=\frac{\log(|v_{h_j}^n-v(t_n)|/|v_{h_{j+1}}^n-v(t_n)|)}{\log(h_j/h_{j+1})},\quad v=\phi,\,p_m,\,\Auc,\,p_c,
\end{eqnarray*}
where $|\cdot|$ denotes the $L^2$ and $H^1$ norm errors with $\|\cdot\|$ and $\|\cdot\|_1$, $v_{h_j}$ is the numerical solution with spatial mesh size $h_j$.
Tables~\ref{table_2_time_BJ} and \ref{table_3_time_BJ} list the $L^2$- and~$H^1$-norm errors of the  phase variable, pressure and  velocity of the
 designed  decoupled linearized  numerical schemes, in which
 a uniform time partition $\Delta t=2.5\times 10^{-4}$ is used.
The numerical results in the two tables clearly show the optimal convergence rates for constructed numerical scheme for all presented error norms in space.

\begin{figure}[!ht]
\centering
\setlength{\abovecaptionskip}{-0.1cm}
\setlength{\belowcaptionskip}{-0.2cm}
\includegraphics[width=3.0in]{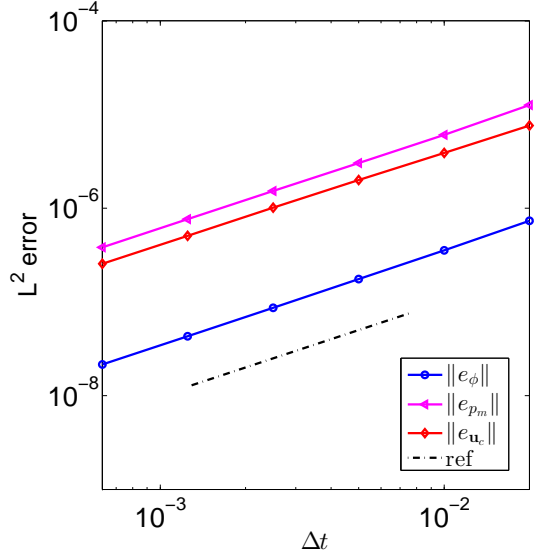}
\caption{Log-Log  plots of the $L^2$  error  norms with   different time step size $\Delta t$.} \label{fig:TimeOrder}
\end{figure}

To illustrate the order  of convergence with respect to the time step $\Delta t$, we introduce the following convergence rate for $L^2$-norm error,
\begin{eqnarray*}
  \mbox{order}_{\Delta t}=\frac{\log(\|v_h^{\Delta t}-v_h^{\Delta t/2}\|/\|v_h^{\Delta t/2}-v_h^{\Delta t/4}\|)}{\log(2)}, \quad v=\phi,\,p_m,\,\Auc.
\end{eqnarray*}
The $L^2$-norm errors are shown in Figure \ref{fig:TimeOrder} with fixing spatial mesh size $h=\frac{1}{32}$ and varying partition $\Delta t=0.02/2^k$, $k=0,1,\ldots,5$, which indicates that the proposed numerical method can achieve the first order accuracy in time for variables $\phi$, $p_m$ and $\Auc$.
\\

\begin{table}[hbt]
\begin{center}
\begin{tabular}{c|c|c|c|c|c|c}
\hline  h   &  $\|e_{\Auc}\|$    &order    &$\|e_{\Auc}\|_1$   &order    &  $\|e_{p_c}\|$    &order  \\[0.5ex]\hline
1/4   &6.3163E-3   &   &5.4761E-2    &       &1.0353E-1   &         \\
  1/8    &6.7679E-4   &3.22   &6.5007E-3      &3.07  &3.6981E-2  &1.84    \\
  1/16   &8.6286E-5   & 2.97    &1.2233E-3      &   2.41   &1.0093E-2  &1.87  \\
  1/32  &1.2159E-5   &2.83    &3.0139E-4       & 2.02    &2.6242E-3 &1.94   \\\hline
\end{tabular}
\end{center}
\caption{The order of convergence   in space for    error norms  for $\Auc$ and $p_c$.} \label{table_2_time_BJ}
\end{table}

\begin{table}[hbt]
\begin{center}
\begin{tabular}{c|c|c|c|c|c|c|c|c}
\hline  h   &  $\|e_{\phi}\|$    &order    &$\|e_{\phi}\|_1$   &order    &  $\|e_{p_m}\|$    &order  &$\|e_{p_m}\|_1$   &order\\[0.5ex]\hline
1/4    &2.9077E-2    &         &5.8721E-1     &        &9.9369E-2      &    &8.2818E-1   &\\
  1/8   &3.2129E-3  &3.18    &1.5886E-1     &1.89    &2.6547E-2    &1.90     &4.5968E-1  &0.849\\
  1/16    &3.6121E-4  & 3.15   &4.0798E-2      &1.96 &5.8678E-3    &2.18   & 2.3371E-1  &0.976\\
  1/32  &4.3424E-5   & 3.06  &1.0277E-2    &1.99   &1.4238E-3    &2.04     &1.1772E-1   &0.989 \\\hline
\end{tabular}
\end{center}
\caption{The order of convergence in space for  error norms for $\phi$ and $p_m$.} \label{table_3_time_BJ}
\end{table}

For the diffuse interface problem,  adaptive mesh refinement is preferable for the computation of different dynamics, due to the fact at leat four  grid elements  are required for accuracy over the width of interface \cite{JKim_KKang_JLowengrub_JCP_2004}.
 Therefore, we use an  adaptive mesh strategy in this simulation.
 In the following numerical experiments, the root-level mesh is taken to be uniform with $h=\frac{1}{32}$.  Starting with this base mesh, mesh refinement is performed.

\noindent
\textbf{Example 2: Shape relaxation and energy dissipation.}
We consider the evolution of a square shaped circle bubble in the domain~$\Omega=[0, 1]\times[0, 2]$ with~$\Omega_c=[0,1]\times[0,1]$~and~$\Omega_m=[1,2]\times[0,1]$.
The parameters are chosen~$M=0.1$,~$\gamma=0.01$,~$\epsilon=0.02$,~$\nu=1$ and~$\mathbb{K}=0.05\mathbb{I}$.
The initial velocity, pressure and chemical potential are set to zero.

\begin{figure}[!ht]
\centering
\setlength{\abovecaptionskip}{-0.1cm}
\setlength{\belowcaptionskip}{-0.2cm}
\includegraphics[width=2.5in]{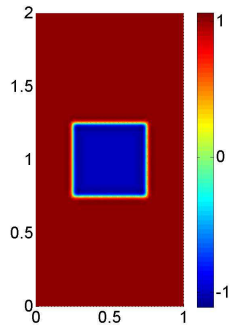}
\caption{Contour plots of the initial bubble.} \label{fig:initialphiEx2}
\end{figure}

\begin{figure}[ht]
\centering
\subfigure
{
\label{fig5:subfig:ea}
{\includegraphics[width=2.0in]{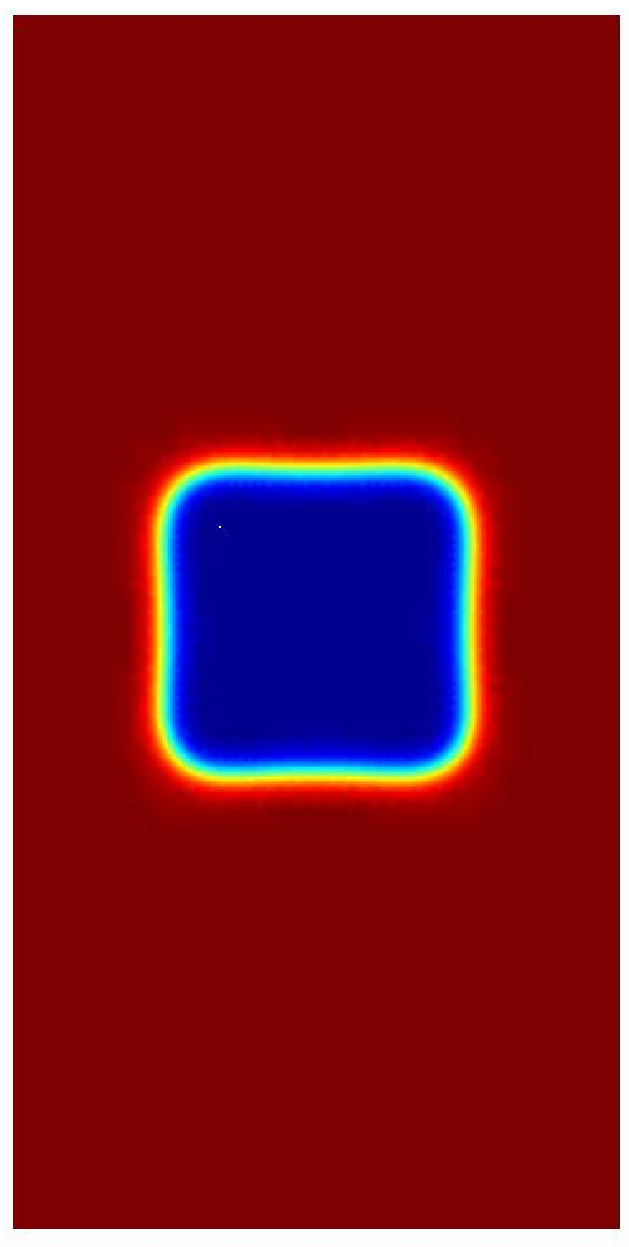}}}
\hskip -0.7in
\subfigure
{
\label{fig5:subfig:eb} 
{\includegraphics[width=2.0in]{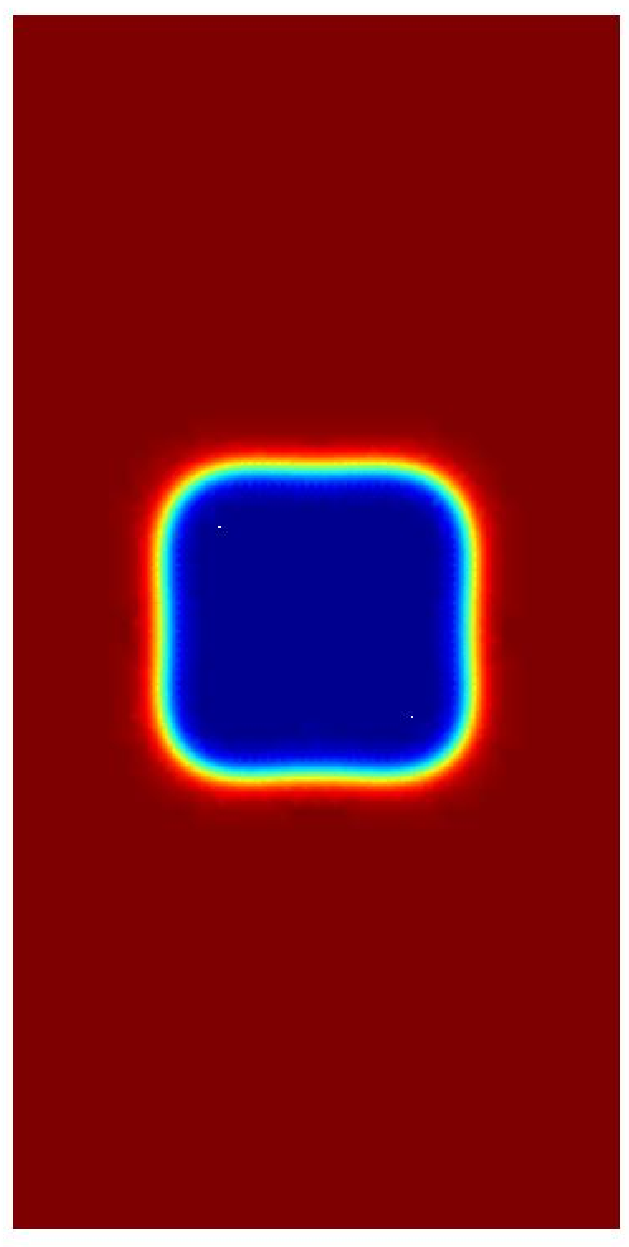}}}
\hskip -0.7in
\subfigure
{
\label{fig5:subfig:ec} 
{\includegraphics[width=2.0in]{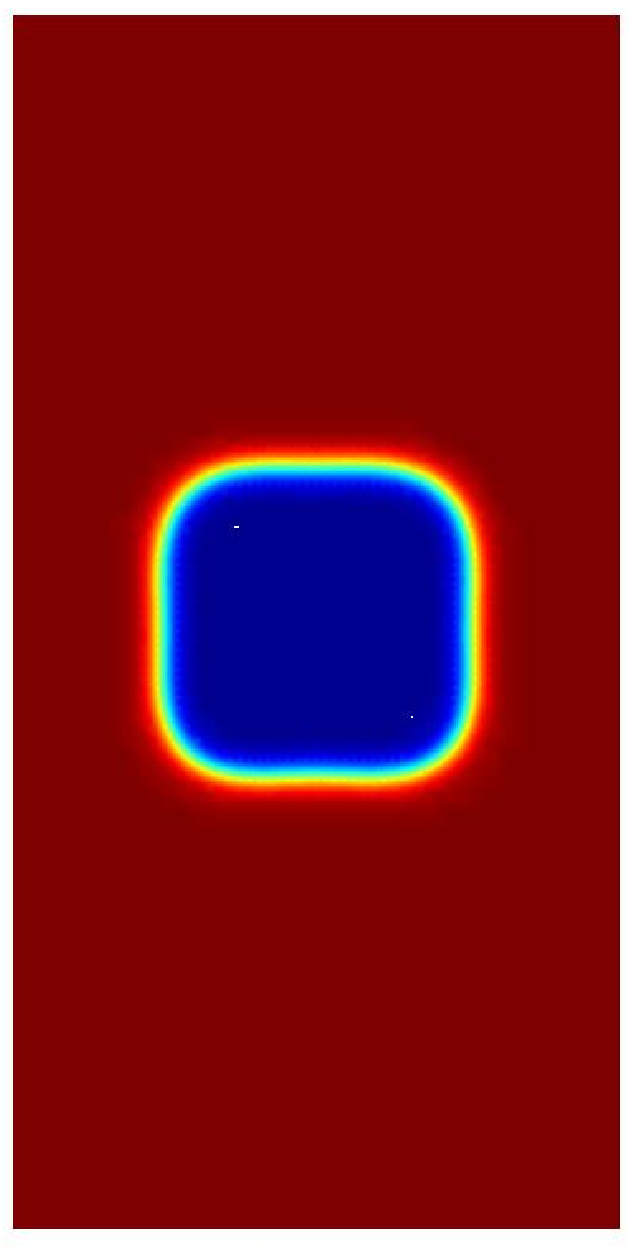}}}
\hskip -0.7in
\subfigure
{\label{fig5:subfig:ee} 
{\includegraphics[width=2.0in]{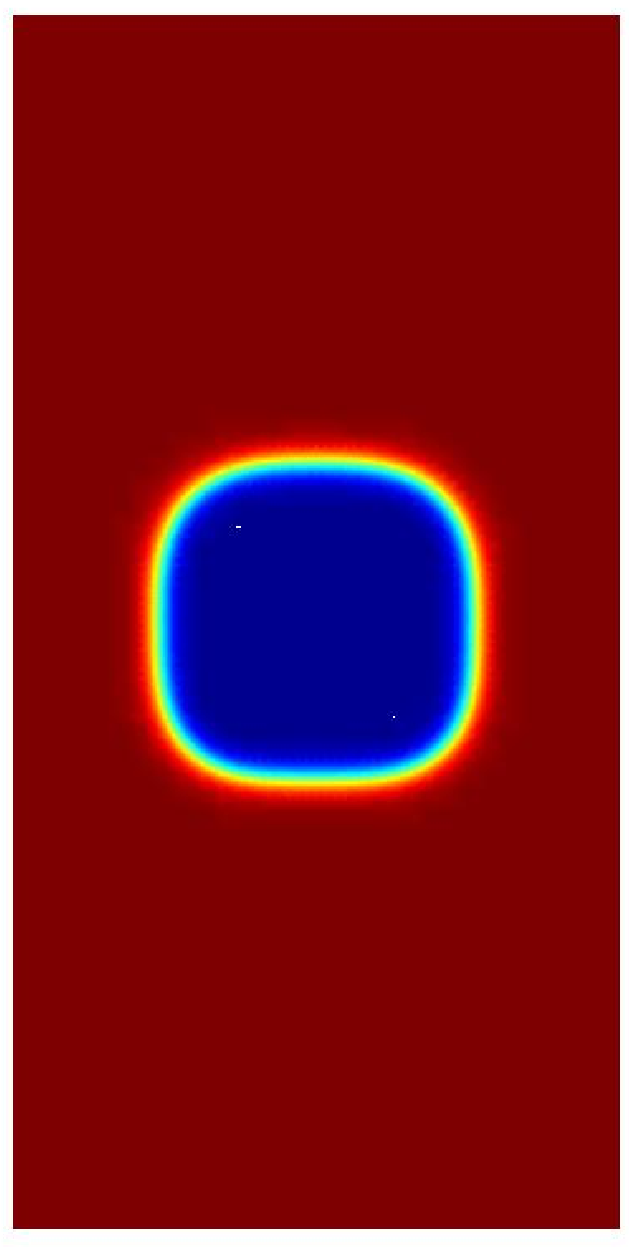}}}
\vskip -0.15in
\subfigure
{\label{fig5:subfig:ef} 
{\includegraphics[width=2.0in]{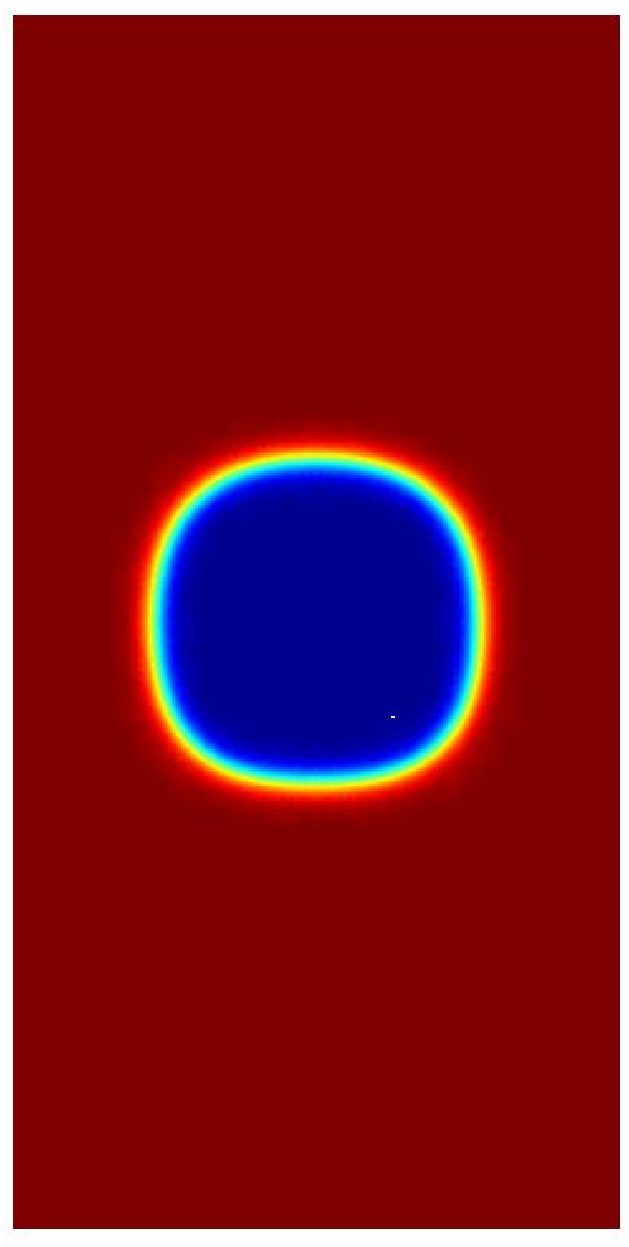}}}
\hskip -0.7in
\subfigure
{\label{fig5:subfig:eg} 
{\includegraphics[width=2.0in]{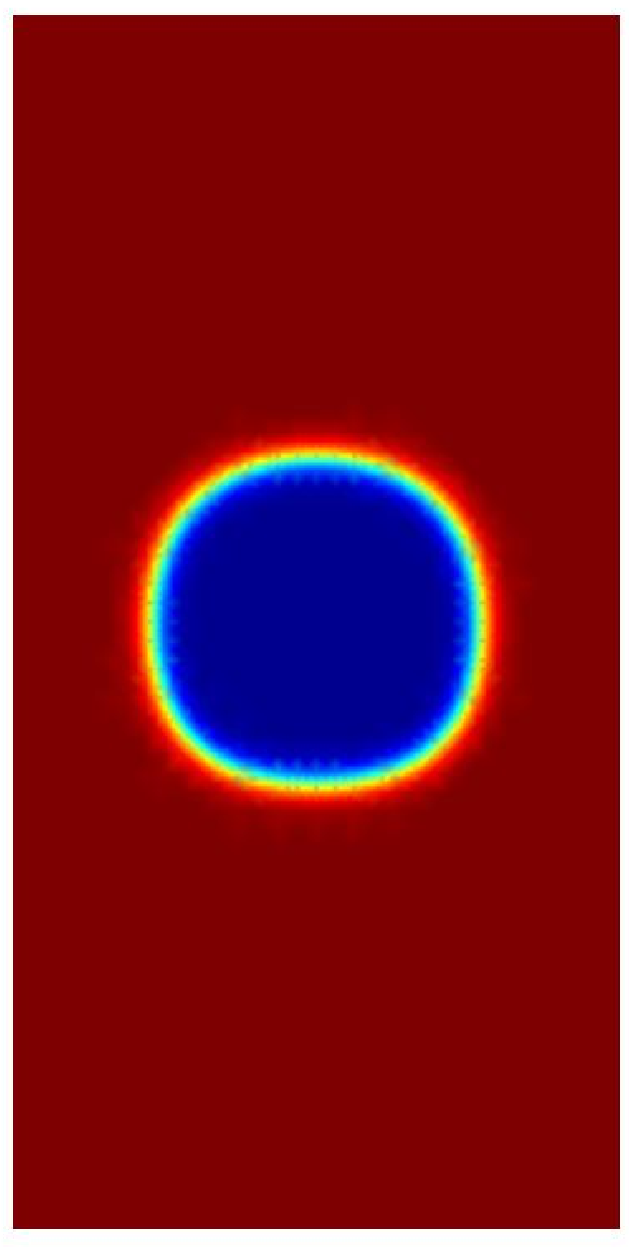}}}
\hskip -0.7in
\subfigure
{\label{fig5:subfig:edb} 
{\includegraphics[width=2.0in]{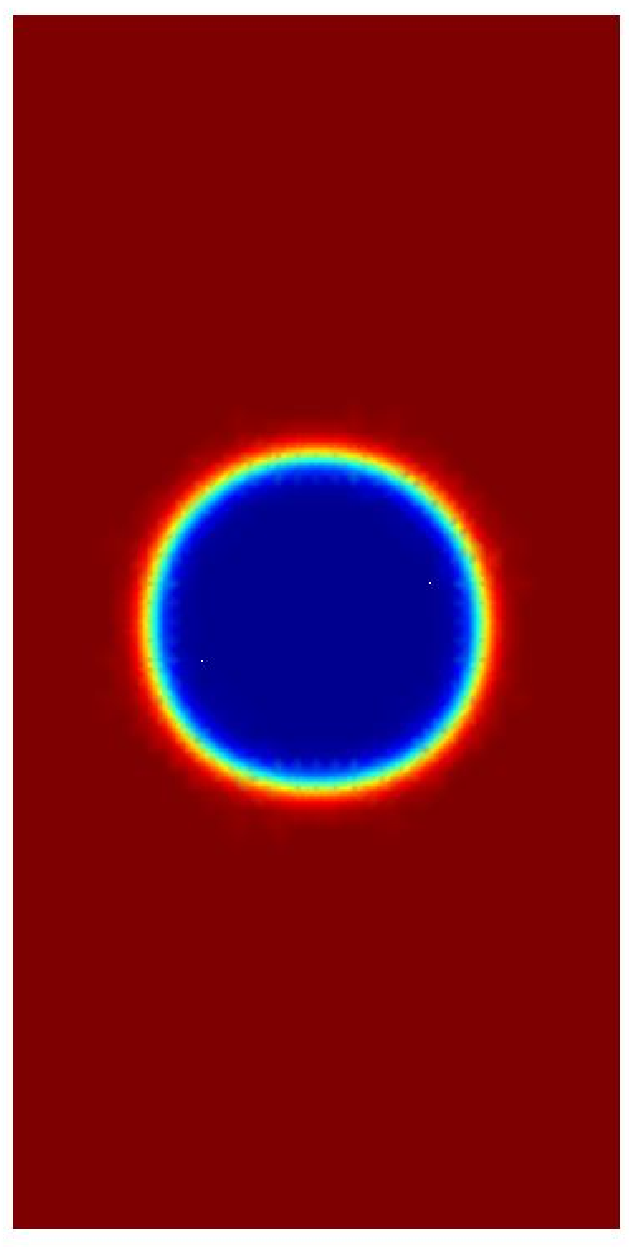}}}
\hskip -0.7in
\subfigure
{\label{fig5:subfig:elc} 
{\includegraphics[width=2.0in]{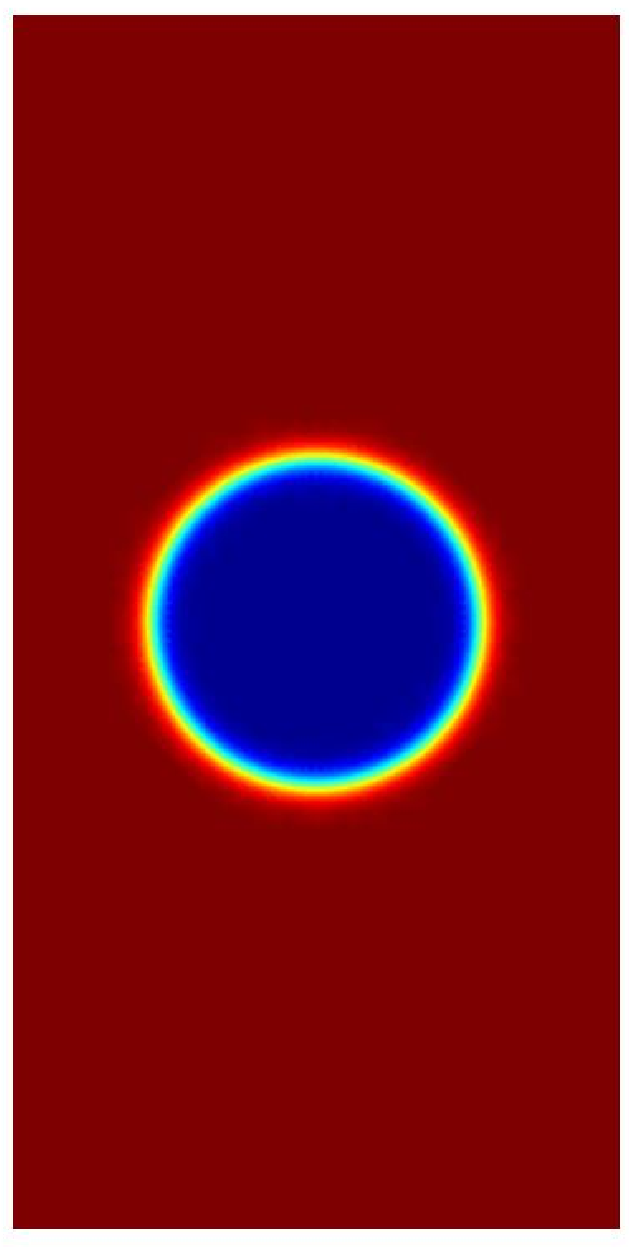}}}
\vskip -0.25in
\caption{The dynamics of a square shape bubble with density ratio 1:50.
All the sub-figures are indexed from left to right row by row as follows: : (a)~$t=0.1$, (b)~$t=0.2$, (c)~$t=0.3$, (d)~$t=0.5$, (e)~$t=0.8$, (f)~$t=1.0$, (g)~$t=1.5$, (h)~$t=10.0$.}
\label{fig:Energy150Ex2}
\end{figure}

\begin{figure}[ht]
\centering
\setlength{\abovecaptionskip}{-0.1cm}
\setlength{\belowcaptionskip}{-0.5cm}
\subfigure[$\rho_1:\rho_2=1:5$]
{\label{fhfg:subfig:a}
{\includegraphics[width=2.5in]{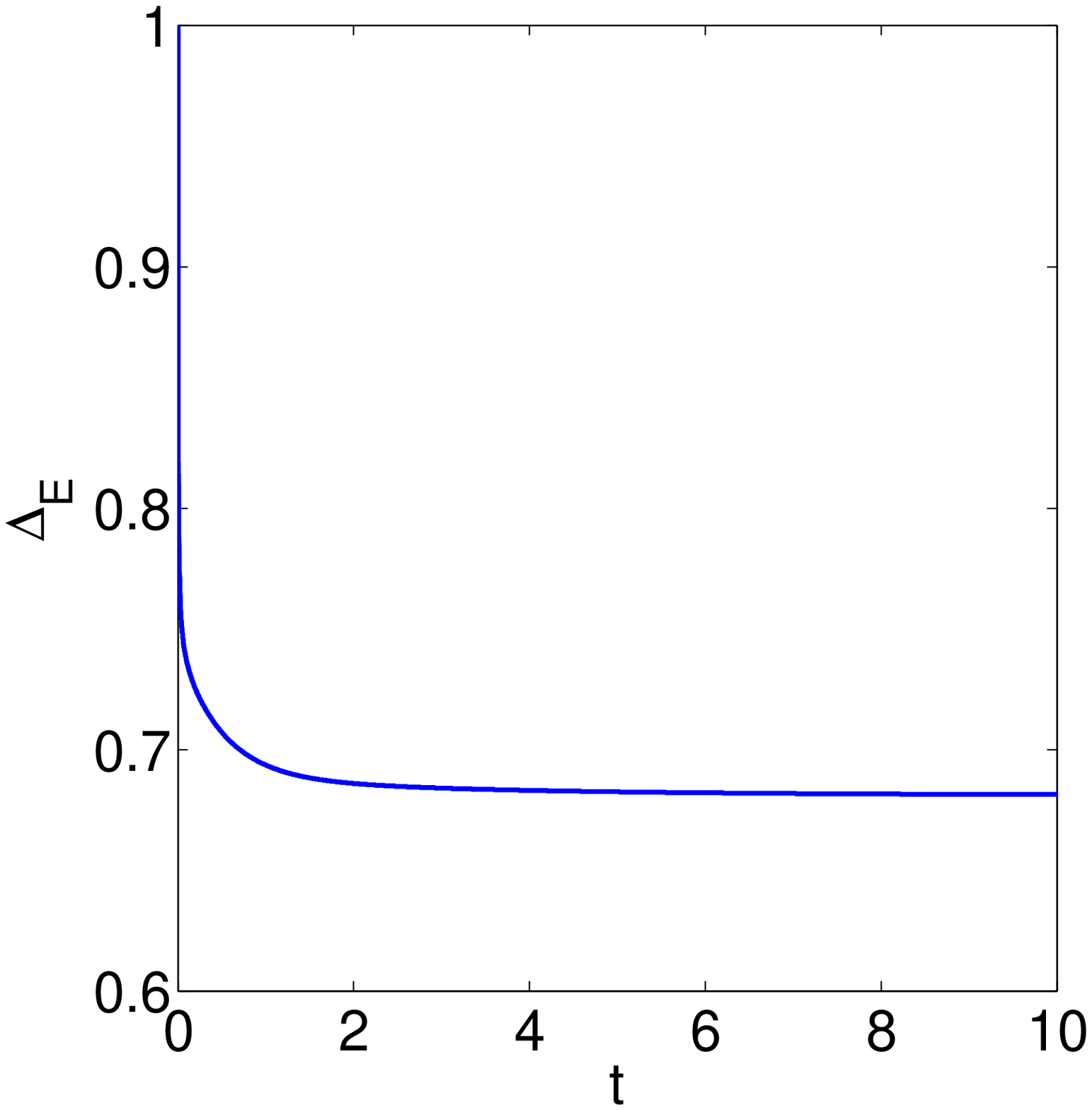}}}
\subfigure[$\rho_1:\rho_2=1:50$]
{\label{fhfg:subfig:b}
{\includegraphics[width=2.5in]{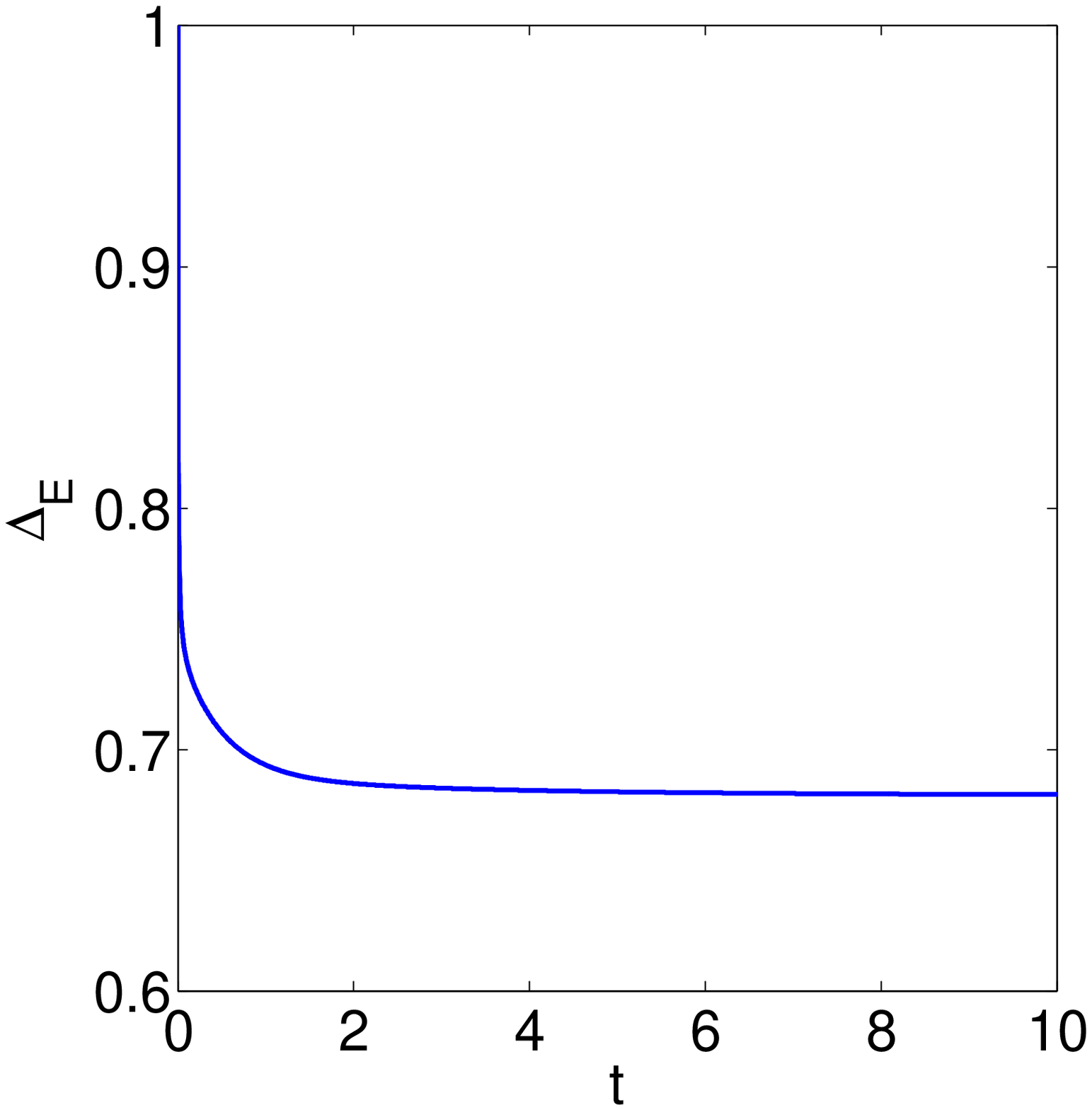}}}
\caption{The evolution of discrete energy of two numerical schemes.} \label{fig:Energy1}
\end{figure}

A uniform time partition with the time step-size~$\Delta t=0.005$~is used in this simulation.
Figure~\ref{fig:initialphiEx2} shows the  initial shape of the bubble.
 Figure~\ref{fig:Energy150Ex2} shows   the dynamic of  the square relaxing to a circular shape under the effect of surface tension  for the density ratio $\rho_1:\rho_2=1:50$ by using the proposed  decoupled numerical method.
The corresponding relative discrete energy~$\Delta_E=E^{n}/E^0$~is presented in Figure~\ref{fig:Energy1}.
We can easily observe that the discrete energy  is non-increasing for different density ratio cases, which is consistent with the theoretical result,
and validates the interface conditions~\eqref{phase_interface_con_1}-\eqref{phase_interface_con_4}.

\begin{figure}[!ht]
\centering
\setlength{\abovecaptionskip}{-0.1cm}
\setlength{\belowcaptionskip}{-0.2cm}
\subfigure[Initial phase function]{
\label{figIni:subfig:a}
\includegraphics[width=2.5in]{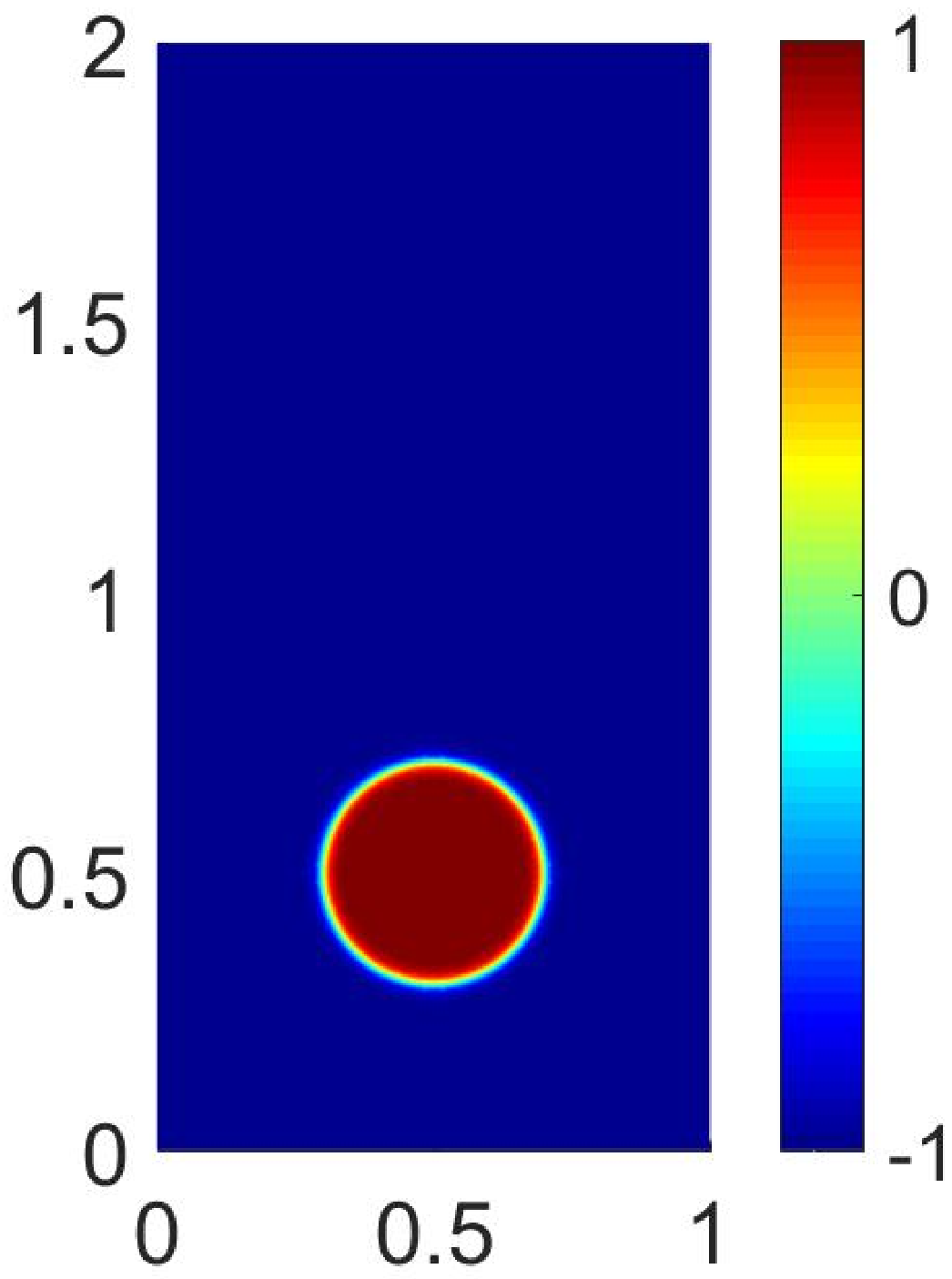}}
\subfigure[Initial adaptive mesh]{
\label{figIniA:subfig:b}
\includegraphics[width=2.5in]{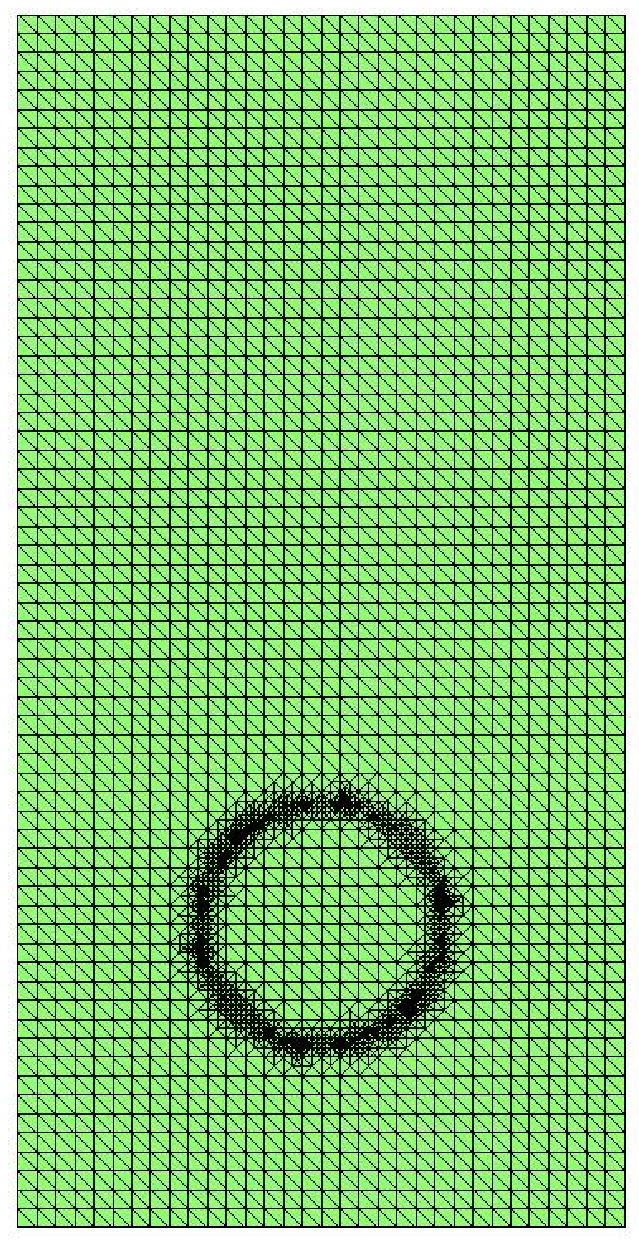}}
\caption{Contour plots of the initial bubble.} \label{fig:initialphiEx4}
\end{figure}

\noindent
\textbf{Example 3: Buoyancy-driven flow.}
In this experiment, we simulate a light bubble rises in a heavier medium in order to validate the efficiency of proposed numerical method with respect to different density variations.
Here, the karst geometry is modelled by a long tube~$\Omega = [0, 1]\times[0, 2]$~with the conduit~$\Omega_c = [0, 1]\times[0, 1]$~and porous media~$\Omega_m = [0, 1]\times[1, 2]$. The interface boundary is at~$[0,1]\times\{1\}$.
We set~$M=0.01$,~$\gamma=0.01$,~$\epsilon=0.01$,~$\nu=1$, and~$\mathbb{K}=0.05\mathbb{I}$.
The initial velocity and pressure are set to be zero and initial phase function is given by
\begin{eqnarray}\label{Buoyancy_initialphi}
\phi_c^0(x,y)=\tanh\left((0.2-\sqrt{(x-0.5)^2+(y-0.5)^2})/(\sqrt{2}\epsilon)\right).
\end{eqnarray}
Figure \ref{fig:initialphiEx4} shows the initial position of the bubble.

\begin{figure}[!htp]
\centering
\subfigure[$t=1.0$]{
\label{figr1:subfig:d} 
\includegraphics[width=2in]{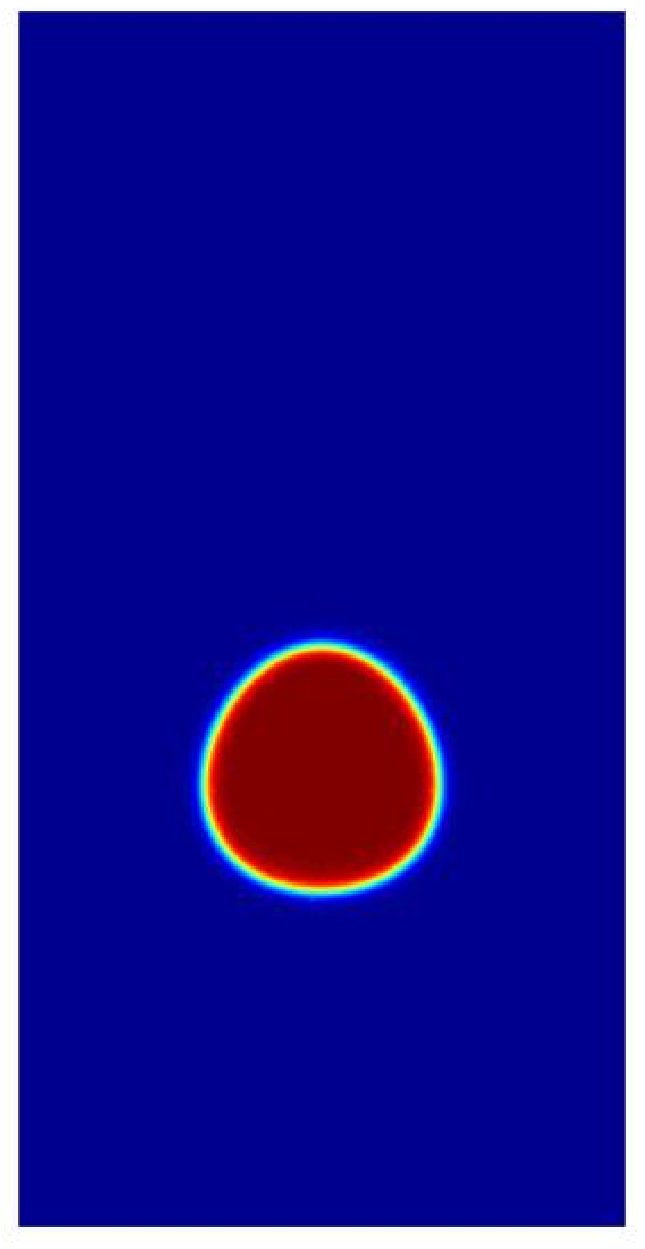}}
\hskip -0.7in
\subfigure[$t=2.0$]{
\label{figsr1:subfig:d} 
\includegraphics[width=2in]{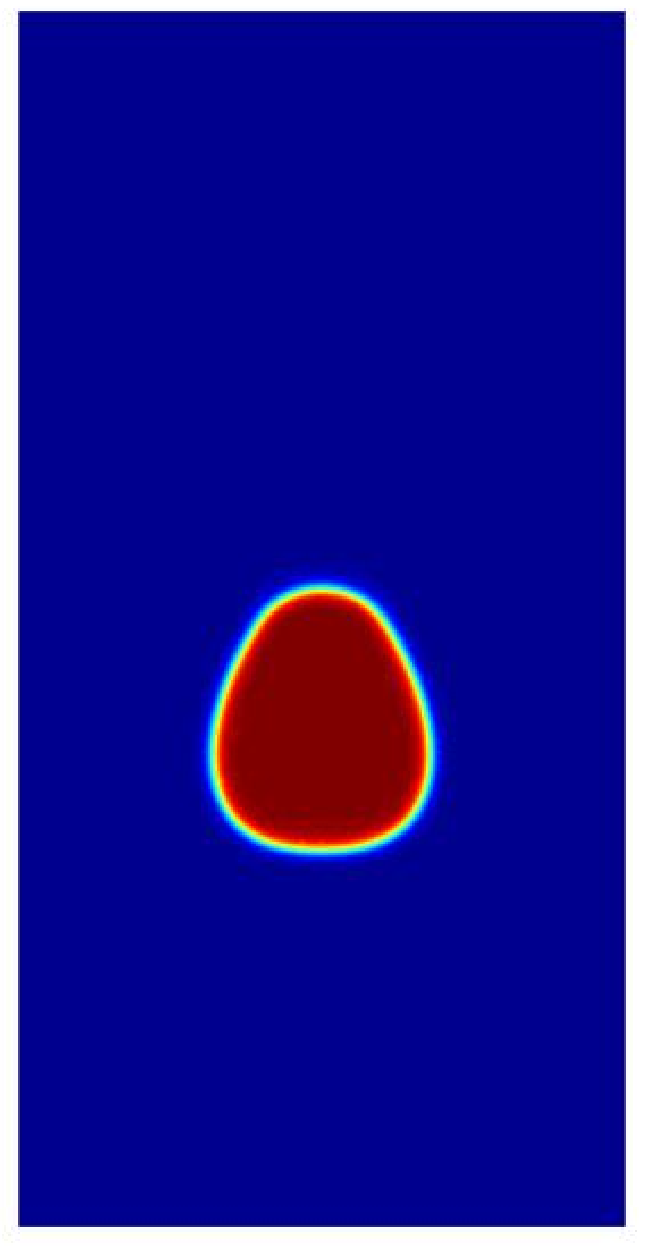}}
\hskip -0.7in
\subfigure
{\label{figr1:subfig:f} 
\includegraphics[width=2in]{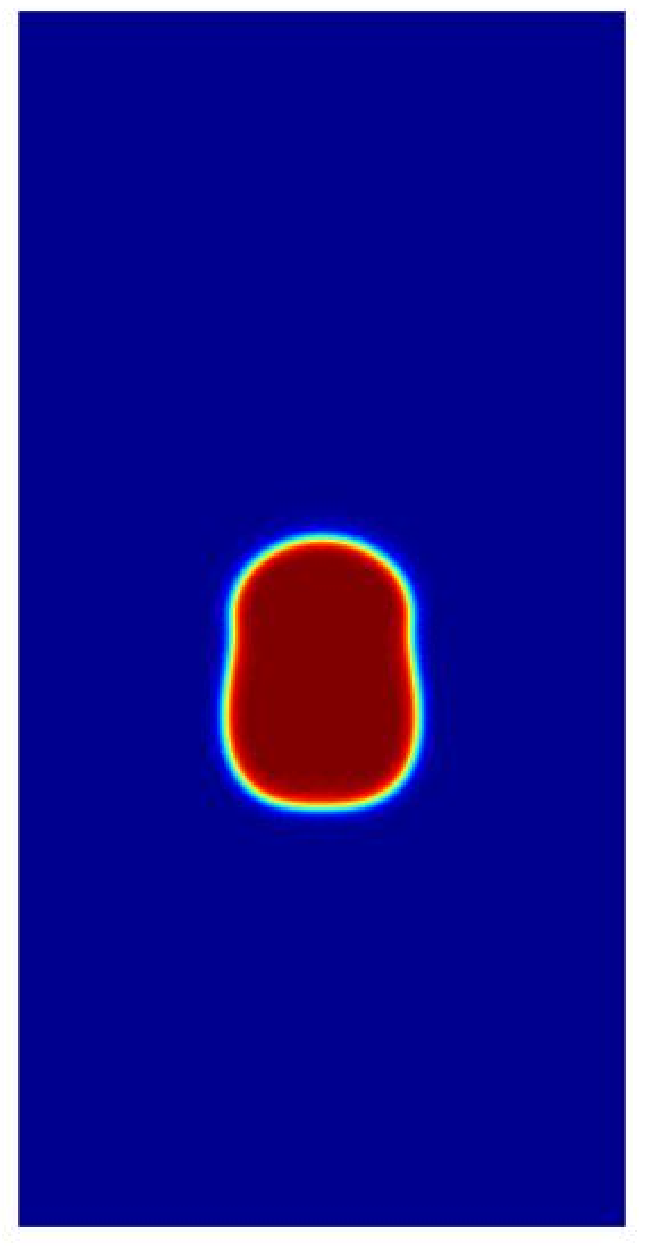}}
\hskip -0.7in
\subfigure[$t=4.0$]{
\label{figr1:subfig:g} 
\includegraphics[width=2in]{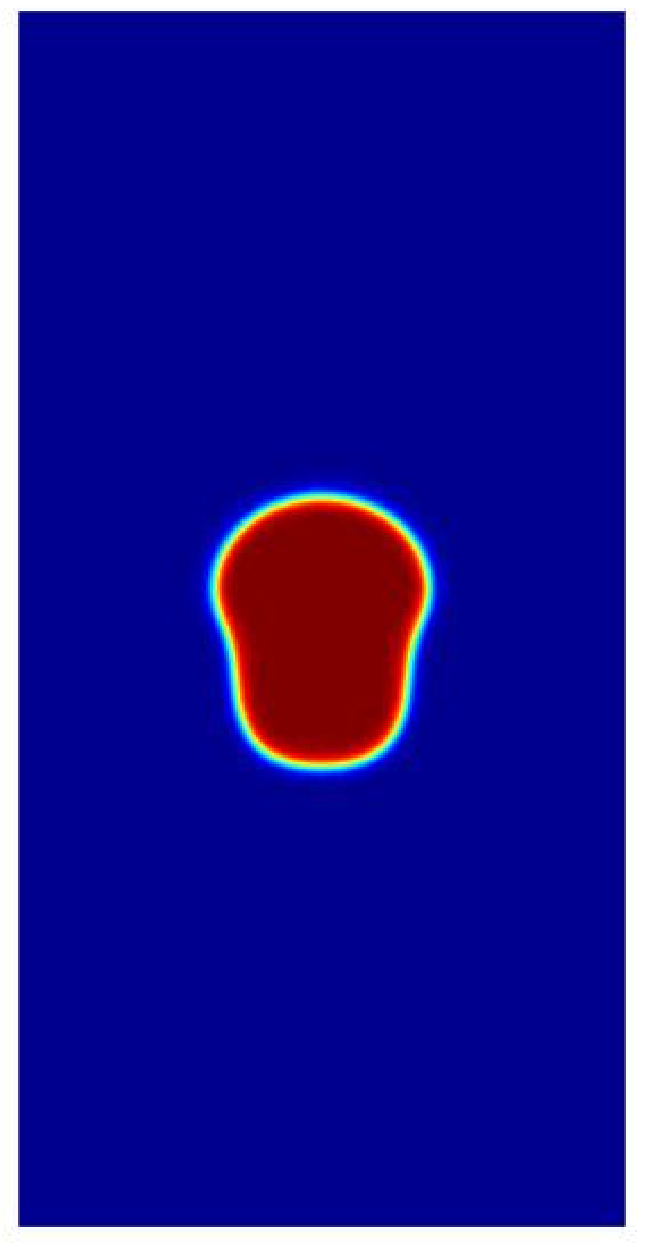}}
\vskip -0.40in
\subfigure
{\label{figr1:subfig:h} 
\includegraphics[width=2in]{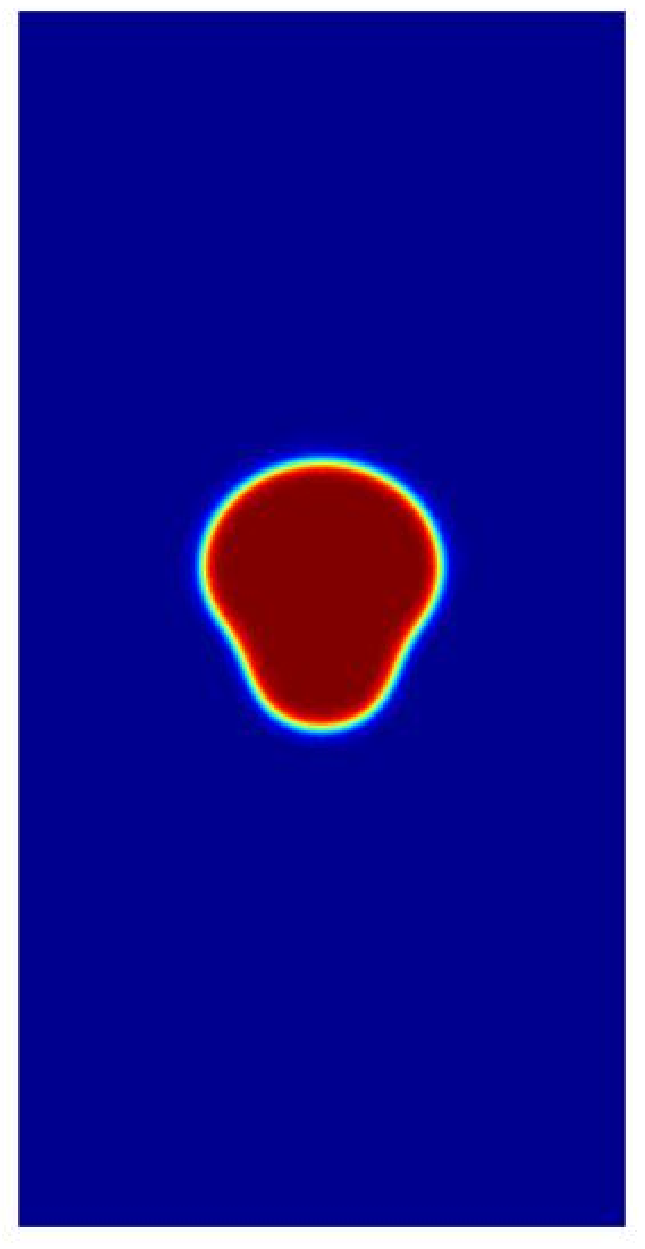}}
\hskip -0.7in
\subfigure
{\label{figr1:subfig:ll} 
\includegraphics[width=2.00in]{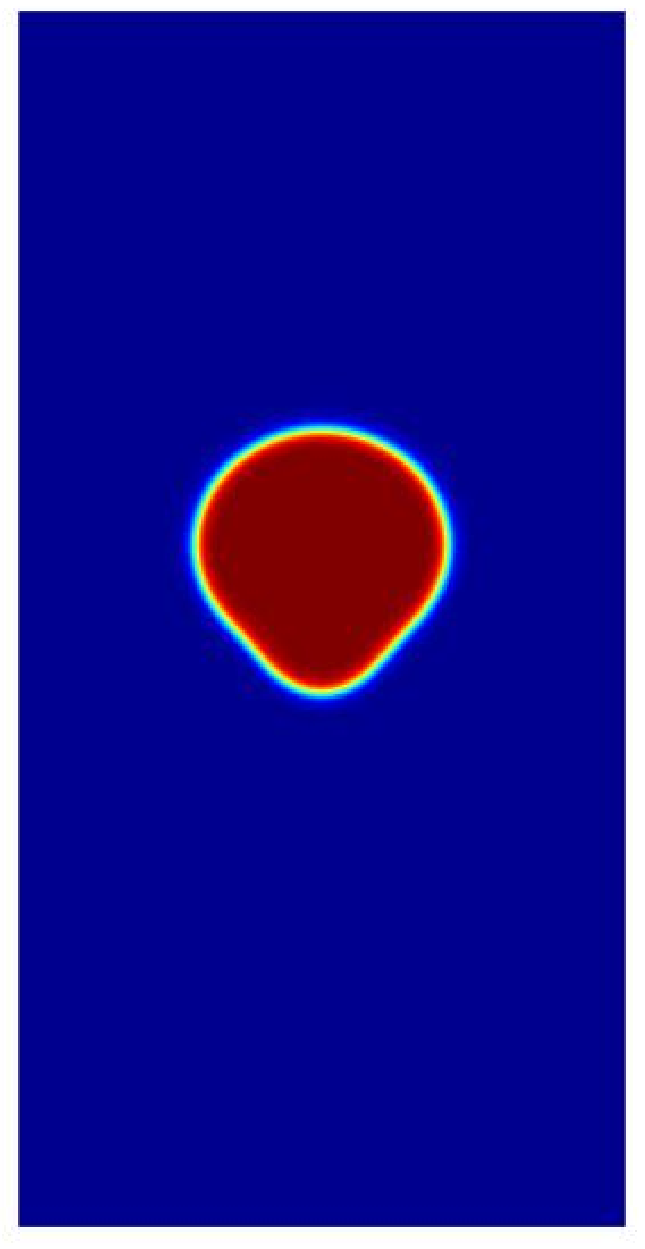}}
\hskip -0.7in
\subfigure
{\label{figr2:subfig:kk} 
\includegraphics[width=2.00in]{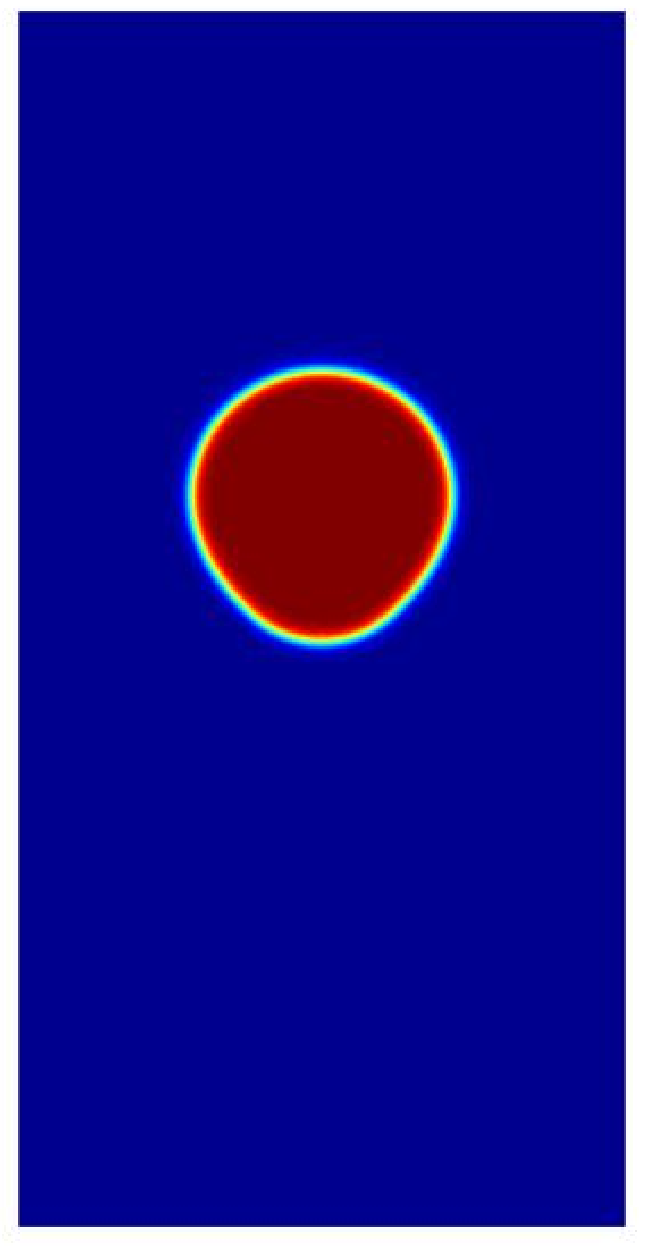}}
\hskip -0.7in
\subfigure
{\label{figr2:subfig:k} 
\includegraphics[width=2.00in]{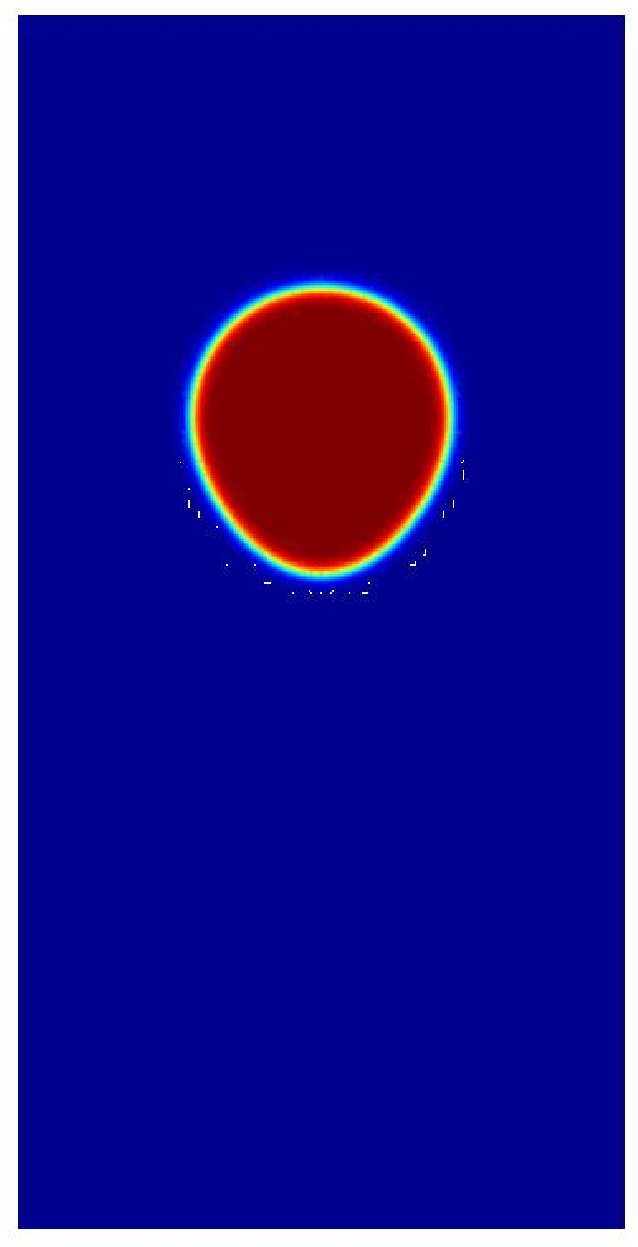}}
\vskip -0.25in
\caption{The evolution of arising drop with density ratio 1:5. All the sub-figures are indexed from left to right row by row as follows: : (a)~$t=1.0$, (b)~$t=2.0$, (c)~$t=3.0$, (d)~$t=4.0$, (e)~$t=5.0$, (f)~$t=6.0$, (g)~$t=8.0$, (h)~$t=11.0$.} \label{fig:graphr2}
\end{figure}

\begin{figure}[!htp]
\centering
\subfigure{
\label{Lfigr1:subfig:d} 
\includegraphics[width=2in]{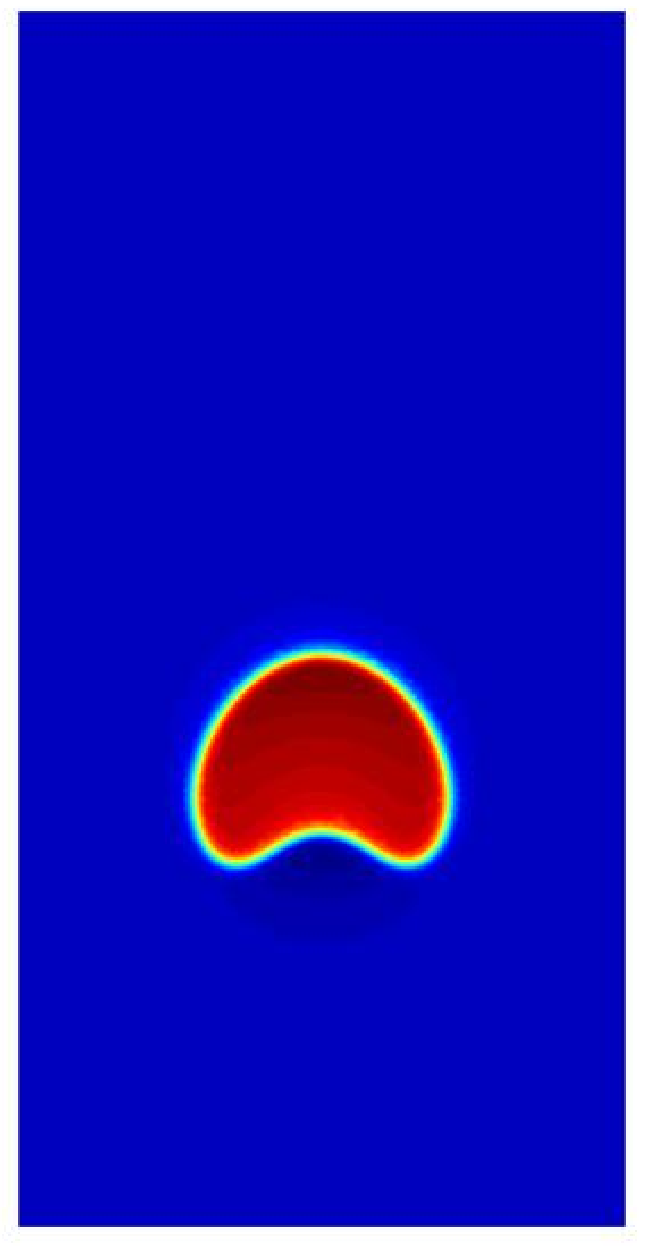}}
\hskip -0.70in
\subfigure{
\label{Lfigr1:subfig:f} 
\includegraphics[width=2in]{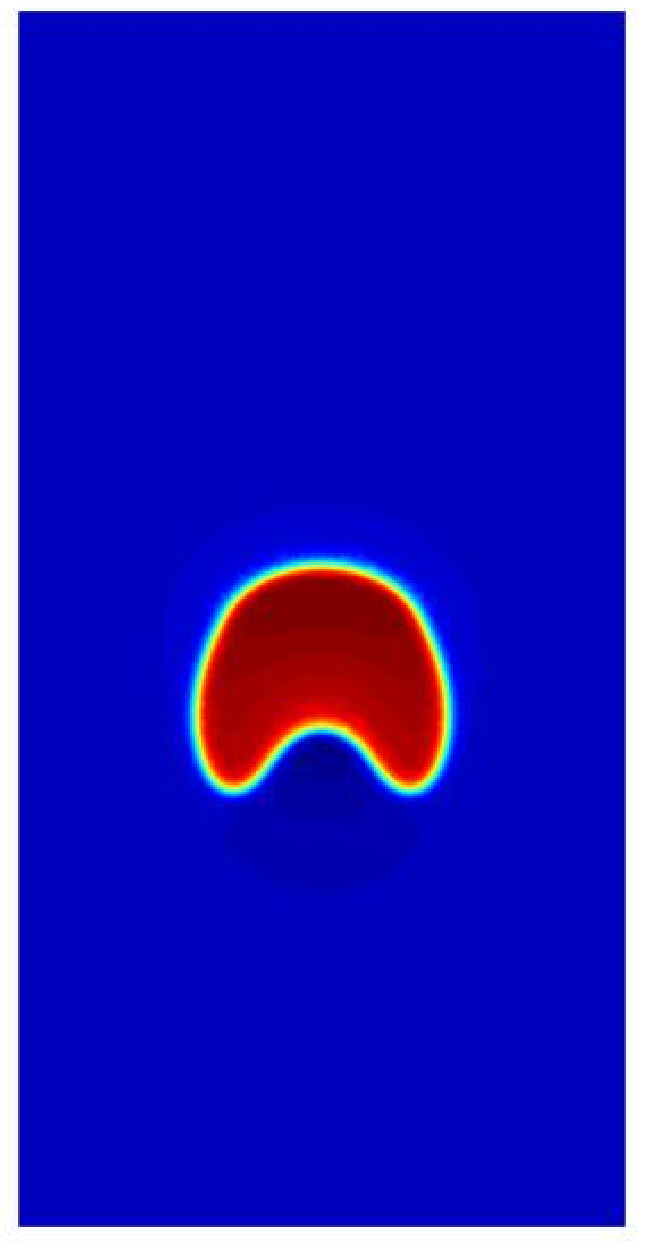}}
\hskip -0.7in
\subfigure
{\label{Lfigr2:subfig:k} 
\includegraphics[width=2in]{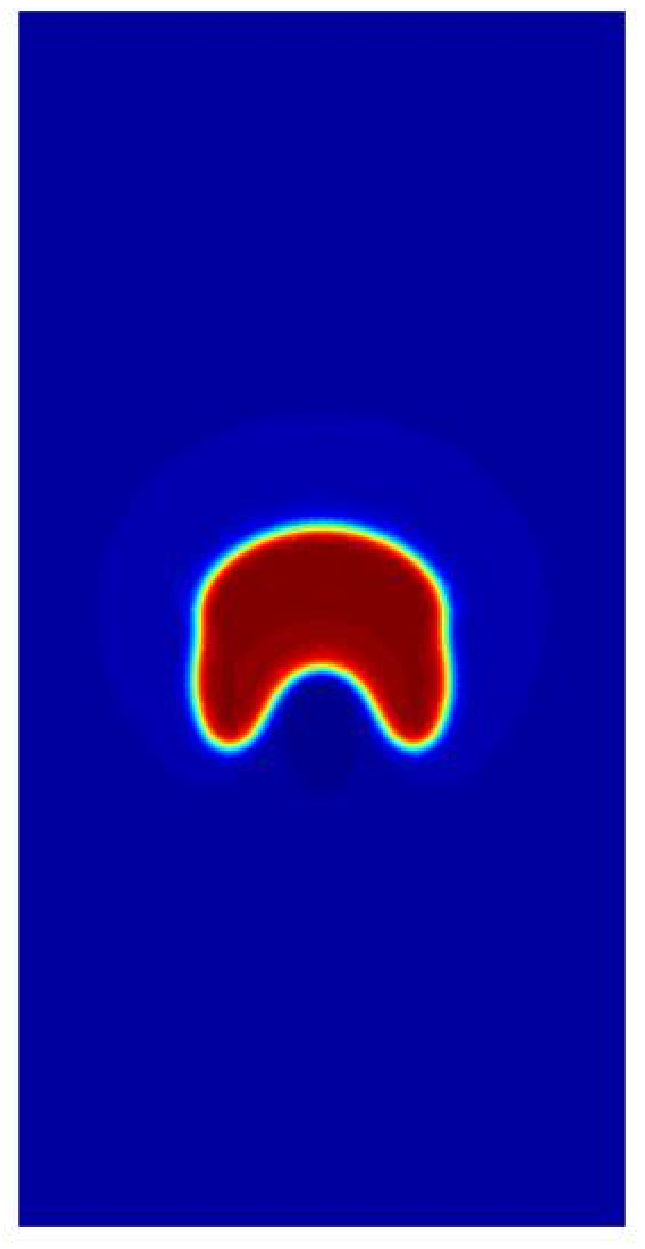}}
\hskip -0.7in
\subfigure
{\label{Lfigr1:subfig:g} 
\includegraphics[width=2in]{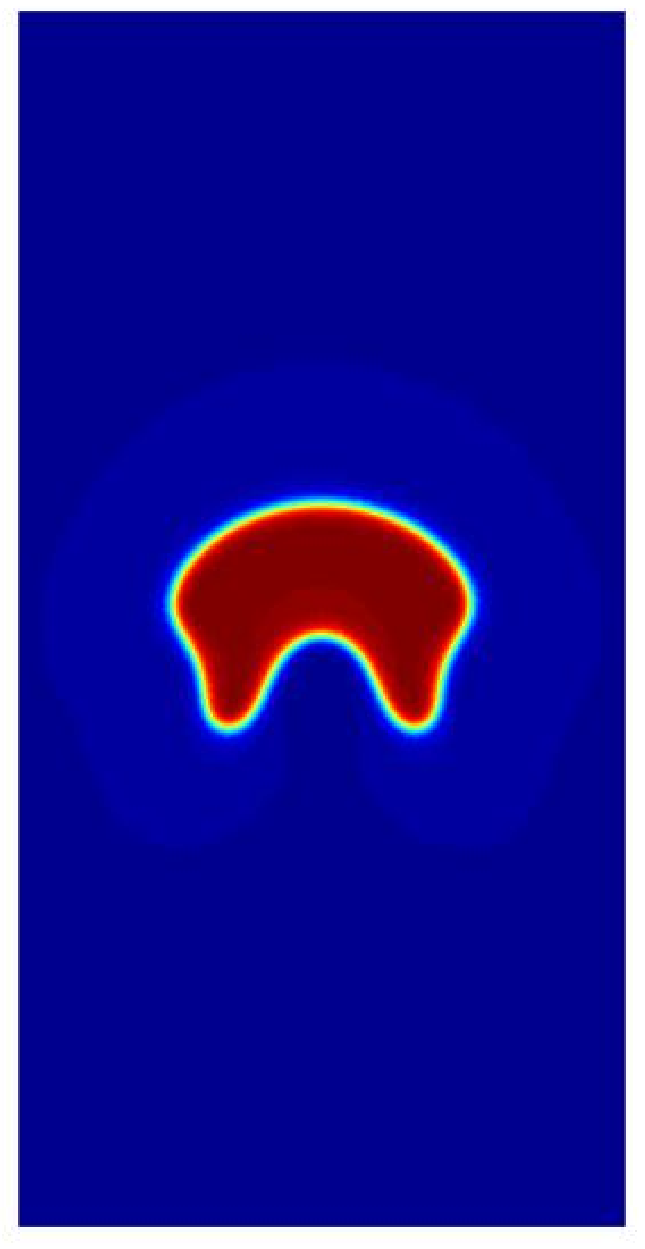}}
\vskip -0.40in
{\label{Lfigr1:subfig:ll} 
\includegraphics[width=2in]{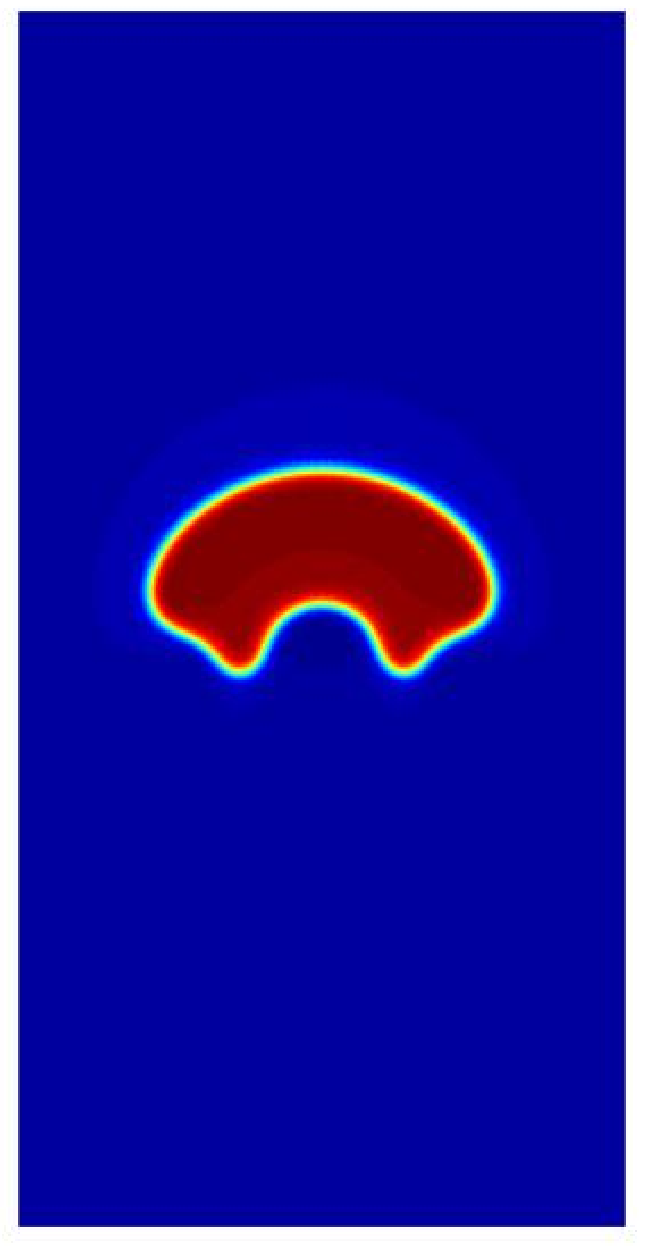}}
\hskip -0.7in
\label{Lfig:graphr1}
\subfigure
{\label{Lfigr2:subfig:b} 
\includegraphics[width=2in]{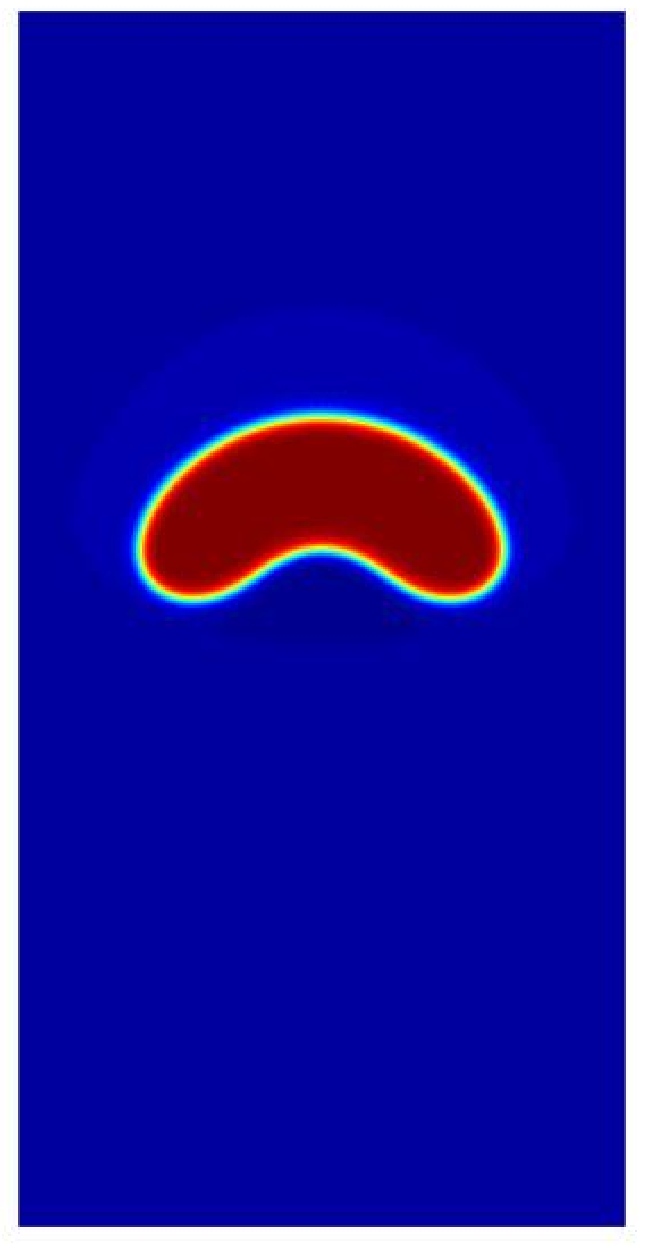}}
\hskip -0.7in
\subfigure
{\label{Lfigr1:subfig:h} 
\includegraphics[width=2in]{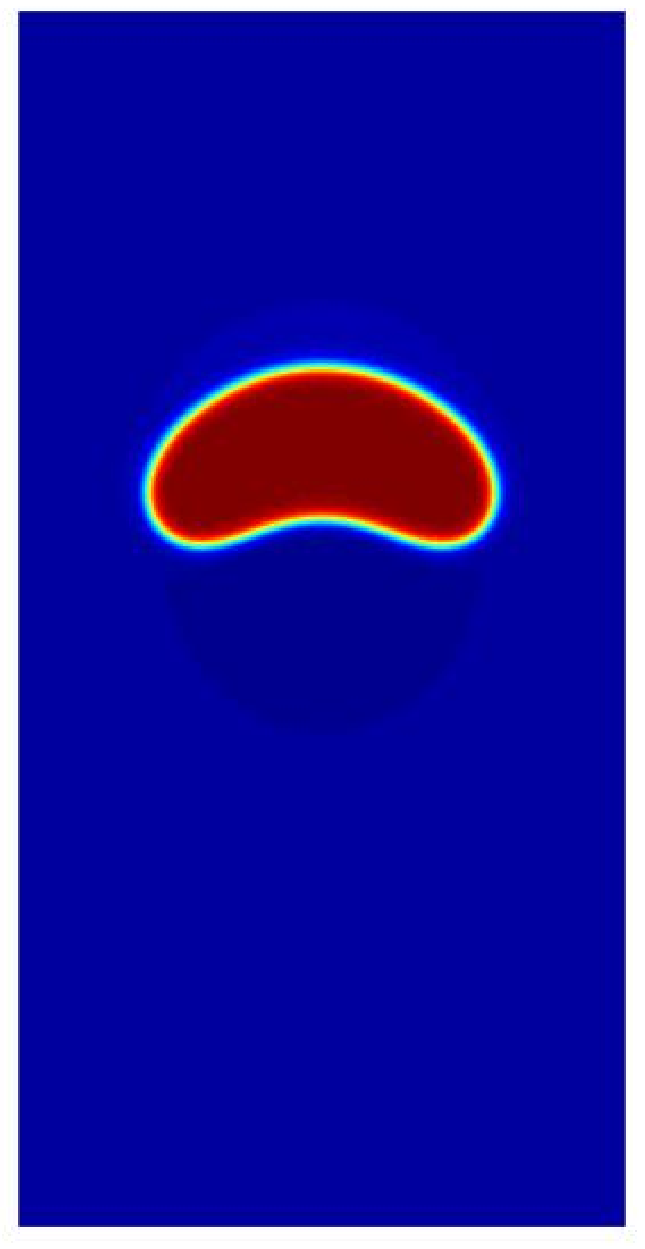}}
\hskip -0.70in
\subfigure
{
\label{Lfigr1:subfig:b} 
\includegraphics[width=2in]{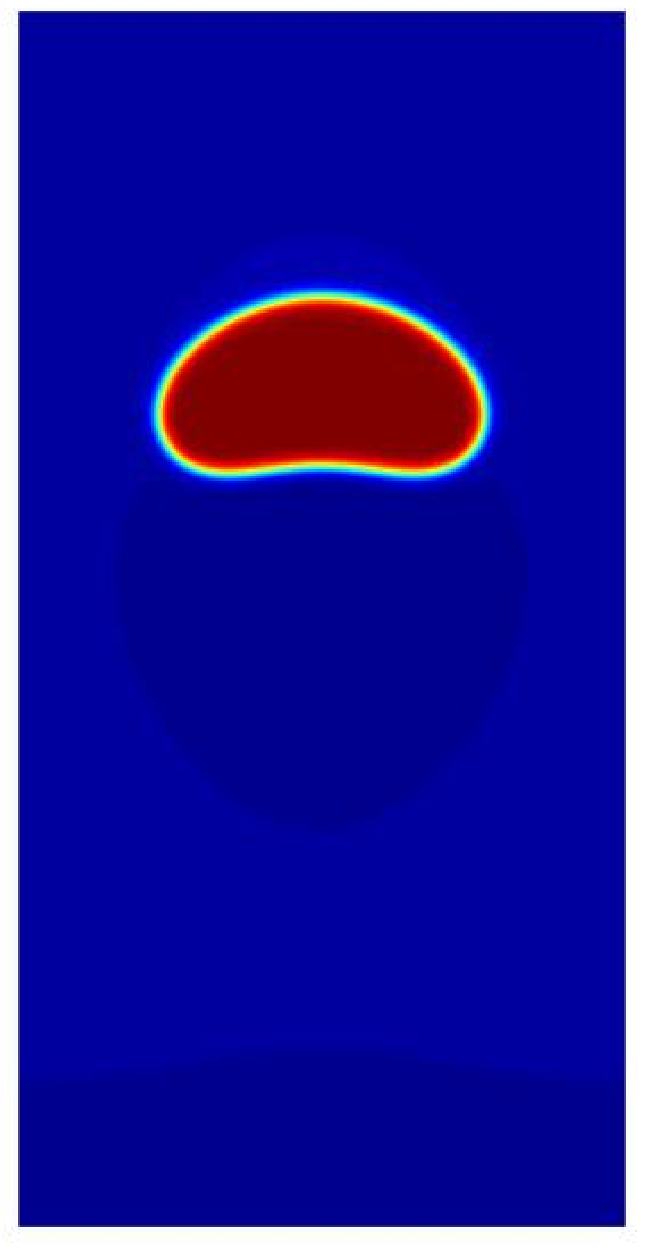}}
\vskip -0.25in
\caption{The evolution of arising drop with density ratio 1:50. All the sub-figures are indexed from left to right row by row as follows: : (a)~$t=1.0$, (b)~$t=1.25$, (c)~$t=1.5$, (d)~$t=2.0$, (e)~$t=2.5$, (f)~$t=3.0$, (g)~$t=4.0$, (h)~$t=5.0$.} \label{Lfig:graphr2}
\end{figure}

We test two cases with density ratios~$1:5$ and $1:50$, respectively.
Figure \ref{fig:graphr2} shows
several snapshots of the droplet passing through the interface  under the influence of
 buoyancy with a density ratio of $\rho_1:\rho_2=1:5$.
As the bubble rises in the conduit domain, it deforms into an ellipsoid. When it passes through the domain interface, one can clearly see an interface separating the bubble in conduit and in the porous medium.
The shape evolution of the rising bubble is shown in Figure \ref{Lfig:graphr2}
for the density ratio~$\rho_1:\rho_2=1:50$.
We can observe that the droplet quickly deforms into a heart-like shape as
compared with those in Figure \ref{fig:graphr2}.
As the droplet moves through the interface, the interface separates the bubble in conduit and matrix as presented in Figures \ref{figr1:subfig:f} and \ref{Lfigr2:subfig:k}. The smooth and excepted shape change of the droplet further validates physically faithful interface conditions when the droplet across the interface.  The tail is seen as it leaves the interface in Figures  \ref{figr1:subfig:ll} and \ref{Lfigr1:subfig:g}.
The tail is eventually smoothed out by the surface tension effect when it completely enters the porous medium as shown in Figures \ref{figr2:subfig:k} and \ref{Lfigr1:subfig:h}.

\begin{figure}[htp]
\centering
\subfigure
{\label{MLfigr1:subfig:f} 
\includegraphics[width=2in]{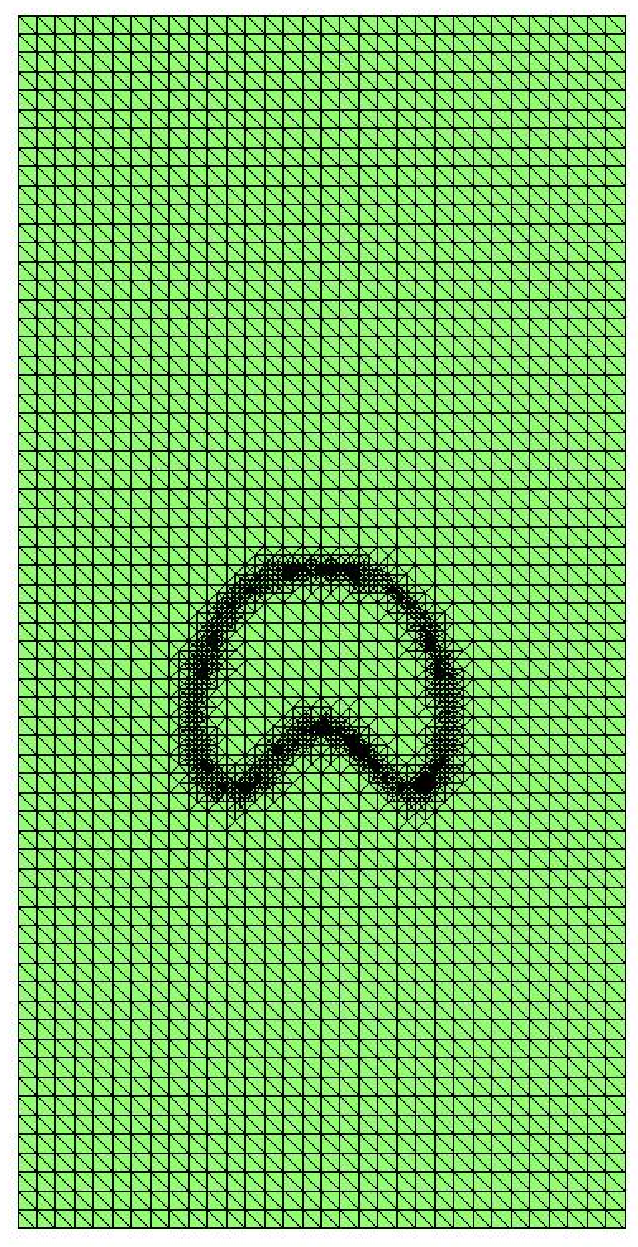}}
\hskip -0.7in
\subfigure
{\label{MLfigr2:subfig:k} 
\includegraphics[width=2in]{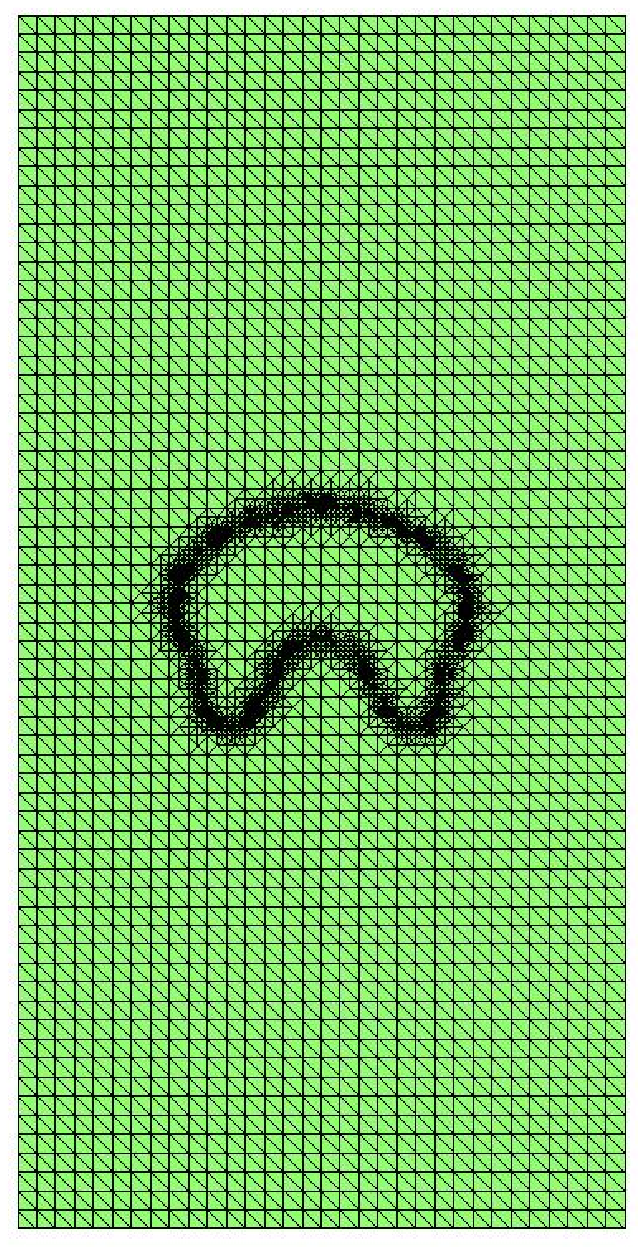}}
\hskip -0.70in
\subfigure
{\label{MLfigr1:subfig:g} 
\includegraphics[width=2in]{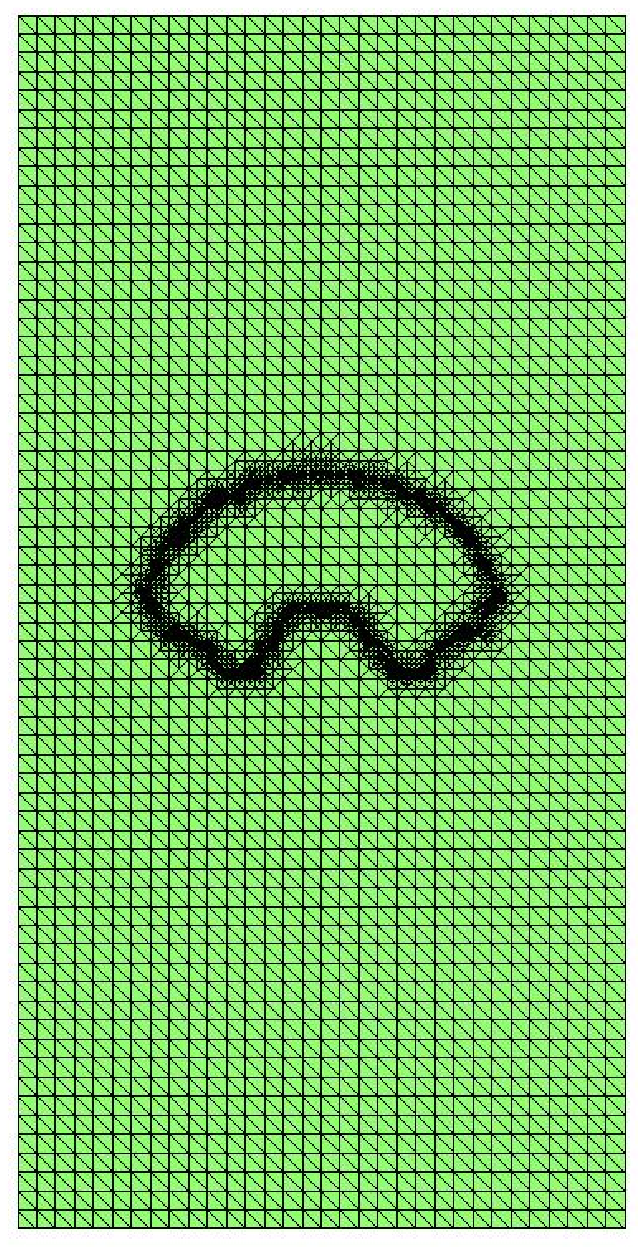}}
\hskip -0.70in
\subfigure
{\label{MLfigr1:subfig:h} 
\includegraphics[width=2in]{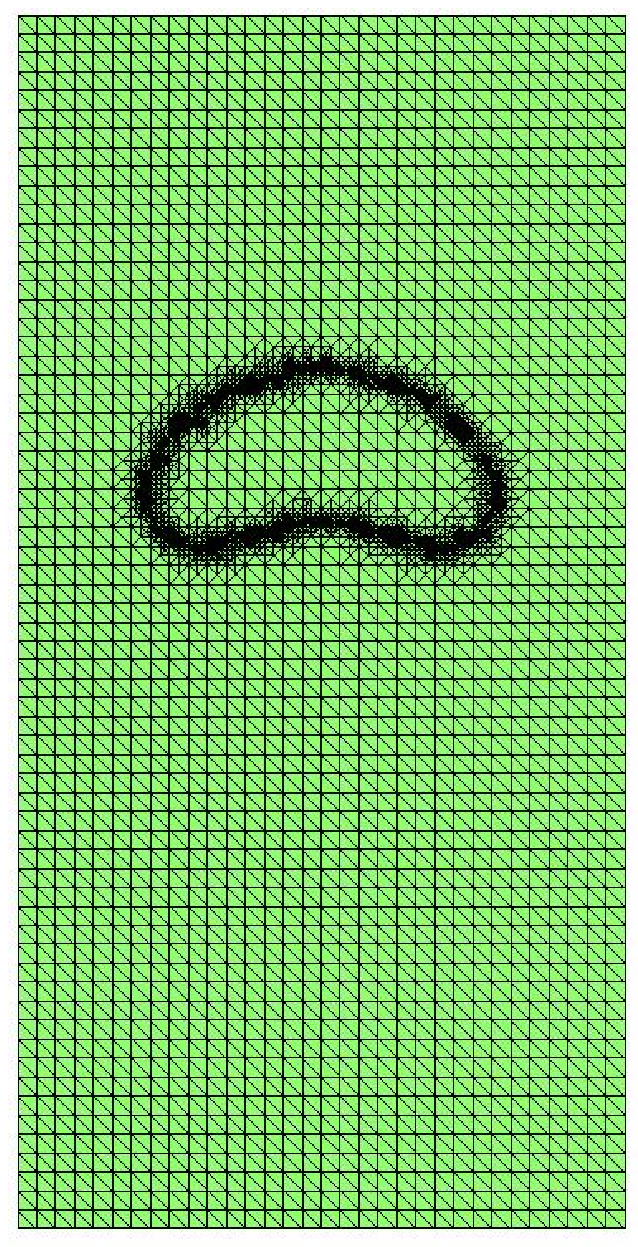}}
\vskip -0.25in
\caption{Adaptive mesh for  arising drop with density ratio 1:50. All the sub-figures are indexed from left to right row by row as follows: : (a)~$t=1.25$, (b)~$t=2.0$, (c)~$t=2.5$, (d)~$t=4.0$.} \label{MLfig:graphr2}
\end{figure}

 Additionally, we plot typical mesh refinement in Figures \ref{figIniA:subfig:b} and \ref{MLfig:graphr2} for this example. Once again, we observe that the mesh is properly refined near the interfacial region. All of these reasonable observations validate the interface conditions, the mathematical model and the numerical method proposed in  this article.

\section{Conclusions}
In this paper,  a new Cahn-Hilliard-Navier-Stokes-Darcy (CHNSD) model and its decoupled numerical schemes  are developed for two-phase
flows of different densities and viscosities in superposed fluid and porous layers.
Moreover,  the unconditionally energy stability  is proposed and analyzed for a  time-stepping method combining with interface conditions.
The novel decoupled numerical scheme is designed  by introducing the artificial compressibility method and pressure stabilization strategy.
The energy law  is proposed and analyzed for the corresponding fully discretization  in the framework of  the finite element method for spatial discretizaion.
Therefore,  only a sequence of linear equations is needed to solve at each discrete time level for the computation of the new decoupled linear numerical method.
The features of the proposed methods, such as
accuracy, energy dissipation, and applicability  for challenging model scenarios,
are demonstrated by the various numerical experiments.

\section*{Acknowledgement}
Gao is partially supported by the NSFC grant 11901461, the Natural Science Foundation of Shaanxi province 2019JQ-024 and and China Postdoctoral Science Foundation 2020M673464. Han acknowledges support from the NSF grant DMS-1912715. He is partially supported by the Alexander von Humboldt Foundation, the NSF grants  DMS-1722647 and DMS-1818642.

\bibliographystyle{plain}
\bibliography{gaoyli_ref_papers,gaoyli_ref_books,xiaoming_ref_papers,xiaoming_ref_books,ruede_articles,multiphase-2018}

\def\cprime{$'$}
\begin{thebibliography}{10}

\bibitem{HAbels_HGarcke_GGrun_2012}
H.~Abels, H.~Garcke, and G.~Gr\"{u}n.
\newblock Thermodynamically consistent, frame indifferent diffuse interface
  models for incompressible two-phase flows with different densities.
\newblock {\em Math. Models Methods Appl. Sci.}, 22(3):1150013, 2012.

\bibitem{TArbogast_MGomez_1}
T.~Arbogast and M.~Gomez.
\newblock {A discretization and multigrid solver for a Darcy-Stokes system of
  three dimensional vuggy porous media}.
\newblock {\em Comput. Geosci.}, 13(3):331--348, 2009.

\bibitem{BLWW2013}
A.~Baskaran, J.~S. Lowengrub, C.~Wang, and S.~M. Wise.
\newblock Convergence analysis of a second order convex splitting scheme for
  the modified phase field crystal equation.
\newblock {\em SIAM J. Numer. Anal.}, 51(5):2851--2873, 2013.

\bibitem{GBeavers_DJoseph_1}
G.~Beavers and D.~Joseph.
\newblock Boundary conditions at a naturally permeable wall.
\newblock {\em J. Fluid Mech.}, 30:197--207, 1967.

\bibitem{YBoubendir_STlupova_2}
Y.~Boubendir and S.~Tlupova.
\newblock {Domain decomposition methods for solving Stokes-Darcy problems with
  boundary integrals}.
\newblock {\em SIAM J. Sci. Comput.}, 35(1):B82--B106, 2013.

\bibitem{FBoyer_3}
F.~Boyer.
\newblock {A theoretical and numerical model for the study of incompressible
  mixture flows}.
\newblock {\em Comput. Fluids}, 31(1):41--68, 2002.

\bibitem{MCai_MMu_JXu_2}
M.~Cai, M.~Mu, and J.~Xu.
\newblock {Numerical solution to a mixed Navier-Stokes/Darcy model by the
  two-grid approach}.
\newblock {\em SIAM J. Numer. Anal.}, 47(5):3325--3338, 2009.

\bibitem{JCamano_GNGatica_ROyarzua_RRuizBaier_PVenegas_1}
J.~Camano, G.~N. Gatica, R.~Oyarzua, R.~Ruiz-Baier, and P.~Venegas.
\newblock {New fully-mixed finite element methods for the Stokes-Darcy
  coupling}.
\newblock {\em Comput. Methods Appl. Mech. Engrg.}, 295:362--395, 2015.

\bibitem{YCao_MGunzburger_XMHe_XWang_1}
Y.~Cao, M.~Gunzburger, X.-M. He, and X.~Wang.
\newblock {Robin-Robin} domain decomposition methods for the steady
  {Stokes-Darcy model with Beaver-Joseph} interface condition.
\newblock {\em Numer. Math.}, 117(4):601--629, 2011.

\bibitem{YCao_MGunzburger_XMHe_XWang_2}
Y.~Cao, M.~Gunzburger, X.-M. He, and X.~Wang.
\newblock {Parallel, non-iterative, multi-physics domain decomposition methods
  for time-dependent Stokes-Darcy systems}.
\newblock {\em Math. Comp.}, 83(288):1617--1644, 2014.

\bibitem{YCao_MGunzburger_XHu_FHua_XWang_WZhao_1}
Y.~Cao, M.~Gunzburger, X.~Hu, F.~Hua, X.~Wang, and W.~Zhao.
\newblock {Finite element approximation for Stokes-Darcy flow with
  Beavers-Joseph interface conditions}.
\newblock {\em SIAM. J. Numer. Anal.}, 47(6):4239--4256, 2010.

\bibitem{MBCardenas_WRR_2015}
M.~Bayani Cardenas.
\newblock Hyporheic zone hydrologic science: {A} historical account of its
  emergence and a prospectus.
\newblock {\em Water Resour. Res.}, 51:3601--3616, 2015.

\bibitem{JChen_SSun_XWang_1}
J.~Chen, S.~Sun, and X.~Wang.
\newblock {A numerical method for a model of two-phase flow in a coupled free
  flow and porous media system}.
\newblock {\em J. Comput. Phys.}, 268:1--16, 2014.

\bibitem{WChen_DHan_XWang_2017}
W.~Chen, D.~Han, and X.~Wang.
\newblock Uniquely solvable and energy stable decoupled numerical schemes for
  the {C}ahn-{H}illiard-{S}tokes-{D}arcy system for two-phase flows in karstic
  geometry.
\newblock {\em Numer. Math.}, 137(1):229--255, 2017.

\bibitem{QChen_XYang_JShen_JCP_2017}
Q.~Cheng, X.~Yang, and J.~Shen.
\newblock Efficient and accurate numerical schemes for a hydro-dynamically
  coupled phase field diblock copolymer model.
\newblock {\em J. Comput. Phys.}, 341:44--60, 2017.

\bibitem{PChidyagwai_BRiviere_1}
P.~Chidyagwai and B.~Rivi\`{e}re.
\newblock {On the solution of the coupled Navier-Stokes and Darcy equations}.
\newblock {\em Comput. Methods Appl. Mech. Engrg.}, 198(47-48):3806--3820,
  2009.

\bibitem{AJChorin_2}
A.~J. Chorin.
\newblock A numerical method for solving incompressible viscous flow problems.
\newblock {\em J. Comput. Phys}, 2:12, 1967.

\bibitem{NCondette_CMMelcher_ESuli_2011}
N.~Condette, C.~Melcher, and E.~S\"{u}li.
\newblock Spectral approximation of pattern-forming nonlinear evolution
  equations with double-well potentials of quadratic growth.
\newblock {\em Math. Comp.}, 80:205--223, 2011.

\bibitem{VDeCaria_WLayton_MMcLaughlin_CMAME_2017}
V.~DeCaria, W.~Layton, and M.~McLaughlin.
\newblock A conservative, second order, unconditionally stable artificial
  compression method.
\newblock {\em Comput. Methods Appl. Mech. Engrg.}, 325:733--747, 2017.

\bibitem{HDing_PDMSpelt_CShu_1}
H.~Ding, P.~D.~M. Spelt, and C.~Shu.
\newblock {Diffuse interface model for incompressible two-phase flows with
  large density ratios}.
\newblock {\em J. Comput. Phys.}, 226(2):2078--2095, 2007.

\bibitem{MDiscacciati_LGerardo-Giorda_1}
M.~Discacciati and L.~Gerardo-Giorda.
\newblock {Optimized Schwarz methods for the Stokes-Darcy coupling}.
\newblock {\em IMA J. Numer. Anal.}, 38(4):1959-1983, 2018.

\bibitem{MDiscacciati_EMiglio_AQuarteroni_1}
M.~Discacciati, E.~Miglio, and A.~Quarteroni.
\newblock Mathematical and numerical models for coupling surface and
  groundwater flows.
\newblock {\em Appl. Numer. Math.}, 43(1-2):57--74, 2002.

\bibitem{MDiscacciati_AQuarteroni_AValli_1}
M.~Discacciati, A.~Quarteroni, and A.~Valli.
\newblock {Robin-Robin domain decomposition methods for the Stokes-Darcy
  coupling}.
\newblock {\em SIAM J. Numer. Anal.}, 45(3):1246--1268, 2007.

\bibitem{ElSt1993}
C.~M. Elliott and A.~M. Stuart.
\newblock The global dynamics of discrete semilinear parabolic equations.
\newblock {\em SIAM J. Numer. Anal.}, 30(6):1622--1663, 1993.

\bibitem{JAFiordilino_WLayton_YRong_CMAME_2018}
J.~A. Fiordilino, W.~Layton, and Y.~Rong.
\newblock An efficient and modular {grad-div} stabilization.
\newblock {\em Comput. Methods Appl. Mech. Engrg.}, 335:327--346, 2018.

\bibitem{LFormaggia_AQuarteroni_AVeneziani_2009}
L.~Formaggia, A.~Quarteroni, and A.~Veneziani.
\newblock {\em Cardiovascular mathematics: modeling and simulation of the
  circulatory system}.
\newblock Springer-Verlag, New York, 2009.

\bibitem{MGao_XWang_1}
M.~Gao and X.~Wang.
\newblock {A gradient stable scheme for a phase field model for the moving
  contact line problem}.
\newblock {\em J. Comput. Phys.}, 231(4):1372--1386, 2012.

\bibitem{YGao_XHe_LMei_XYang_2018}
Y.~Gao, X.~He, L.~Mei, and X.~Yang.
\newblock Decoupled, linear, and energy stable finite element method for the
  {Cahn-Hilliard-Navier-Stokes-Darcy} phase field model.
\newblock {\em SIAM J. Sci. Comput.}, 40(1):B110--B137, 2018.

\bibitem{VGirault_BRiviere_1}
V.~Girault and B.~Rivi\`{e}re.
\newblock {DG approximation of coupled Navier-Stokes and Darcy equations by
  Beaver-Joseph-Saffman interface condition}.
\newblock {\em SIAM J. Numer. Anal}, 47(3):2052--2089, 2009.

\bibitem{GlSw2004}
J.~G. Gluyas and R.~E. Swarbrick.
\newblock {\em Petroleum Geology}.
\newblock Blackwell publishing, 2004.

\bibitem{Grun2013}
G.~Gr{\"u}n.
\newblock On convergent schemes for diffuse interface models for two-phase flow
  of incompressible fluids with general mass densities.
\newblock {\em SIAM J. Numer. Anal.}, 51(6):3036--3061, 2013.

\bibitem{JLGuermond_LQuartapelle_2000}
J.-L. Guermond and L.~Quartapelle.
\newblock A projection {FEM} for variable density incompressible flows.
\newblock {\em J. Comput. Phys.}, 165(1):167--188, 2000.

\bibitem{MGunzburger_XMHe_BLi_SINN_2018}
M.~Gunzburger, X.-M. He, and B.~Li.
\newblock {On Ritz projection and multi-step backward differentiation schemes
  in decoupling the Stokes-Darcy model}.
\newblock {\em SIAM J. Numer. Anal.}, 56(1):397--427, 2018.

\bibitem{GLLW2017}
Z.~Guo, P.~Lin, J.~Lowengrub, and S.~M. Wise.
\newblock Mass conservative and energy stable finite difference methods for the
  quasi-incompressible {N}avier-{S}tokes-{C}ahn-{H}illiard system: primitive
  variable and projection-type schemes.
\newblock {\em Comput. Methods Appl. Mech. Engrg.}, 326:144--174, 2017.

\bibitem{MLHadji_AAssala_FZNouri_1}
M.~L. Hadji, A.~Assala, and F.~Z. Nouri.
\newblock {A posteriori error analysis for Navier-Stokes equations coupled with
  Darcy problem}.
\newblock {\em Calcolo}, 52(4):559--576, 2015.

\bibitem{DHan_JSC_2016}
D.~Han.
\newblock A decoupled unconditionally stable numerical scheme for the
  {Cahn-Hilliard-Hele-Shaw} system.
\newblock {\em J. Sci Comput.}, 66(3):1102--1121, 2016.

\bibitem{DHan_DSun_XWang_1}
D.~Han, D.~Sun, and X.~Wang.
\newblock {Two-phase flows in karstic geometry}.
\newblock {\em Math. Methods Appl. Sci.}, 37(18):3048--3063, 2014.

\bibitem{DHan_XWang_HWu_1}
D.~Han, X.~Wang, and H.~Wu.
\newblock {Existence and uniqueness of global weak solutions to a
  Cahn-Hilliard-Stokes-Darcy system for two phase incompressible flows in
  karstic geometry}.
\newblock {\em J. Differential Equations}, 257(10):3887--3933, 2014.

\bibitem{XHe_NJiang_CQiu_IJNME_2019}
X.-M. He, N.~Jiang, and C.~Qiu.
\newblock An artificial compressibility ensemble algorithm for a stochastic
  {Stokes-Darcy} model with random hydraulic conductivity and interface
  conditions.
\newblock {\em Int. J. Numer. Methods Eng.}, pages 1--28, 2019.

\bibitem{XMHe_JLi_YLin_JMing_1}
X.-M. He, J.~Li, Y.~Lin, and J.~Ming.
\newblock {A domain decomposition method for the steady-state
  Navier-Stokes-Darcy model with Beavers-Joseph interface condition}.
\newblock {\em SIAM J. Sci. Comput.}, 37(5):S264--S290, 2015.

\bibitem{JHou_MQiu_XMHe_CGuo_MWei_BBai_1}
J.~Hou, M.~Qiu, X.-M. He, C.~Guo, M.~Wei, and B.~Bai.
\newblock {A dual-porosity-Stokes model and finite element method for coupling
  dual-porosity flow and free flow}.
\newblock {\em SIAM J. Sci. Comput.}, 38(5):B710--B739, 2016.

\bibitem{DKay_RWelford_1}
D.~Kay and R.~Welford.
\newblock {Efficient numerical solution of Cahn-Hilliard-Navier-Stokes fluids
  in 2D}.
\newblock {\em SIAM J. Sci. Comput.}, 29(6):2241--2257, 2007.

\bibitem{JKim_KKang_JLowengrub_JCP_2004}
J.~Kim, K.~Kang, and J.~Lowengrub.
\newblock Conservative multigrid methods for {Cahn-Hilliard} fluids.
\newblock {\em J. Comput. Phys.}, 193:511--543, 2004.

\bibitem{AGLamorgese_DMolin_RMauri_MJM_2011}
A.~G. Lamorgese, D.~Molin, and R.~Mauri.
\newblock Phase field approach to multiphase flow modeling.
\newblock {\em Milan J. Math.}, 79(2):597--642, 2011.

\bibitem{WJLayton_FSchieweck_IYotov_1}
W.~J. Layton, F.~Schieweck, and I.~Yotov.
\newblock Coupling fluid flow with porous media flow.
\newblock {\em SIAM J. Numer. Anal.}, 40(6):2195--2218, 2002.

\bibitem{HGLee_JLowengrub_JGoodman_PFI_2002}
H.~G. Lee, J.~Lowengrub, and J.~Goodman.
\newblock {Modeling pinchoff and reconnection in a Hele-Shaw cell. I. The
  models and their calibration}.
\newblock {\em Phys. Fluids}, 14(2):492--513, 2002.

\bibitem{HGLee_JLowengrub_JGoodman_PFII_2002}
H.~G. Lee, J.~Lowengrub, and J.~Goodman.
\newblock {Modeling pinchoff and reconnection in a Hele-Shaw cell. II. Analysis
  and simulation in the nonlinear regime}.
\newblock {\em Phys. Fluids}, 14(2):514--545, 2002.

\bibitem{RLi_YGao_JChen_LZhang_XMHe_ZChen_1}
R.~Li, Y.~Gao, J.~Chen, L.~Zhang, X.-M. He, and Z.~Chen.
\newblock {Discontinuous finite volume element method for a coupled
  Navier-Stokes-Cahn-Hilliard phase field model}.
\newblock {\em Adv. Comput. Math.}, 46:\#25, 2020.

\bibitem{FLin_XMHe_XWen_1}
F.~Lin, X.-M. He, and X.~Wen.
\newblock {Fast, unconditionally energy stable large time stepping method for a
  new Allen-Cahn type square phase-field crystal model}.
\newblock {\em Appl. Math. Lett.}, 92:248--255, 2019.

\bibitem{KLipnikov_DVassilev_IYotov_1}
K.~Lipnikov, D.~Vassilev, and I.~Yotov.
\newblock {Discontinuous Galerkin and mimetic finite difference methods for
  coupled Stokes-Darcy flows on polygonal and polyhedral grids}.
\newblock {\em Numer. Math.}, 126(2):321--360, 2014.

\bibitem{YLiu_WChen_CWang_SMWise_NM_2017}
Y.~Liu, W.~Chen, C.~Wang, and S.M. Wise.
\newblock Error analysis of a mixed finite element method for a
  {Cahn-Hilliard-Hele-Shaw} system.
\newblock {\em Numer. Math.}, 135(3):679--709, 2017.

\bibitem{JLowengrub_LTruskinovsky_1998}
J.~Lowengrub and L.~Truskinovsky.
\newblock Quasi-incompressible {Cahn-Hilliard} fluids and topological
  transitions.
\newblock {\em R. Soc. Lond. Proc. Ser. A Math. Phys. Eng. Sci.},
  454(1978):2617--2654, 1998.

\bibitem{MdAAlMahbub_XMHe_NJNasu_CQiu_HZheng_1}
Md. A.~Al Mahbub, X.-M. He, N.~J. Nasu, C.~Qiu, and H.~Zheng.
\newblock {Coupled and decoupled stabilized mixed finite element methods for
  non-stationary dual-porosity-Stokes fluid flow model}.
\newblock {\em Int. J. Numer. Meth. Eng.}, 120(6):803--833, 2019.

\bibitem{AMarquez_SMeddahi_FJSayas_2}
A.~M\'{a}rquez, S.~Meddahi, and F.~J. Sayas.
\newblock {Strong coupling of finite element methods for the Stokes-Darcy
  problem}.
\newblock {\em IMA J. Numer. Anal.}, 35(2):969--988, 2015.

\bibitem{JMatusick_PZanbergen_GRA_2007}
J.~Matusick and P.~Zanbergen.
\newblock Comparative study of groundwater vulnerability in a karst aquifer in
  central florida.
\newblock {\em Geophy. Res. Abst.}, 9:1--1, 2007.

\bibitem{MMoraiti_1}
M.~Moraiti.
\newblock {On the quasistatic approximation in the Stokes-Darcy model of
  groundwater-surface water flows}.
\newblock {\em J. Math. Anal. Appl.}, 394(2):796--808, 2012.

\bibitem{MMu_JXu_1}
M.~Mu and J.~Xu.
\newblock A two-grid method of a mixed {Stokes-Darcy} model for coupling fluid
  flow with porous media flow.
\newblock {\em SIAM J. Numer. Anal.}, 45(5):1801--1813, 2007.

\bibitem{CQiu_XMHe_JLi_YLin_1}
C.~Qiu, X.-M. He, J.~Li, and Y.~Lin.
\newblock {A domain decomposition method for the time-dependent
  Navier-Stokes-Darcy model with Beavers-Joseph interface condition and
  defective boundary condition}.
\newblock {\em J. Comput. Phys.}, 411:\#109400, 2020.

\bibitem{BRiviere_1}
B.~Rivi\`{e}re.
\newblock Analysis of a discontinuous finite element method for the coupled
  {Stokes and Darcy} problems.
\newblock {\em J. Sci. Comput.}, 22/23:479--500, 2005.

\bibitem{BRiviere_IYotov_1}
B.~Rivi\`{e}re and I.~Yotov.
\newblock Locally conservative coupling of {Stokes and Darcy} flows.
\newblock {\em SIAM J. Numer. Anal.}, 42(5):1959--1977, 2005.

\bibitem{JShen_CWang_XWang_SMWise_SINA_2012}
J.~Shen, C.~Wang, X.~Wang, and S.M. Wise.
\newblock Second-order convex splitting schemes for gradient flows with
  {E}hrlich-{S}chwoebel type energy: application to thin film epitaxy.
\newblock {\em SIAM J. Numer. Anal.}, 50(1):105--125, 2012.

\bibitem{JShen_JXu_JYang_JCP_2018}
J.~Shen, J.~Xu, and J.~Yang.
\newblock The scalar auxiliary variable ({SAV}) approach for gradient flows.
\newblock {\em J. Comput. Phys.}, 353:407--416, 2018.

\bibitem{JShen_JXu_JYang_SIRV_2019}
J.~Shen, J.~Xu, and J.~Yang.
\newblock A new class of efficient and robust energy stable schemes for
  gradient flows.
\newblock {\em SIAM Rev.}, 61(3):474--506, 2019.

\bibitem{JShen_XYang_1}
J.~Shen and X.~Yang.
\newblock {A phase-field model and its numerical approximation for two-phase
  incompressible flows with different densities and viscosities}.
\newblock {\em SIAM J. Sci. Comput.}, 32(3):1159--1179, 2010.

\bibitem{ShenYang2010}
J.~Shen and X.~Yang.
\newblock Numerical approximations of {Allen-Cahn} and {Cahn-Hilliard}
  equations.
\newblock {\em Discrete Contin. Dyn. Syst.}, 28:1169--1691, 2010.

\bibitem{JShen_XYang_SINN_2015}
J.~Shen and X.~Yang.
\newblock {Decoupled, energy stable schemes for phase-field models of two-phase
  incompressible flows}.
\newblock {\em SIAM J. Numer. Anal.}, 53(1):279--296, 2015.

\bibitem{SKFStoter_PMuller_elat_CMAME_2017}
S.K.F. Stoter, P.~M\"uller, L.~Cicalese, M.~Tuveri, D.~Schillinger, and
  T.~J.~R. Hughes.
\newblock A diffuse interface method for the {N}avier-{S}tokes/{D}arcy
  equations: perfusion profile for a patient-specific human liver based on
  {MRI} scans.
\newblock {\em Comput. Methods Appl. Mech. Engrg.}, 321:70--102, 2017.

\bibitem{RTemam_BSMF_1968}
R.~Temam.
\newblock Une m\'ethode d'approximation de la solution des \'equations de
  {N}avier-{S}tokes.
\newblock {\em Bull. Soc. Math. France}, 96:115--152, 1968.

\bibitem{STlupova_RCortez_1}
S.~Tlupova and R.~Cortez.
\newblock Boundary integral solutions of coupled {Stokes and Darcy} flows.
\newblock {\em J. Comput. Phys.}, 228(1):158--179, 2009.

\bibitem{KTuber_DPocza_CHebling_JPS_2003}
K.~Tuber, D.~Pocza, and C.~Hebling.
\newblock Visualization of water buildup in the cathode of a transparent {PEM}
  fuel cell.
\newblock {\em J. Power Sources}, 124(2):403--414, 2003.

\bibitem{DVassilev_CWang_IYotov_1}
D.~Vassilev, C.~Wang, and I.~Yotov.
\newblock {Domain decomposition for coupled Stokes and Darcy flows}.
\newblock {\em Comput. Methods Appl. Mech. Engrg.}, 268:264--283, 2014.

\bibitem{CXu_CChen_XYang_XMHe_1}
C.~Xu, C.~Chen, X.~Yang, and X.-M. He.
\newblock {Numerical approximations for the hydrodynamics coupled binary
  surfactant phase field model: second order, linear, unconditionally energy
  stable schemes}.
\newblock {\em Commun. Math. Sci.}, 17(3):835--858, 2019.

\bibitem{Xu06}
C.~Xu and T.~Tang.
\newblock Stability analysis of large time-stepping methods for epitaxial
  growth models.
\newblock {\em SIAM J. Numer. Anal.}, 44(4):1759--1779, 2006.

\bibitem{YYan_WChen_CWang_SMWise_CICP_2018}
Y.~Yan, W.~Chen, C.~Wang, and S.M. Wise.
\newblock A second-order energy stable {BDF} numerical scheme for the
  {C}ahn-{H}illiard equation.
\newblock {\em Commun. Comput. Phys.}, 23(2):572--602, 2018.

\bibitem{NNYanenko_1971}
N.~N. Yanenko.
\newblock {\em The Method of Fractional Steps. The Solution of Problems of
  Mathematical Physics in Several Variables}.
\newblock Springer-Verlag, New York, 1971.

\bibitem{JYang_SMao_XMHe_XYang_YHe_1}
J.~Yang, S.~Mao, X.-M. He, X.~Yang, and Y.~He.
\newblock {A diffuse interface model and semi-implicit energy stable finite
  element method for two-phase magnetohydrodynamic flows}.
\newblock {\em Comput. Meth. Appl. Mech. Eng.}, 356:435--464, 2019.

\bibitem{XYang_DHan_JCP_2017}
X.~Yang and D.~Han.
\newblock Linearly first- and second-order, unconditionally energy stable
  schemes for the phase field crystal equation.
\newblock {\em J. Comput. Phys.}, 330:13--22, 2017.

\bibitem{XYang_LJu_CMAME_2017}
X.~Yang and L.~Ju.
\newblock Linear and unconditionally energy stable schemes for the binary
  fluid-surfactant phase field model.
\newblock {\em Comput. Methods Appl. Mech. Engrg.}, 318:1005--1029, 2017.

\bibitem{XYang_JZhao_XMHe_1}
X.~Yang, J.~Zhao, and X.-M. He.
\newblock {Linear, second order and unconditionally energy stable schemes for
  the viscous Cahn-Hilliard equation with hyperbolic relaxation using the
  invariant energy quadratization method}.
\newblock {\em J. Comput. Appl. Math.}, 343(1):80--97, 2018.

\bibitem{JZhao_XYang_YGong_QWang_CMAME_2017}
J.~Zhao, X.~Yang, Y.~Gong, and Q.~Wang.
\newblock A novel linear second order unconditionally energy stable scheme for
  a hydrodynamic {$\bold{Q}$}-tensor model of liquid crystals.
\newblock {\em Comput. Methods Appl. Mech. Engrg.}, 318:803--825, 2017.

\bibitem{ZhuCST1999}
J.~Zhu, L.~Q. Chen, J.~Shen, and V.~Tikare.
\newblock Morphological evolution during phase separation and coarsening with
  strong inhomogeneous elasticity.
\newblock {\em Model. Simul. Mater. Sci. Eng.}, 9(6):499--511, 2001.

\end{thebibliography}

\end{document}